\title{Braided Subfactors, Spectral Measures, \\ Planar algebras and Calabi-Yau algebras \\ associated to $SU(3)$ modular invariants}
\author{
        David E. Evans and Mathew Pugh \\ \\
        \small School of Mathematics, 
        \small Cardiff University, \\
        \small Senghennydd Road, 
        \small Cardiff, CF24 4AG, 
        \small Wales, U.K.
}
\date{\today}

\documentclass[12pt]{article}
\usepackage{amssymb}
\usepackage{graphicx}
\usepackage[all]{xy}
\usepackage[usenames]{color}

\begin{document}
\maketitle

\begin{abstract}
Braided subfactors of von Neumann algebras provide a framework for studying two dimensional conformal field theories and their modular invariants. We review this in the context of $SU(3)$ conformal field theories through corresponding $SU(3)$ braided subfactors and various subfactor invariants including spectral measures for the nimrep graphs, $A_2$-planar algebras and almost Calabi-Yau algebras.
\end{abstract}

\section{Introduction}

A statistical mechanical two dimensional lattice model at criticality can give rise to a conformal
quantum field theory when the correlation length diverges. The lattice model itself can be described via  a configuration space of a distribution of edges of a graph on two dimensional lattice with Boltzmann weights to describe local interactions. Integrable or solvable models can be constructed  from $SU(n)$ data with  Boltzmann weights satisfying a Yang-Baxter equation coming from the representation theory  of Hecke algebras. Corresponding conformal field theories can be described using $SU(n)$ loop group models either in a vertex operator algebra framework or in algebraic quantum field theory. We will illustrate
the use of subfactor theory to understand modular invariant partition functions in critical statistical mechanical models and conformal invariant field theories using a detailed exposition of $SU(3)$ models.

\subsection{Modular invariant partition functions} \label{sect:MI}

Modular invariant partition functions arise as continuum limits in statistical mechanics, i.e. letting the lattice spacing tend to zero whilst simultaneously approaching the critical temperature which can produce conformally invariant field theories. We will focus on the fundamental role of the modular invariants  in conformal field theory in the case of loop group models arising from $SU(3)$ at level $k$. The case of $SU(2)$ at level 2 is the familiar case of the Ising model.

These modular invariant partition functions have the form:
$$Z(\tau) = \sum_{\lambda, \mu \in \mathcal{A}} Z_{\lambda, \mu} \chi_{\lambda}(\tau) \chi_{\mu}(\tau)^{\ast}.$$
Here $\mathcal{A}$ denotes the irreducible symmetries in the Verlinde algebra of positive energy representations of the loop group of $SU(n)$, and  $\chi_{\lambda} = \mathrm{tr}(q^{L_0 - c/24})$ denote the character of an irreducible representation $\lambda \in \mathcal{A}$, which is the trace in the positive energy representation, where $q = e^{2 \pi i \tau}$, $\mathrm{Im}(\tau)>0$. Moreover $L_0$ is the conformal Hamiltonian which is the infinitesimal generator of the rotation group on the circle. The characters are transformed linearly under the action of $SL(2;\mathbb{Z})$, via $\left( \begin{array}{cc} 0 & -1 \\ 1 & 0 \end{array} \right): \tau \rightarrow -1/\tau$, $\left( \begin{array}{cc} 1 & 1 \\ 0 & 1 \end{array} \right): \tau \rightarrow \tau+1$, or more precisely, $\chi_{\lambda}(-1/\tau) = \sum_{\mu} S_{\lambda,\mu} \chi_{\mu}(\tau)$, $\chi_{\lambda}(\tau+1) = \sum_{\mu} T_{\lambda,\mu} \chi_{\mu}(\tau)$, where $S$ is a symmetric unitary matrix which remarkably diagonalizes the fusion rules (see (\ref{eqn:verlinde_formula})), with $S_{\lambda,0} \geq S_{0,0} > 0$, and $T$ is a diagonal finite order matrix.

The problem of the classification of modular invariants is to first find all non-negative integer valued matrices $Z$ such that $ZS = SZ$, $ZT = TZ$, subject to the constraint $Z_{0,0} = 1$ (which reflects the physical concept of the uniqueness of the vacuum state). The non-negative integer requirement on the entries of $Z$ comes from the understanding of the entries as multiplicities of the decomposition of the underlying Hilbert space. However we are interested in the structure which glues together the two chiral copies in the quadratic expression for the partition function and the corresponding conformal field theory. This extra structure is well described through the theory of braided subfactors and alpha-induction.
The Verlinde algebra $\mathcal{A}$ of the positive energy representations, or the irreducible characters is represented by a system of braided endomorphisms on a factor $N$. The modular invariant $Z$ subfactor will be represented by an inclusion $N \subset M$, which will carry this gluing information.

The trivial modular invariant, given by $Z_{\lambda, \mu} = \delta_{\lambda, \mu}$ or $Z = \sum_{\lambda} |\chi_{\lambda}|^2$ is always a solution. There may also be permutation invariants $Z = \sum_{\lambda} \chi_{\lambda} \chi_{\omega(\lambda)}^{\ast}$, where $\omega$ is a permutation of the labels which preserves the fusion rules and $\omega(0) = 0$. For a Rational Conformal Field Theory (RCFT), after extending the chiral algebras maximally, the partition function is at most a permutation matrix $Z_{\tau, \tau'}^{\mathrm{ext}} = \delta_{\tau, \omega(\tau')}$, where $\tau$, $\tau'$ label the representations of an extended chiral algebra and $\omega$ is now a permutation of these labels (see \cite{moore/seiberg:1989}). The extended characters $\chi_{\tau}^{\mathrm{ext}}$ can be decomposed in terms of the original characters $\chi_{\lambda}$ as $\chi_{\tau}^{\mathrm{ext}} = \sum_{\lambda} b_{\tau, \lambda} \chi_{\lambda}$ for some non-negative coefficients $b_{\tau, \lambda}$, and $Z_{\lambda, \mu} = \sum_{\tau} b_{\tau, \lambda} b_{\omega(\tau), \mu}$. Modular invariants for which $\omega$ is trivial are called type I, i.e. $Z^{\mathrm{ext}} = \sum_{\tau} |\chi_{\tau}^{\mathrm{ext}}|^2$, whereas those corresponding to non-trivial $\omega$ are called type II, i.e. $Z^{\mathrm{ext}} = \sum_{\tau} \chi_{\tau}^{\mathrm{ext}} (\chi_{\omega(\tau)}^{\mathrm{ext}}){^\ast}$. The matrices $Z_{\lambda, \mu}$ for these type I invariants are symmetric, whereas for the type II invariants only ``vacuum coupling'' is necessarily symmetric: $Z_{0,\lambda} = Z_{\lambda,0}$ for all $\lambda$.
There are also modular invariants which do not even have the symmetry $Z_{0,\lambda} = Z_{\lambda,0}$, see \cite{bockenhauer/evans:2000}. It can happen that two different chiral extensions $\mathcal{A} \subset \mathcal{A}^{\mathrm{ext}}_{\pm}$ of the original chiral algebra $\mathcal{A}$ are compatible such that a given coupling matrix has to be interpreted as an ``automorphism'' invariant with respect to the enhanced heterotic algebra $\mathcal{A}^{\mathrm{ext}}_+ \otimes \mathcal{A}^{\mathrm{ext}}_-$. The labels of the left and right sectors are generically different, and thus $\omega$ is no longer a permutation but rather a bijection, ``isomorphism of fusion rules'' is a more precise notion than ``automorphism'', and the distinction between diagonal and permutation invariants no longer makes sense for a maximally extended heterotic symmetry algebra.
The notion of classifying type was modified in \cite{evans/pugh:2009ii}, where type instead referred to an inclusion $N \subset M$ of factors, as will be discussed in the Section \ref{sect:braided_subfactors}.

In the cases we are interested in, the factor $N$ is obtained as a local factor $N=N(I_{\circ})$ of a conformally covariant quantum field theoretical net of factors $\{N(I)\}$ indexed by proper intervals $I \subset \mathbb{R}$ on the real line arising from current algebras defined in terms of local loop group representations, and that the $N$-$N$ system is obtained as restrictions of Doplicher-Haag-Roberts morphisms (cf. \cite{haag:1992}) to $N$.
Taking two copies of such a net and placing the real axes on the light cone, then this defines a local conformal net $\{A(\mathcal{O})\}$, indexed by double cones $\mathcal{O}$ on two-dimensional Minkowski space (cf. \cite{rehren:2000} for such constructions).
A braided subfactor $N \subset M$, determining in turn two subfactors $N \subset M_{\pm}$ obeying chiral locality, will provide two local nets of subfactors $\{N(I) \subset M_{\pm}(I) \}$.
Arranging $M_+(I)$ and $M_-(J)$ on the two light cone axes defines a local net of subfactors $\{A(\mathcal{O})\subset A_{\mathrm{ext}}(\mathcal{O})\}$ in Minkowski space.
The embedding $M_+\otimes M_-^{\mathrm{op}} \subset B$ gives rise to another net of subfactors $\{A_{\mathrm{ext}}(\mathcal{O}) \subset B(\mathcal{O})\}$, where the conformal net $\{B(\mathcal{O})\}$ satisfies locality.
As shown in \cite{rehren:2000}, there exist a local conformal two-dimensional quantum field theory such that the coupling matrix $Z$ describes its restriction to
the tensor products of its chiral building blocks $N(I)$.
There are chiral extensions $N(I) \subset M_+(I)$ and $N(I) \subset M_-(I)$ for left and right chiral nets which are indeed maximal and should
be regarded as the subfactor version of left- and right maximal extensions of the chiral algebra.

We begin in Section \ref{sect:subfactors} by considering subfactors. In Sections \ref{sect:braided_subfactors}-\ref{sect:SU(3)MI} we focus on braided subfactors and their relation to modular invariant partition functions.
In particular, in Section \ref{sect:SU(3)MI} we describe the modular invariants for $SU(3)$ at level $k$. These are labelled by the $SU(3)$ $\mathcal{ADE}$ graphs which appear as nimreps. The finite subgroups of $SU(3)$ are discussed in Section \ref{sect:subgroups-SU(3)}, including their relation to the $SU(3)$ $\mathcal{ADE}$ graphs.
In the remainder of the paper we then review a number of other invariants for braided subfactors, focusing mainly on those braided subfactors which realise the $SU(3)$ modular invariants. In Section \ref{sect:spectral_measures} we review spectral measures for the $SU(3)$ $\mathcal{ADE}$ graphs and the McKay graphs of subgroups of $SU(3)$. Then in Section \ref{sect:diagrammatic} we review various related pictorial descriptions associated to the $SU(3)$ $\mathcal{ADE}$ graphs, including a diagrammatic presentation of the $A_2$-Temperley-Lieb algebra, the $A_2$-Temperley-Lieb category, and the $A_2$-planar algebras of \cite{evans/pugh:2009iii}. Finally in Section \ref{sect:almostCYalg} we review our work on the almost Calabi-Yau algebra, which is isomorphic to the graded algebra $\Sigma \cong \Pi$ of Section \ref{sect:module_categories}, associated to a finite $SU(3)$ $\mathcal{ADE}$ graph which carries a cell system.

\subsection{Subfactors} \label{sect:subfactors}

In the last twenty years, a very fruitful circle of ideas has developed linking the theory of subfactors of von Neumann algebras with modular invariants in conformal field theory. Subfactors have been studied through their paragroups, planar algebras and have serious contact with free probability theory. The understanding and classification of modular invariants is significant for conformal field theory and their underlying statistical mechanical models. These areas are linked through the use of braided subfactors and $\alpha$-induction which in particular for $SU(2)$ subfactors and $SU(2)$ modular invariants invokes $ADE$ classifications on both sides.

A group acting outerly on a factor can be recovered from the inclusion of its fixed point algebra.
A general subfactor encodes a more sophisticated symmetry or a way of handling non group like symmetries including but going beyond quantum groups \cite{evans/kawahigashi:1998}.
The classification of subfactors was initiated by Jones \cite{jones:1983} who found that the minimal symmetry to understand the inclusion is through the Temperley-Lieb algebra.
This arises from the representation theory of $SU(2)$ or dually certain representations of Hecke algebras.
There are a number of invariants (encoding the symmetry) one can assign to a subfactor, and under certain circumstances they are complete at least for hyperfinite subfactors.
Popa \cite{popa:1995} axiomatized the inclusions of relative commutants in the Jones tower, and Jones \cite{jones:planar} showed that this was equivalent to his planar algebra description. Here one is naturally forced to work with nonamenable factors through free probabilistic constructions e.g. \cite{guionnet/jones/shlyakhtenko:2007}.

Let $A$ and $B$ be type $\mathrm{III}$ von Neumann factors. A unital $\ast$-homomorphism $\rho:A\rightarrow B$ is called a $B$-$A$ morphism.
The positive number $d_{\rho}=[B:\rho(A)]^{1/2}$ is called the statistical dimension of $\rho$; here $[B:\rho(A)]$ is the Jones-Kosaki index \cite{jones:1983, kosaki:1986} of the subfactor $\rho(A)\subset B$.
Some $B$-$A$ morphism $\rho'$ is called equivalent to $\rho$ if $\rho'=\mathrm{Ad}(u)\circ\rho$ for some unitary $u\in B$.
The equivalence class $[\rho]$ of $\rho$ is called the $B$-$A$ sector of $\rho$. If $\rho$ and $\sigma$ are $B$-$A$ morphisms with finite statistical dimensions, then the vector space of intertwiners
$$\mathrm{Hom}(\rho,\sigma)=\{ t\in B: t\rho(a)=\sigma(a)t \,, \,\, a\in A \}$$
is finite-dimensional, and we denote its dimension by $\langle \rho, \sigma \rangle$.
A $B$-$A$ morphism $\rho$ with finite statistical dimension is called irreducible if $\langle \rho,\rho \rangle=1$, i.e. if $\mathrm{Hom}(\rho,\rho) = \mathbb{C} \mathbf{1}_B$.
Then, if $\langle \rho, \tau \rangle \neq 0$ for some (possibly reducible) $B$-$A$ morphism $\tau$, then $[\rho]$ is called an irreducible subsector of $[\tau]$ with multiplicity $\langle \rho, \tau \rangle$.
An irreducible $A$-$B$ morphism $\overline{\rho}$ is a conjugate morphism of the irreducible $\rho$ if and only if $[\overline{\rho}\rho]$ contains the trivial sector $[\mathrm{id}_A]$ as a subsector, and then $\langle \rho\overline{\rho}, \mathrm{id}_B \rangle = 1 = \langle \overline{\rho}\rho, \mathrm{id}_A \rangle$ automatically \cite{izumi:1991}.

For an $M$-$M$ endomorphism $\rho$ where the subfactor $N = \rho(M) \subset M$ has finite index, the inclusion $\rho(M) \subset M$ can be extended upwards (called the \emph{Jones tower}) or downwards (called the \emph{Jones tunnel}):
$$\begin{array}{ccccccccccc}
\cdots \subset & \rho\overline{\rho}(M) & \subset & \rho(M) & \subset & M & \subset & M \otimes_N M & \subset & M \otimes_N M \otimes_N M & \subset \cdots \\
& \| & & \| & & \| & & \| & & \| \\
& M_{-2} & & N = M_{-1} & & M_0 & & M_1 & & M_2
\end{array}$$
The tunnel construction is generally used when working with type $\mathrm{III}$ factors, whereas the tower construction is generally used when working with type $\mathrm{II}_1$ factors, using the Jones basic construction $N \subset M \subset M_1 = \langle M,e \rangle$ to adjoin an extra projection $e$ arising from the projection of $M$ onto $N$. In the type $\mathrm{III}$ setting the projection $e$ arises from the conditional expectation of $M$ onto $N$.
These Jones projections satisfy the Temperley-Lieb relations of integrable statistical mechanics. They are contained in the tower of relative commutants of any finite index subfactor and are in some sense the minimal symmetries.

For a subfactor with finite index, the relative commutants $M_k' \cap M_l$, for $M_k \subset M_l$, are finite dimensional. We are interested in the case of finite depth subfactors, where the inclusions in the sequence of commuting squares below are all described by finite graphs:
\begin{equation} \label{eqn:seq_rel_commutant}
\begin{array}{cccccccc}
M' \cap M_k & \subset & M' \cap M_{k+1} & \subset & M' \cap M_{k+2} & \subset & \cdots \qquad \longrightarrow & A \\
\cap & & \cap & & \cap & & & \cap \\
N' \cap M_k & \subset & N' \cap M_{k+1} & \subset & N' \cap M_{k+2} & \subset & \cdots \qquad \longrightarrow & B
\end{array}
\end{equation}
The inclusions in the above commuting squares are given by two graphs, the principal graph and dual principal graph.
The principal graph for the subfactor $N \subset M$ has vertices for the $N$-$N$ and $N$-$M$ bimodules, and an edge
between $Y$ and $Z$ for each copy of $Z$ appearing inside $Y \otimes_N {}_N M_M$ or $Y \otimes_M {}_M M_N$ as appropriate.
The dual principal graph is defined similarly, but with vertices for the $M$-$N$ and $M$-$M$ bimodules.
The $\mathrm{II}_1$ factors $A$, $B$ are given by the GNS-completion of $\bigcup_{k\geq 0} M' \cap M_k$, $\bigcup_{k\geq 0} N' \cap M_k$ respectively, with respect to the trace on (\ref{eqn:seq_rel_commutant}). The subfactor $A \subset B$ has finite depth as it has the same standard invariant as $N \subset M$. It is a hyperfinite type $\mathrm{II}_1$ subfactor, so the best one could expect is that $A \otimes N \subset B \otimes N \cong N \subset M$. Indeed, this is the case by Popa. For a precise formulation see \cite{popa:1995ii} and references therein.

Consider a subfactor $N \subset M$ with principal graph $\mathcal{G}$.
The relative commutants $(\rho\overline{\rho} \cdots \rho M)' \cap M$ in the tunnel, and $\mathrm{End}({}_N \underbrace{M \otimes_N M \otimes_N \cdots \otimes_N M}_{k \textrm{\scriptsize{ copies of M}}} {}_M)$ in the tower, can be identified with $\bigoplus_v \ast \cdot \mathrm{End} \left( (\mathbb{C}\mathcal{G})_{2k-1} \right) \cdot v$, where $\ast$ is the distinguished vertex of $\mathcal{G}$ with lowest Perron-Frobenius weight, the direct sum is over all vertices $v$ of $\mathcal{G}$, and $(\mathbb{C}\mathcal{G})_{n}$ denotes the space of all paths of length $n$ on $\mathcal{G}$.
Similarly, the relative commutants $(\rho\overline{\rho} \cdots \overline{\rho} M)' \cap M$ and $\mathrm{End}({}_N \underbrace{M \otimes_N M \otimes_N \cdots \otimes_N M}_{k \textrm{\scriptsize{ copies of M}}} {}_N)$ can be identified with $\bigoplus_v \ast \cdot \mathrm{End} \left( (\mathbb{C}\mathcal{G})_{2k} \right) \cdot v$.
The path Hilbert space $\mathbb{C}\mathcal{G}$ is a graded algebra where multiplication $(\mathbb{C}\mathcal{G})_i \times (\mathbb{C}\mathcal{G})_j \rightarrow (\mathbb{C}\mathcal{G})_{i+j}$ of two paths $x \in (\mathbb{C}\mathcal{G})_i$, $y \in (\mathbb{C}\mathcal{G})_j$ is given by concatenation of paths $xy$, and is defined to be zero if $r(x) \neq s(y)$, where $s(x)$, $r(x)$ denote the source, range vertices of the path $x$ respectively.
The endomorphisms $\mathrm{End} \left( (\mathbb{C}\mathcal{G})_n \right)$ on $(\mathbb{C}\mathcal{G})_n$ are the matrices with rows and columns labeled by the paths of length $n$ on $\mathcal{G}$.
The $\mathrm{End} \left( (\mathbb{C}\mathcal{G})_n \right)$ have an algebra structure given by matrix multiplication.
Thus there are two different notions of path algebra of $\mathcal{G}$. In the theory of operator algebras, the path algebra of $\mathcal{G}$ is usually $\bigoplus_{n \geq 0} \mathrm{End} \left( (\mathbb{C}\mathcal{G})_n \right)$ \cite[Section 2.9]{evans/kawahigashi:1998}. In this paper however, we work with the graded algebra $\mathbb{C}\mathcal{G}$.

\begin{figure}[tb]
\begin{minipage}[t]{7.9cm}
\begin{center}
  \includegraphics[width=60mm]{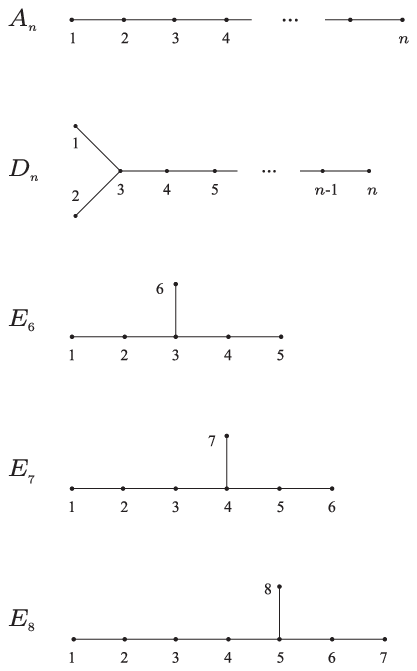}
 \caption{$ADE$ Dynkin diagrams} \label{fig-ADE}
\end{center}
\end{minipage}
\hfill
\begin{minipage}[t]{7.9cm}
\begin{center}
  \includegraphics[width=60mm]{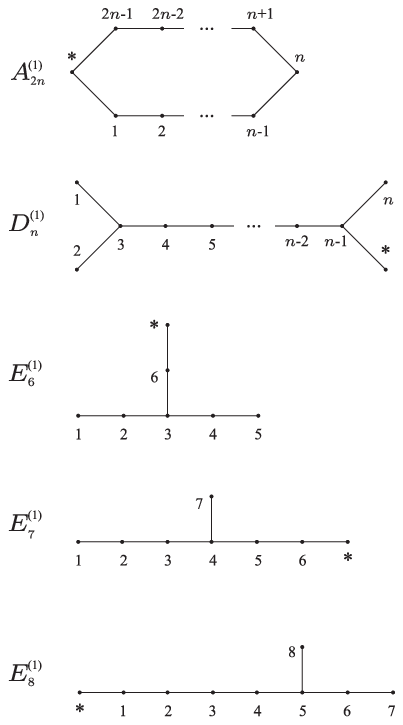}
 \caption{Affine $ADE$ Dynkin diagrams} \label{fig-AffineADE}
\end{center}
\end{minipage}
\end{figure}

\begin{figure}[tb]
\begin{center}
  \includegraphics[width=45mm]{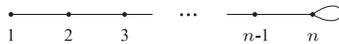}
 \caption{Tadpole graph $T_n$, $n \geq 1$} \label{fig-T}
\end{center}
\end{figure}

The subfactors with index less than or equal to 4 have been classified for a long time \cite{goodman/de_la_harpe/jones:1989, ocneanu:1988, izumi:1991}.
For index less than 4 they have an $ADE$ classification, where the $ADE$ graphs, the simply laced Dynkin diagrams illustrated in Figure \ref{fig-ADE}, appear as the principal graph.
There are two infinite families, $A_n$, $n \geq 1$, and $D_n$, $n \geq 4$, and three exceptional graphs $E_6$, $E_7$ and $E_8$. The graph $D_n$ may be obtained from $A_{2n-3}$ by a $\mathbb{Z}_2$-orbifold procedure, and vice-versa.
There is also a separate family of graphs, the tadpole graphs $T_n$, $n \geq 1$, illustrated in Figure \ref{fig-T}, which are obtained by taking a $\mathbb{Z}_2$-orbifold of an even $A$ graph $A_{2n}$, and vice-versa.
However the Dynkin diagrams $D_{\mathrm{odd}}$ and $E_7$ cannot appear as the principal graph of a subfactor \cite{ocneanu:1988, izumi:1991}.
For both $\mathcal{G} = E_6$ and $\mathcal{G} = E_8$ there are up to isomorphism two non-isomorphic subfactors with principal graph $\mathcal{G}$, due to the existence of two inequivalent connections, which are the complex conjugate of each other.

These $ADE$ graphs appear in a variety of other contexts. They classify the simply-laced Lie groups, which are the Lie groups such that all the nonzero roots of the corresponding Lie algebra have equal length.
They also classify the quotient Kleinian singularities $\mathbb{C}^2/\Gamma$ for finite subgroups $\Gamma \subset SU(2)$ \cite{reid:2002}.
Along with the tadpole graphs $T_n$ they classify the symmetric non-negative integer matrices of norm $<2$, that is, any such matrix must be the adjacency matrix of one of the $ADET$ graphs. The tadpole graphs may be ruled out if one imposes the condition of bipartiteness or 2-colourability.

For the subfactor $R = \otimes_{\mathbb{N}} M_2 \subset R \otimes M_2$ at index 4, of type II or III depending on how the factor $R$ is completed, we can form the tower $R \subset R \otimes M_2 \subset R \otimes M_2 \otimes M_2 \subset \cdots $, and the relative commutants yield the tower $\mathbb{C} \subset M_2 \subset M_2 \otimes M_2 \subset \cdots $.
Taking the fixed point algebra for the natural action of a subgroup $G \subset SU(2)$ one obtains the subfactor $R^G \subset (R \otimes M_2)^G$. Then we have the tower $R^G \subset (R \otimes M_2)^G \subset (R \otimes M_2 \otimes M_2)^G \subset \cdots $, and the relative commutants yield the tower $\mathbb{C} \subset M_2^G \subset (M_2 \otimes M_2)^G \subset \cdots $ \cite{wassermann:1988} (see also \cite{goodman/de_la_harpe/jones:1989}).
The finite subgroups of $SU(2)$, the binary polyhedral groups, have themselves an $ADE$ classification, given in Table \ref{table:ADE-affineADE}.
The McKay graphs of the representation theory of the finite subgroups of $SU(2)$ are the simply laced affine Dynkin diagrams, illustrated in Figure \ref{fig-AffineADE}, where $\ast$ denotes the identity representation. For each affine Dynkin diagram, the corresponding $ADE$ Dynkin diagram is obtained from the affine diagram by deleting the vertex $\ast$ and all edges attached to it.
The principal graph for the subfactor $R^G \subset (R \otimes M_2)^G$ is the affine Dynkin diagram associated to the finite subgroup $G \subset SU(2)$.
This subfactor was constructed in \cite{goodman/de_la_harpe/jones:1989}. See \cite{izumi/kawahigashi:1993} for a more general construction.
When $G = SU(2)$ itself, the fixed point algebra $(\otimes^p M_2)^G$ is isomorphic to the Temperley-Lieb algebra generated by operators $U_i \in \mathrm{End}(\otimes^p \mathbb{C}^2)$, for $i=1,\ldots,p-1$, and the principal graph of $R^{SU(2)} \subset (R \otimes M_2)^{SU(2)}$ is the infinite Dynkin diagram $A_{\infty}$. In the case where $G = \mathbb{T}$ is the infinite subgroup of $SU(2)$ given by the torus, the principal graph of $R^{\mathbb{T}} \subset (R \otimes M_2)^{\mathbb{T}}$ is the doubly-infinite Dynkin diagram $A_{\infty,\infty}$. In the case where $G = \mathbb{T} \rtimes \mathbb{Z}_2$, where the $\mathbb{Z}_2$ action sends $\mathbb{T} \ni t \mapsto \overline{t}$, , the principal graph of $R^{\mathbb{T} \rtimes \mathbb{Z}_2} \subset (R \otimes M_2)^{\mathbb{T} \rtimes \mathbb{Z}_2}$ is the infinite Dynkin diagram $D_{\infty}$.

\begin{table}[tb]
\begin{center}
\begin{tabular}{|c|c|c|c|} \hline
Dynkin Diagram $\mathcal{G}$ & Type & Subgroup $\Gamma \subset SU(2)$
& $|\Gamma|$ \\
\hline\hline $A_l$ & I & binary cyclic, $\mathbb{Z}_{l+1}$ & $l+1$ \\
\hline $D_{2k}$ & I & binary dihedral, $BD_{2k} = Q_{2k-2}$ & $8k-8$ \\
\hline $D_{2k+1}$ & II & binary dihedral, $BD_{2k+1} = Q_{2k-1}$ & $8k-4$ \\
\hline $E_6$ & I & binary tetrahedral, $BT = BA_4$ & 24 \\
\hline $E_7$ & II & binary octahedral, $BO = BS_4$ & 48 \\
\hline $E_8$ & I & binary icosahedral, $BI = BA_5$ & 120 \\
\hline
\end{tabular}\\
\caption{The $ADE$ classification of the finite subgroups of $SU(2)$} \label{table:ADE-affineADE}
\end{center}
\end{table}

The principal graph of any subfactor at index 4 must be one of these affine Dynkin diagrams, the infinite Dynkin graphs $A_{\infty}$, $D_{\infty}$, or the doubly infinite Dynkin graph $A_{\infty,\infty}$.
The classification of subfactors has now been extended to index up to 5 \cite{haagerup:1994, bisch:1998, asaeda/haagerup:1999, asaeda/yasuda:2009, bigelow/morrison/peters/snyder:2009, morrison/snyder:2011, morrison/penneys/peters/snyder:2011, izumi/jones/morrison/snyder:2011, penneys/tener:2011}.
In our work we are interested in a related graph, called the nimrep graph (see Section \ref{sect:braided_subfactors}), where, unlike with the principal graphs, the Dynkin diagrams $D_{\mathrm{odd}}$ and $E_7$ do appear.

\subsection{Braided subfactors and modular invariants} \label{sect:braided_subfactors}

As our starting point, the Verlinde algebra is realised by systems of endomorphisms ${}_N \mathcal{X}_N$ of a hyperfinite type $\mathrm{III}$ factor $N$.
That is, ${}_N \mathcal{X}_N$ denotes a finite system of finite index irreducible endomorphisms of a factor $N$ in the sense that different elements of ${}_N \mathcal{X}_N$ are not unitary equivalent, for any $\lambda \in {}_N \mathcal{X}_N$ there is a representative $\overline{\lambda} \in {}_N \mathcal{X}_N$ of the conjugate sector $[\overline{\lambda}]$, and ${}_N \mathcal{X}_N$ is closed under composition and subsequent irreducible decomposition.
In the case of Wess-Zumino-Witten (WZW) models associated to $SU(n)$ at level $k$, the Verlinde algebra is a non-degenerately braided system of endomorphisms ${}_N \mathcal{X}_N$, labelled by the positive energy representations of the loop group of $SU(n)_k$ on a type $\mathrm{III}_1$ factor $N$, with fusion rules $\lambda \mu = \bigoplus_{\nu} N_{\lambda \nu}^{\mu} \nu$ which exactly match those of the positive energy representations \cite{wassermann:1998}. The fusion matrices $N_{\lambda} = [N_{\rho \lambda}^{\sigma}]_{\rho,\sigma}$ are a family of commuting normal matrices which give a representation themselves of the fusion rules of the positive energy representations of the loop group of $SU(n)_k$, $N_{\lambda} N_{\mu} = \sum_{\nu} N_{\lambda \nu}^{\mu} N_{\nu}$.
This family $\{ N_{\lambda} \}$ of fusion matrices can be simultaneously diagonalised:
\begin{equation} \label{eqn:verlinde_formula}
N_{\lambda} = \sum_{\sigma} \frac{S_{\sigma, \lambda}}{S_{\sigma,0}} S_{\sigma} S_{\sigma}^{\ast},
\end{equation}
where $0$ is the trivial representation, and the eigenvalues $S_{\sigma, \lambda}/S_{\sigma,0}$ and eigenvectors $S_{\sigma} = [S_{\sigma, \mu}]_{\mu}$ are described by the statistics $S$ matrix.
Moreover, there is equality between the statistics $S$- and $T$- matrices and the Kac-Peterson modular $S$- and $T$- matrices which perform the conformal character transformations \cite{kac:1990}, thanks to \cite{frohlich/gabbiani:1990, fredenhagen/rehren/schroer:1992, wassermann:1998}.

The key structure in the conformal field theory is the modular invariant partition function $Z$. In the subfactor setting this is realised by
a braided subfactor $N \subset M$ where trivial (or permutation) invariants in the ambient factor $M$ when restricted to $N$ yield $Z$. This would mean that the dual canonical endomorphism is in $\Sigma({}_N \mathcal{X}_N)$, i.e. decomposes as a finite linear combination of endomorphisms in ${}_N \mathcal{X}_N$.
Indeed if this is the case for the inclusion $N \subset M$, then the process of $\alpha$-induction allows us to analyse the modular invariant,
providing two extensions of $\lambda$ on $N$ to endomorphisms $\alpha^{\pm}_{\lambda}$ of $M$, such that the matrix $Z_{\lambda,\mu} = \langle \alpha_{\lambda}^+, \alpha_{\mu}^- \rangle$ is a modular invariant \cite{bockenhauer/evans/kawahigashi:1999, bockenhauer/evans:2000, evans:2003}.

\begin{figure}[tb]
\begin{center}
  \includegraphics[width=160mm]{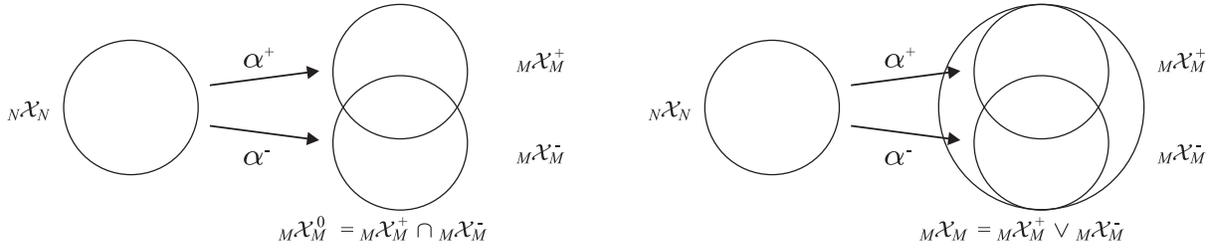}
 \caption{$\alpha$-induction}
\end{center}
\end{figure}

Let ${}_N \mathcal{X}_M$, ${}_M \mathcal{X}_M$ denote a system of endomorphisms consisting of a choice of representative endomorphism of each irreducible subsector of sectors of the form $[\lambda \overline{\iota}]$, $[\iota \lambda \overline{\iota}]$ respectively, for each $\lambda \in {}_N \mathcal{X}_N$, where $\iota: N \hookrightarrow M$ is the inclusion map which we may consider as an $M$-$N$ morphism, and $\overline{\iota}$ is a representative of its conjugate $N$-$M$ sector.
The action of the system ${}_N \mathcal{X}_N$ on the $N$-$M$ sectors ${}_N \mathcal{X}_M$ produces a \emph{nimrep} (non-negative matrix integer representation of the fusion rules) $G_{\lambda} G_{\mu} = \sum_{\nu} N_{\lambda \nu}^{\mu} G_{\nu}$,
whose spectrum reproduces exactly the diagonal part of the modular invariant, i.e.
\begin{equation} \label{eqn:verlinde_formulaG}
G_{\lambda} = \sum_{\sigma} \frac{S_{\sigma,\lambda}}{S_{\sigma,0}} \psi_{\sigma} \psi_{\sigma}^{\ast},
\end{equation}
with the spectrum of $G_{\lambda}$ given by \cite[Theorem 4.16]{bockenhauer/evans/kawahigashi:2000}
\begin{equation} \label{eqn:spectrumG}
G_{\lambda} = \{ S_{\mu, \lambda}/S_{\mu,0} \textrm{ with multiplicity } Z_{\mu,\mu} \}
\end{equation}
The labels $\mu$ of the non-zero diagonal elements are called the exponents of $\mathcal{G} = G_{\rho}$, counting multiplicity. A modular invariant for which there exists a nimrep whose spectrum is described by the diagonal part of the invariant is said to be nimble.
The zero-column of the modular invariant $Z$ associated with the subfactor $N \subset M$ determines $\langle \alpha^+_j, \alpha^+_{j'} \rangle$.
Then for all modular invariants with the same zero-column, the sectors $[\alpha^{\pm}_0]$
have the same nimrep graphs. Let $v$ be an isometry which intertwines the identity and the canonical endomorphism $\gamma = \iota \overline{\iota}$. Proposition 3.2 in \cite{bockenhauer/evans:2000} states that the following conditions are equivalent:
\begin{itemize}
\item[1.] $Z_{\lambda,0} = \langle \theta, \lambda \rangle$ for all $\lambda \in {}_N \mathcal{X}_N$.
\item[2.] $Z_{0,\lambda} = \langle \theta, \lambda \rangle$ for all $\lambda \in {}_N \mathcal{X}_N$.
\item[3.] Chiral locality holds: $\varepsilon^+(\theta,\theta)v^2 = v^2$.
\end{itemize}
The chiral locality condition, which can be expressed in terms of the single inclusion $N \subset M$ and the braiding, expresses local commutativity (locality) of the extended net, if $N \subset M$ arises from a net of subfactors \cite{longo/rehren:1995}. Chiral locality holds if and only if the dual canonical endomorphism is visible in the vacuum row, $[\theta] = \bigoplus_{\lambda} Z_{0,\lambda} [\lambda]$ (and hence in the vacuum column also).
An inclusion $N \subset M$ is type I if and only if one of the above equivalent conditions 1-3 hold, otherwise it is type II \cite{evans/pugh:2009ii}.
It was shown in \cite{evans/pugh:2009ii} that it is possible for a type I modular invariant to be realized by a type II inclusion, and suggests that care needs to be taken with the type I, II labelling of modular invariants. The situation is in fact more subtle, with the possibility of a modular invariant being neither type I nor type II, as explained in Section \ref{sect:MI}. However, such a modular invariant did not appear in our treatment of $SU(3)$ modular invariants in \cite{evans/pugh:2009ii}.

Behrend, Pearce, Petkova and Zuber \cite{behrend/pearce/petkova/zuber:2000} (see also \cite{zuber:2002}) systematically proposed nimreps as a framework for boundary conformal field theory. The $N$-$M$ system ${}_N \mathcal{X}_M$ corresponds to boundary fields in their language, and the $M$-$M$ system ${}_M \mathcal{X}_M$ to defect lines.

The systems ${}_N \mathcal{X}_N$, ${}_N \mathcal{X}_M$, ${}_M \mathcal{X}_M$ are (the irreducible objects of) tensor categories of endomorphisms with the Hom-spaces as their morphisms. Thus ${}_N \mathcal{X}_N$ gives a braided modular tensor category, and ${}_N \mathcal{X}_M$ a module category.
The structure of the module category ${}_N \mathcal{X}_M$ is the same as a tensor functor $F$ from ${}_N \mathcal{X}_N$ to the category $\mathrm{Fun}({}_N \mathcal{X}_M,{}_N \mathcal{X}_M)$ of additive functors from ${}_N \mathcal{X}_M$ to itself, see \cite{ostrik:2003}. That is, $F$ is essentially the module category ${}_N \mathcal{X}_M$.
Ostrik \cite{kirillov/ostrik:2002, ostrik:2003} took up a categorical description of subfactor $\alpha$-induction, see \cite[Remark 14]{ostrik:2003}, and this was taken further by Fjelstad, Fr\"{o}hlich, Fuchs, Schweigert and Runkel as a categorical framework for conformal field theory. See \cite{fuchs/runkel/schweigert:2010} for a review.

\subsection{The generalized Temperley-Lieb algebras} \label{sect:Hecke}

The Jones basic construction $M_{i-1} \subset M_i \subset M_{i+1}$ is through adjoining an extra projection $e_i$ arising from the projection or conditional expectation of $M_i$ onto $M_{i-1}$. These projections satisfy the Temperley-Lieb relations of integrable statistical mechanics. They are contained in the tower of relative commutants of any finite index subfactor and are in some sense the minimal symmetries.
The Temperley-Lieb algebra has a realization from $SU(2)$, from the fixed point algebras of quantum $SU(2)$ on the Pauli algebra and special representations of Hecke algebras of type $A$. We now review these Hecke algebras.

The subfactor $R^G \subset (R \otimes M_2)^G$ in the $SU(2)$ setting of Section \ref{sect:subfactors} can be generalized to $SU(n)$, $n \geq 2$.
Let $M_n = \mathrm{End}(\mathbb{C}^n)$. By Weyl duality, the fixed point algebra of $\otimes^m M_n$ under the product adjoint action of $SU(n)$ is generated by a representation $\sigma \rightarrow g_{\sigma}$ on $\otimes^m \mathbb{C}^n$ of the group ring of the symmetric, or permutation, group $S_m$. This algebra is generated by unitary operators $g_j$, $j=1,\ldots,m-1$, which represent transpositions $(j,j+1)$, satisfying the relations
\begin{eqnarray}
g_j^2 & = & 1, \label{eqn:Hecke-1} \\
g_i g_j & = & g_j g_i, \quad |i-j|>1, \label{eqn:Hecke-2} \\
g_i g_{i+1} g_i & = & g_{i+1} g_i g_{i+1}, \label{eqn:braiding_relation}
\end{eqnarray}
and the vanishing of the antisymmetrizer
\begin{equation}\label{SU(N)condition}
\sum_{\sigma \in S_m} \mathrm{sgn}(\sigma) g_{\sigma} = 0.
\end{equation}
Writing $g_j = 1 - U_j$, these unitary generators and relations are equivalent to the self-adjoint generators $\mathbf{1}, U_j$, $j=1, \ldots, m-1$, and relations
\begin{center}
\begin{minipage}[b]{15cm}
 \begin{minipage}[b]{2cm}
  \begin{eqnarray*}
  \textrm{H1:}\\
  \textrm{H2:}\\
  \textrm{H3:}
  \end{eqnarray*}
 \end{minipage}
 \hspace{2cm}
 \begin{minipage}[b]{7cm}
  \begin{eqnarray*}
  U_i^2 & = & \delta U_i,\\
  U_i U_j & = & U_j U_i, \quad |i-j|>1,\\
  U_i U_{i+1} U_i - U_i & = & U_{i+1} U_i U_{i+1} - U_{i+1},
  \end{eqnarray*}
 \end{minipage}
\end{minipage}
\end{center}
where $\delta = 2$, and the analogue of (\ref{SU(N)condition}).

There is a $q$-version of this algebra, which is a representation of a Hecke algebra. This is the centralizer of a representation of the quantum group $SU(n)_q$ (or the universal enveloping algebra), with a deformation of (\ref{eqn:Hecke-1}) to
\begin{equation} \label{eqn:Hecke-1'}
(q^{-1} - g_j)(q + g_j) = 0.
\end{equation}
The invertible generators $g_j$, $j=1,\ldots,m-1$, satisfy the relations (\ref{eqn:Hecke-2}), (\ref{eqn:braiding_relation}), (\ref{eqn:Hecke-1'}) and the vanishing of the $q$-antisymmetrizer \cite{di_francesco/zuber:1990}
\begin{equation}\label{SU(N)q condition}
\sum_{\sigma \in S_m} (-q)^{|I_{\sigma}|} g_{\sigma} = 0,
\end{equation}
where $g_{\sigma} = \prod_{i \in I_{\sigma}} g_i$ if $\sigma = \prod_{i \in I_{\sigma}} \tau_{i,i+1}$.
Then writing $g_j = q^{-1} - U_j$, we are interested in the generalized Temperley-Lieb algebra generated by self-adjoint operators $\mathbf{1}, U_j$, $j=1, \ldots, m-1$, satisfying H1-H3 and the analogue of (\ref{SU(N)q condition}), where now $\delta = q+q^{-1}$. For $SU(n)$ at level $k$, $q = e^{2\pi i/(n+k)}$.
In the case $n=2$, (\ref{SU(N)q condition}) reduces for $SU(2)$ to the Temperley-Lieb condition
\begin{equation}
U_i U_{i \pm 1} U_i - U_i = 0. \label{eqn:SU(2)q_condition}
\end{equation}

There are minor errors in a parallel discussion in Section 2 of the published version of \cite{evans/pugh:2009iii} which have been corrected in the arXiv version and in \cite{evans/pugh:2010ii}.

\subsection{$SU(2)$ modular invariants} \label{sect:SU(2)MI}

We first turn briefly to the case of $SU(2)$ at level $k$, where $q$ is a $k+2^{\mathrm{th}}$ root of unity.
The classification of $SU(2)$ modular invariants is due to Cappelli, Itzykson and Zuber \cite{cappelli/itzykson/zuber:1987ii}.
They label the modular invariant with a finite $ADE$ graph $\mathcal{G}$, such that the diagonal part $Z_{\mu,\mu}$ of the invariant is exactly the multiplicity of the eigenvalue $S_{\mu,\rho}/S_{\mu,1}$ of $\mathcal{G}$, where 1, $\rho$ denote the trivial, fundamental representations respectively.
The second column in Table \ref{table:ADE-affineADE} indicates the type of the associated modular invariant.

In the subfactor theory, this is understood in the following way. Suppose $N \subset M$ is a braided subfactor which realises the modular invariant $Z_{\mathcal{G}}$. Evaluating the nimrep $G$ at the fundamental representation $\rho$, we obtain for the inclusion $N \subset M$ a matrix $G_{\rho}$, which is the adjacency matrix for the $ADE$ graph $\mathcal{G}$ which labels the modular invariant and with spectrum given by (\ref{eqn:spectrumG}).
Every $SU(2)$ modular invariant is realised by a braided subfactor, and all nimreps are realised by subfactors \cite{ocneanu:2000ii, ocneanu:2002, xu:1998, bockenhauer/evans:1999i, bockenhauer/evans:1999ii, bockenhauer/evans/kawahigashi:1999, bockenhauer/evans/kawahigashi:2000}, apart from the tadpole nimreps of the orbifolds of the even $A$'s (see e.g. \cite{bockenhauer/evans:2001} for an explanation of the failure of the tadpole nimreps).

\subsection{$SU(3)$ modular invariants} \label{sect:SU(3)MI}

The classification of $SU(3)$ modular invariants was shown to be complete by Gannon \cite{gannon:1994}, and the complete list is reproduced in \cite{evans/pugh:2009ii}.
The four infinite families are the identity invariants $Z_{\mathcal{A}^{(m)}}$, the orbifold invariants $Z_{\mathcal{D}^{(m)}}$, the conjugate invariants $Z_{\mathcal{A}^{(m)\ast}}$, and the conjugate orbifold invariants $Z_{\mathcal{D}^{(m)\ast}}$. There are also six exceptional invariants.
We explicitly write out all the modular invariants at level 5 as an example.
Let $\mathcal{P}^{(m)} = \{ \mu=(\mu_1,\mu_2) \in \mathbb{Z}^2 | \mu_1,\mu_2 \geq 0; \mu_1 + \mu_2 \leq m-3 \}$. These $\mu$ are the admissible representations of the Ka\v{c}-Moody algebra $su(3)^{\wedge}$ at level $k = m-3$. We define the automorphism $A$ of order 3 on the weights $\mu \in \mathcal{P}^{(m)}$ by $A(\mu_1,\mu_2) = (m-3-\mu_1-\mu_2,\mu_1)$.
The invariants at level 5 are $Z_{\mathcal{A}^{(8)}} = \sum |\chi_{\mu}|^2$, $Z_{\mathcal{D}^{(8)}} = \sum \chi_{\mu} \chi_{A^{5(\mu_1-\mu_2)} \mu}^{\ast}$, $Z_{\mathcal{A}^{(8)\ast}} = \sum \chi_{\mu} \chi_{\overline{\mu}}^{\ast}$, $Z_{\mathcal{D}^{(8)\ast}} = \sum \chi_{\mu} \chi_{\overline{A^{5(\mu_1-\mu_2)} \mu}}^{\ast}$, where in each case the sum is over all $\mu \in P^{(8)}_+$, and the two exceptional invariants
\begin{eqnarray}
Z_{\mathcal{E}^{(8)}} & = & |\chi_{(0,0)}+\chi_{(2,2)}|^2 + |\chi_{(0,2)}+\chi_{(3,2)}|^2 + |\chi_{(2,0)}+ \chi_{(2,3)}|^2 + |\chi_{(2,1)}+\chi_{(0,5)}|^2 \nonumber \\
& & \;\; + |\chi_{(3,0)}+\chi_{(0,3)}|^2 + |\chi_{(1,2)}+\chi_{(5,0)}|^2, \label{Z(E8)} \\
Z_{\mathcal{E}^{(8)\ast}} & = & |\chi_{(0,0)}+\chi_{(2,2)}|^2 + (\chi_{(0,2)}+\chi_{(3,2)})(\chi_{(2,0)}^{\ast}+\chi_{(2,3)}^{\ast}) \nonumber \\
& & \;\; + (\chi_{(2,0)}+\chi_{(2,3)})(\chi_{(0,2)}^{\ast}+\chi_{(3,2)}^{\ast}) + (\chi_{(2,1)}+\chi_{(0,5)})(\chi_{(1,2)}^{\ast}+\chi_{(5,0)}^{\ast}) \nonumber \\
& & \;\; + |\chi_{(3,0)}+\chi_{(0,3)}|^2 + (\chi_{(1,2)}+\chi_{(5,0)})(\chi_{(2,1)}^{\ast}+\chi_{(0,5)}^{\ast}). \label{Z(E8star))}
\end{eqnarray}

Ocneanu claimed \cite{ocneanu:2000ii, ocneanu:2002} that all $SU(3)$ modular invariants were realised by subfactors and this was shown in \cite{xu:1998, bockenhauer/evans:1999i, bockenhauer/evans:1999ii, bockenhauer/evans/kawahigashi:1999, bockenhauer/evans:2001, bockenhauer/evans:2002, evans/pugh:2009i, evans/pugh:2009ii}.
The modular invariants arising from $SU(3)_k$ conformal embeddings are (see \cite{evans:2003}):
\begin{itemize}
\item $Z_{\mathcal{D}^{(6)}}$: $SU(3)_3 \subset SO(8)_1$, also realised as an orbifold $SU(3)_3 / \mathbb{Z}_3$,
\item $Z_{\mathcal{E}^{(8)}}$: $SU(3)_5 \subset SU(6)_1$, plus its conjugate $Z_{\mathcal{E}^{(8)\ast}}$,
\item $Z_{\mathcal{E}_1^{(12)}}$: $SU(3)_9 \subset (\mathrm{E}_6)_1$,
\item $Z_{\mathcal{E}^{(24)}}$: $SU(3)_{21} \subset (\mathrm{E}_7)_1$.
\end{itemize}
The Moore-Seiberg invariant $Z_{\mathcal{E}_{5}^{(12)}}$ \cite{moore/seiberg:1989}, an automorphism of the orbifold invariant $Z_{\mathcal{D}^{(12)}} = SU(3)_9 / \mathbb{Z}_3$, is the $SU(3)$ analogue of the $E_7$ invariant for $SU(2)$, which is an automorphism of the orbifold invariant $Z_{D_{10}} = SU(2)_{16}/\mathbb{Z}_2$ (see Section 5.3 of \cite{bockenhauer/evans/kawahigashi:2000} for a realisation by a braided subfactor).
The conjugate Moore-Seiberg invariant $Z_{\mathcal{E}_4^{(12)}}$ is given by a product of the Moore-Seiberg invariant and the conjugate invariant $Z_{\mathcal{A}^{(12)\ast}}$, which were both realised by subfactors in \cite{evans/pugh:2009ii}. Then \cite[Theorem 3.6]{evans/pinto:2003ii} shows that the conjugate Moore-Seiberg invariant $Z_{\mathcal{E}_4^{(12)}}$ is also realised by a subfactor.
In Table \ref{Table:ADE-subgroupsSU(3)}, Type denotes the type of the inclusion found in \cite{evans/pugh:2009ii} which yielded the $\mathcal{ADE}$ graph as a nimrep.
Caution should be taken regarding the type for the graph $\mathcal{E}_4^{(12)}$. Although we have not yet shown that $\mathcal{E}_4^{(12)}$ is the nimrep obtained from an inclusion, it was shown in \cite{evans/pugh:2009ii} that such an inclusion would be of type II.

\begin{figure}[tbp]
\begin{center}
  \includegraphics[width=160mm]{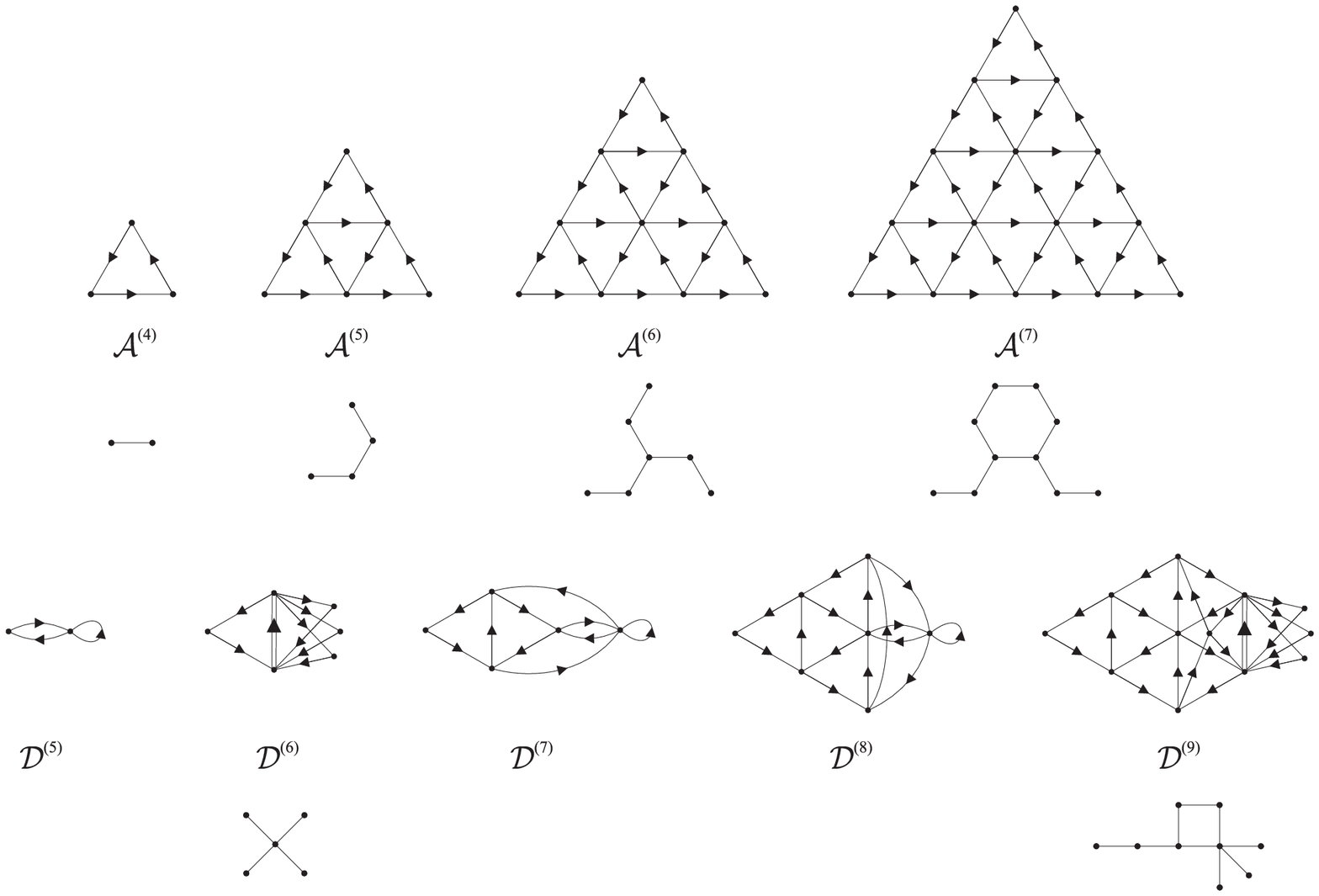}
 \caption{$SU(3)$ graphs: infinite $\mathcal{A}$, $\mathcal{D}$ families and their 0-1 parts} \label{fig-SU(3)AD}
\end{center}
\end{figure}

\begin{figure}[tbp]
\begin{center}
  \includegraphics[width=160mm]{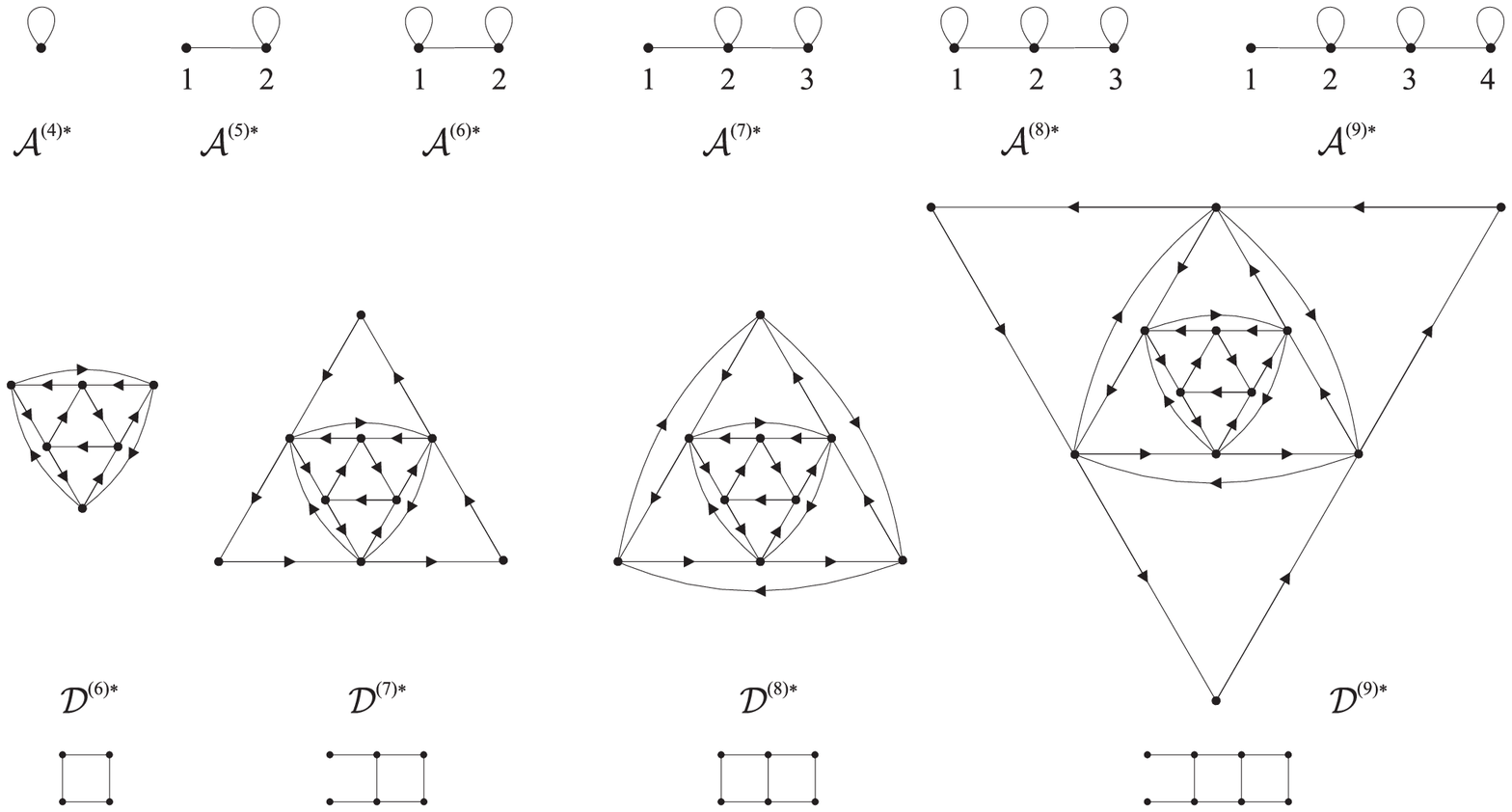}
 \caption{$SU(3)$ graphs: infinite $\mathcal{A}^{\ast}$, $\mathcal{D}^{\ast}$ families and their 0-1 parts} \label{fig-SU(3)ADstar}
\end{center}
\end{figure}

\begin{figure}[tbp]
\begin{center}
  \includegraphics[width=110mm]{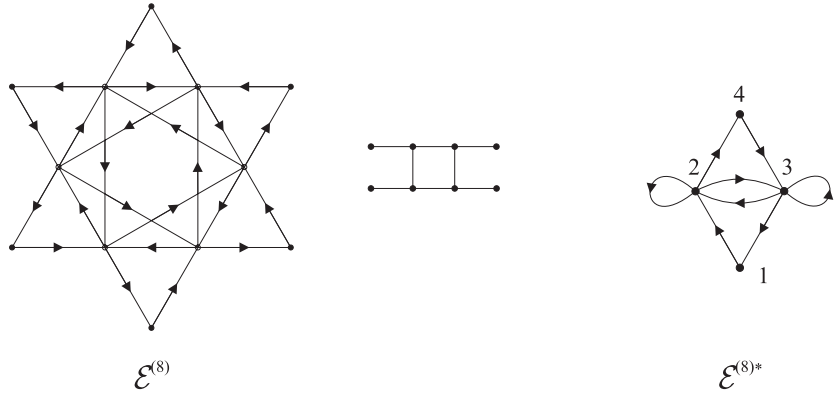}
 \caption{$SU(3)$ graphs: exceptional graphs $\mathcal{E}^{(8)}$, and its 0-1 part, and $\mathcal{E}^{(8)\ast}$} \label{fig-SU(3)E(8)}
\end{center}
\end{figure}

\begin{figure}[tbp]
\begin{center}
  \includegraphics[width=140mm]{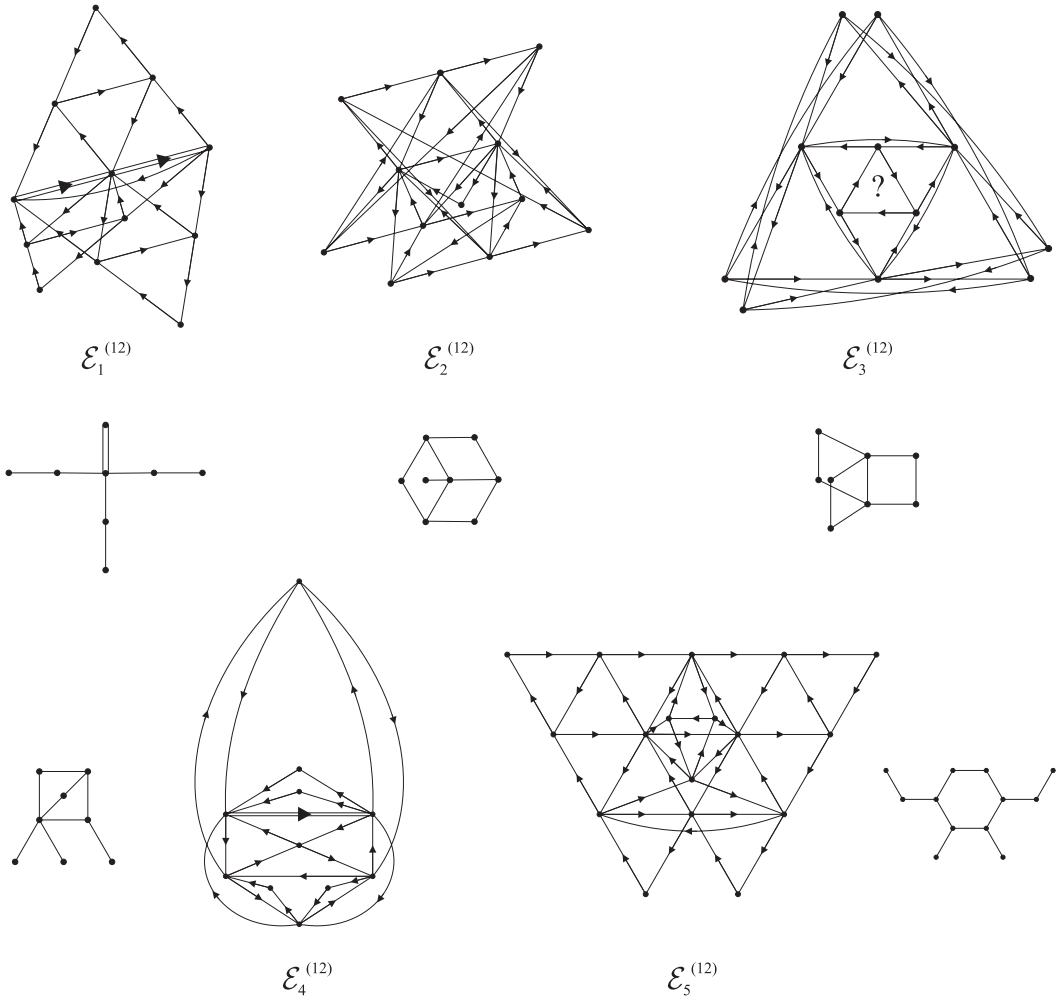}
 \caption{$SU(3)$ graphs: exceptional graphs $\mathcal{E}_i^{(12)}$, $i=1,\ldots,5$, and their 0-1 parts} \label{fig-SU(3)E(12)}
\end{center}
\end{figure}

\begin{figure}[tbp]
\begin{center}
  \includegraphics[width=120mm]{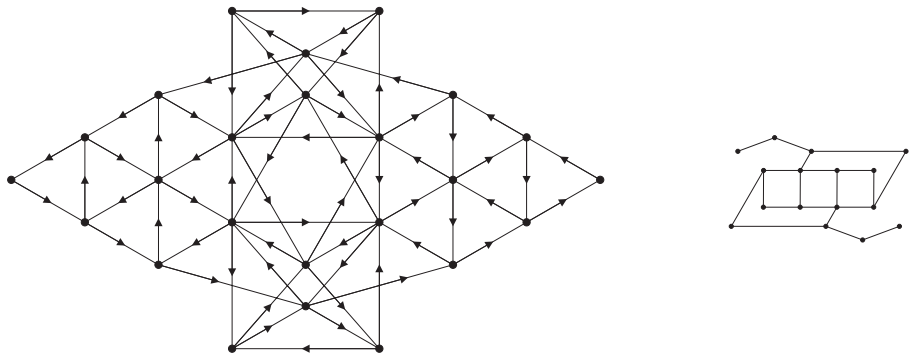}
 \caption{$SU(3)$ graphs: exceptional graph $\mathcal{E}^{(24)}$ and its 0-1 part} \label{fig-SU(3)E(24)}
\end{center}
\end{figure}

Since the $ADE$ graphs can be matched to the affine Dynkin diagrams, di Francesco and Zuber \cite{di_francesco/zuber:1990} were guided to find candidates for classifying graphs for $SU(3)$ modular invariants by first considering the McKay graphs of the finite subgroups of $SU(3)$ to produce a candidate list of $\mathcal{ADE}$ graphs whose spectra described the diagonal part of the modular invariant. They proposed candidates for most of the modular invariants, except for the conjugate invariants $\mathcal{A}^{\ast}$ as they restricted themselves to only look for graphs which are three-colourable.
B\"{o}ckenhauer and Evans \cite{bockenhauer:1999} understood that nimrep graphs for the $SU(3)$ conjugate invariants were not three-colourable. This was also realised simultaneously by Behrend, Pearce, Petkova and Zuber \cite{behrend/pearce/petkova/zuber:2000} and Ocneanu \cite{ocneanu:2002}.
The classification of nimreps is incomplete if one relaxes the condition that the nimrep be compatible with a modular invariant \cite{gannon:2001, graves:2010}.
The complete list of the $\mathcal{ADE}$ graphs are illustrated in Figures \ref{fig-SU(3)AD}-\ref{fig-SU(3)E(24)}.
Where the graph is three-colourable, we have also drawn its 0-1 part, that is the subgraph consisting of the vertices of colour 0 and 1 and the edges connecting these vertices.

There are four infinite families of graphs, $\mathcal{A}^{(m)}$, $m \geq 4$, $\mathcal{D}^{(m)}$, $m \geq 5$, $\mathcal{A}^{(m)\ast}$, $m \geq 4$, and $\mathcal{D}^{(m)\ast}$, $m \geq 6$. There are seven exceptional graphs $\mathcal{E}^{(8)}$, $\mathcal{E}^{(8)\ast}$, $\mathcal{E}_i^{(12)}$, $i=1,2,4,5$, and $\mathcal{E}^{(24)}$.
Here the graphs $\mathcal{D}^{(m)}$ may be obtained from $\mathcal{A}^{(m)}$ by a $\mathbb{Z}_3$-orbifold procedure, and vice-versa.
For $m \neq 3p$, $p \in \mathbb{N}$, the graphs $\mathcal{D}^{(m)}$ are the $SU(3)$ analogues of the tadpole graphs $T_m$ in the $SU(2)$ situation. However, unlike $SU(2)$, these $\mathcal{D}^{(m)}$ are not discarded as $N$-$M$ graphs.
The graphs $\mathcal{D}^{(m)\ast}$, $\mathcal{E}^{(8)\ast}$, $\mathcal{E}_2^{(12)}$ may similarly be obtained from $\mathcal{A}^{(m)\ast}$, $\mathcal{E}^{(8)}$, $\mathcal{E}_1^{(12)}$ respectively by a $\mathbb{Z}_3$-orbifold procedure, and vice-versa. The graphs $\mathcal{E}_1^{(12)}$ and $\mathcal{E}_2^{(12)}$ are isospectral. A third isospectral graph $\mathcal{E}_3^{(12)}$ was proposed by di Francesco and Zuber in \cite{di_francesco/zuber:1990}, however this graph was discarded by Ocneanu \cite{ocneanu:2000ii, ocneanu:2002} as it does not support a cell system. This claim was proven in \cite{coquereaux/isasi/schieber:2010}. In the case where a graph is not three-colourable, the graph corresponding to its 0-1 part could be thought of as the 0-1 part of the unfolded graph, that is, its $\mathbb{Z}_3$-orbifold.
There is a conjugation on the $SU(3)$ $\mathcal{ADE}$ graphs: the conjugation $\tau: {}_N \mathcal{X}_N \rightarrow {}_N \mathcal{X}_N$ on the braided system of endomorphisms of $SU(3)_k$ on a factor $N$, given by the conjugation on the representations of $SU(3)$, induces a conjugation $\tau: {}_N \mathcal{X}_M \rightarrow {}_N \mathcal{X}_M$ such that $G_{\overline{\lambda}} = \tau G_{\lambda} \tau$, where $G_{\lambda} a = \lambda a$ for $\lambda \in {}_N \mathcal{X}_N$, $a \in {}_N \mathcal{X}_M$. The conjugation interchanges the vertices of colours $1 \leftrightarrow 2$.

\subsection{Finite subgroups of $SU(3)$} \label{sect:subgroups-SU(3)}

\begin{figure}[b]
\begin{center}
  \includegraphics[width=145mm]{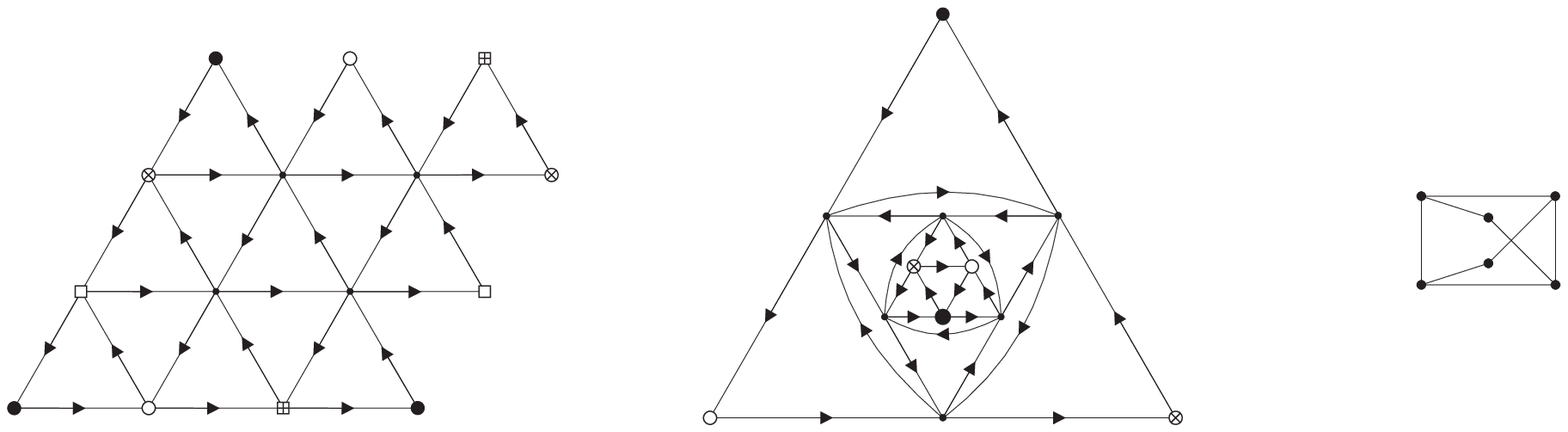}\\
 \caption{Two representations of $\mathbb{Z}_3 \times \mathbb{Z}_3$, and its 0-1 part.}\label{fig:extended-ADstargraph}
\end{center}
\end{figure}

The classification of finite subgroups of $SU(3)$ is due to \cite{miller/blichfeldt/dickson:1916, fairbairn/fulton/klink:1964, bovier/luling/wyler:1981, yau/yu:1993, grimus/ludl:2011}. Clearly any finite subgroup of $U(2)$ is a finite subgroup of $SU(3)$, since we can embed $U(2)$ in $SU(3)$ by sending $g \rightarrow g \oplus \mathrm{det}(g)^{-1} \in SU(3)$, for any $g \subset U(2)$. These subgroups of $SU(3)$ are called type B.
There are three other infinite series of finite groups, called type A, C, D. The groups of type A are the diagonal abelian groups $\mathbb{Z}_m \times \mathbb{Z}_n$, which correspond to an embedding of the two torus $\mathbb{T}^2$ in $SU(3)$ given by
$$(\rho|_{\mathbb{T}^2})(\omega_1,\omega_2) = \left( \begin{array}{ccc} \omega_1 & 0 & 0 \\ 0 & \omega_2^{-1} & 0 \\ 0 & 0 & \omega_1^{-1}\omega_2 \end{array} \right),$$
for $(\omega_1,\omega_2) \in \mathbb{T}^2$. The groups of type C are isomorphic to $(\mathbb{Z}_m \times \mathbb{Z}_n) \rtimes \mathbb{Z}_3$, and those of type D are isomorphic to $(\mathbb{Z}_m \times \mathbb{Z}_n) \rtimes S_3$, where in both cases $n$ is a divisor of $m$ \cite{ludl:2011, grimus/ludl:2011}. The infinite family of groups of type C contains the groups $\Delta (3m^2) \cong (\mathbb{Z}_m \times \mathbb{Z}_m) \rtimes \mathbb{Z}_3$, whilst the infinite family of groups of type D contains the groups $\Delta(6m^2) \cong (\mathbb{Z}_m \times \mathbb{Z}_m) \rtimes S_3$.
The ternary trihedral groups $\Delta (3m^2)$, $\Delta(6m^2)$ were considered in \cite{bovier/luling/wyler:1981, luhn/nasri/ramond:2007, escobar/luhn:2009, ludl:2011, grimus/ludl:2011} and generalize the binary dihedral subgroups of $SU(2)$.
The group $\Delta (3m^2)$ has a presentation generated by the following matrices in $SU(3)$:
$$S_1 = \left( \begin{array}{ccc} 1   &   0     &   0   \\
                                0   & \zeta     &   0   \\
                                0   &   0       & \zeta^{-1}   \end{array} \right),
S_2 = \left( \begin{array}{ccc} \zeta   &   0   &   0   \\
                                0       &  1    &   0   \\
                                0       &   0   & \zeta^{-1}   \end{array} \right)
\quad \textrm{ and } \quad
T = \left( \begin{array}{ccc}   0   &   1   &   0   \\
                                0   &   0   &   1   \\
                                1   &   0   &   0   \end{array} \right),$$
where $\zeta = e^{2 \pi i/m}$.
In this presentation, the action of $\mathbb{Z}_3$ on $\mathbb{Z}_m \times \mathbb{Z}_m \cong \langle S_1, S_2 \rangle$ is given by conjugation by the matrix $T$. Here $\langle g_1, g_2, \ldots, g_s \rangle$ denotes the group generated by the elements $g_1, g_2, \ldots, g_s$.
The group $\Delta (6m^2)$ has a presentation generated by the generators $S_1, S_2$ and $T$ of $\Delta(3n^2)$, and the matrix $Q \in SU(3)$ given by
$$Q = \left( \begin{array}{ccc}   -1  &   0   &   0   \\
                                  0   &   0   &   -1  \\
                                  0   &   -1  &   0   \end{array} \right).$$
In this presentation, the action of $S_3$ on $\mathbb{Z}_m \times \mathbb{Z}_m \cong \langle S_1, S_2 \rangle$ is given by conjugation by the matrices $T, Q$.
This corrects the presentation of these ternary trihedral groups $\Delta (3m^2)$, $\Delta(6m^2)$, and the action of $\mathbb{Z}_3$, $S_3$ respectively, given in \cite{evans/pugh:2010i}.
There are also eight exceptional groups E-L. The generators and relations for these groups are given in \cite{evans/pugh:2010i}. The complete list of finite subgroups of $SU(3)$ is given in Table \ref{Table:ADE-subgroupsSU(3)}.
The notation $\lfloor x \rfloor$ denotes the integer part of $x$. The subgroup $TA_5$ (respectively $TA_6$, $TPSL(2,7)$) is the ternary $A_5$ group (respectively ternary $A_6$ group, ternary $PSL(2,7)$ group), which is the extension of $A_5$ (respectively $A_6$, $PSL(2,7)$) by a cyclic group of order 3.

\begin{center}
\begin{table}[tb]
\begin{tabular}{|c|c|c|c|} \hline
$\mathcal{ADE}$ graph & Type & Subgroup $\Gamma \subset SU(3)$ & $|\Gamma|$ \\
\hline\hline ($ADE$) & - & B: finite subgroups of $SU(2) \subset SU(3)$ & - \\
\hline $\mathcal{A}^{(m)}$ & I & A: $\mathbb{Z}_{m-2} \times \mathbb{Z}_{m-2}$ & $(m-2)^2$ \\
\hline - & - & A: $\mathbb{Z}_m \times \mathbb{Z}_n \;\;$ ($m \neq n \neq 3$) & $mn$ \\
\hline $\mathcal{D}^{(m)} \;$ ($m \equiv 0 \textrm{ mod } 3$) & I & C: $\Delta (3(m-3)^2) = (\mathbb{Z}_{m-3} \times \mathbb{Z}_{m-3}) \rtimes \mathbb{Z}_3$ & $3(m-3)^2$ \\
\hline $\mathcal{D}^{(m)} \;$ ($m \not \equiv 0 \textrm{ mod } 3$) & II & - & - \\
\hline - & - & C: $\Delta (3m^2) = (\mathbb{Z}_m \times \mathbb{Z}_m) \rtimes \mathbb{Z}_3$, & $3m^2$ \\
& & ($n \not \equiv 0 \textrm{ mod } 3$) & \\
\hline - & - & C: $(\mathbb{Z}_m \times \mathbb{Z}_n) \rtimes \mathbb{Z}_3 \;\;$ ($m \neq n)$ & 3mn \\
\hline - & - & D: $\Delta (6m^2) = (\mathbb{Z}_m \times \mathbb{Z}_m) \rtimes S_3$ & $6m^2$ \\
\hline - & - & D: $(\mathbb{Z}_m \times \mathbb{Z}_n) \rtimes S_3 \;\;$ ($m \neq n)$ & 6mn \\
\hline $\mathcal{A}^{(m)\ast}$ & II & - & - \\
\hline $\mathcal{D}^{(m) \ast} \;$ ($m \geq 7$) & II & A: $\mathbb{Z}_{\lfloor (m+1)/2 \rfloor} \times \mathbb{Z}_{3}$ & $3 \lfloor (m+1)/2 \rfloor$ \\
\hline $\mathcal{E}^{(8)}$ & I & E $= \Sigma (36 \times 3) = \Delta (3.3^2) \rtimes \mathbb{Z}_4$ & 108 \\
\hline $\mathcal{E}^{(8)\ast}$ & II & - & - \\
\hline $\mathcal{E}_1^{(12)}$ & I & F $= \Sigma (72 \times 3)$ & 216 \\
\hline $\mathcal{E}_2^{(12)}$ & II & G $= \Sigma (216 \times 3)$ & 648 \\
\hline ($\mathcal{E}_3^{(12)}$) & - & B $\times \mathbb{Z}_3$: $BD_4 \times \mathbb{Z}_3$ & 24 \\
\hline $\mathcal{E}_4^{(12)}$ & (II) & L $= \Sigma (360 \times 3) \cong TA_6$ & 1080 \\
\hline $\mathcal{E}_5^{(12)}$ & II & K $\cong TPSL(2,7)$ & 504 \\
\hline $\mathcal{E}^{(24)}$ & I & - & - \\
\hline - & - & H $= \Sigma (60) \cong A_5$ & 60 \\
\hline - & - & I $= \Sigma (168) \cong PSL(2,7)$ & 168 \\
\hline - & - & J $\cong TA_5$ & 180 \\
\hline
\end{tabular} \\
\caption{Relationship between $\mathcal{ADE}$ graphs and subgroups $\Gamma$ of $SU(3)$.} \label{Table:ADE-subgroupsSU(3)}
\end{table}
\end{center}

\begin{figure}[tbp]
\begin{center}
  \includegraphics[width=90mm]{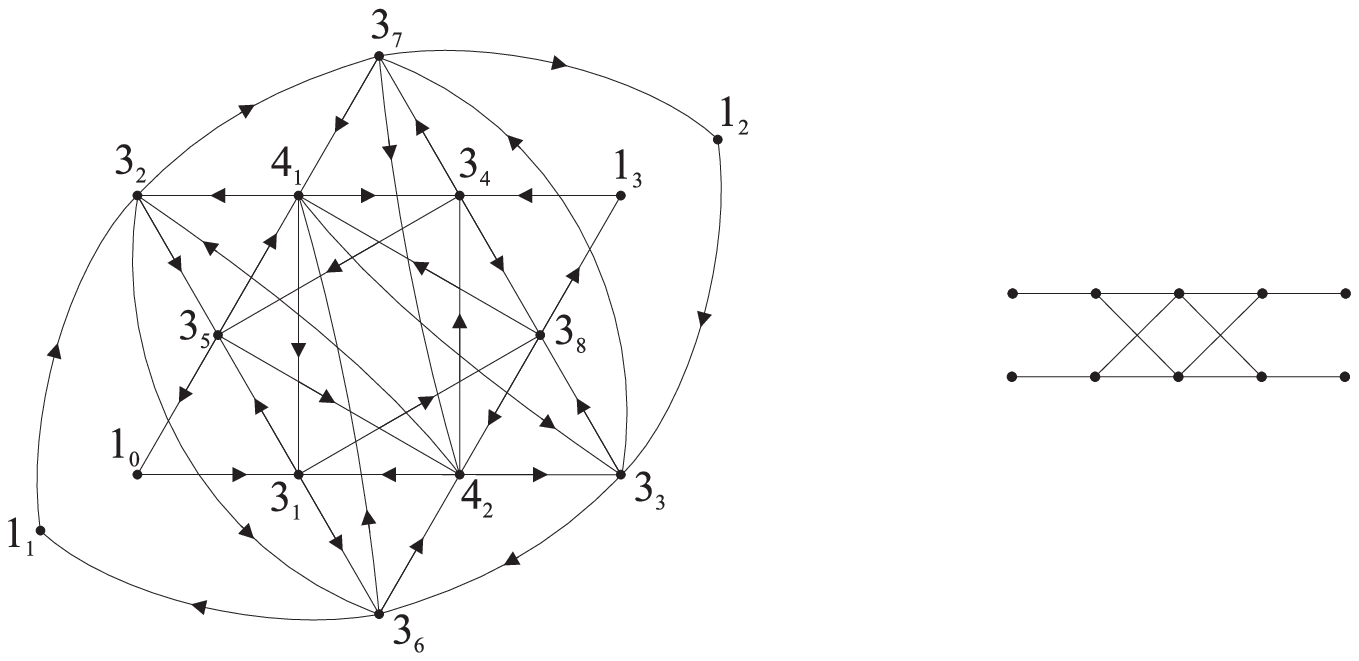}\\
 \caption{E $= \sum (36 \times 3)$ and its 0-1 part}\label{fig:extended-E(8)graph}
\end{center}
\end{figure}

\begin{figure}[tbp]
\begin{center}
  \includegraphics[width=140mm]{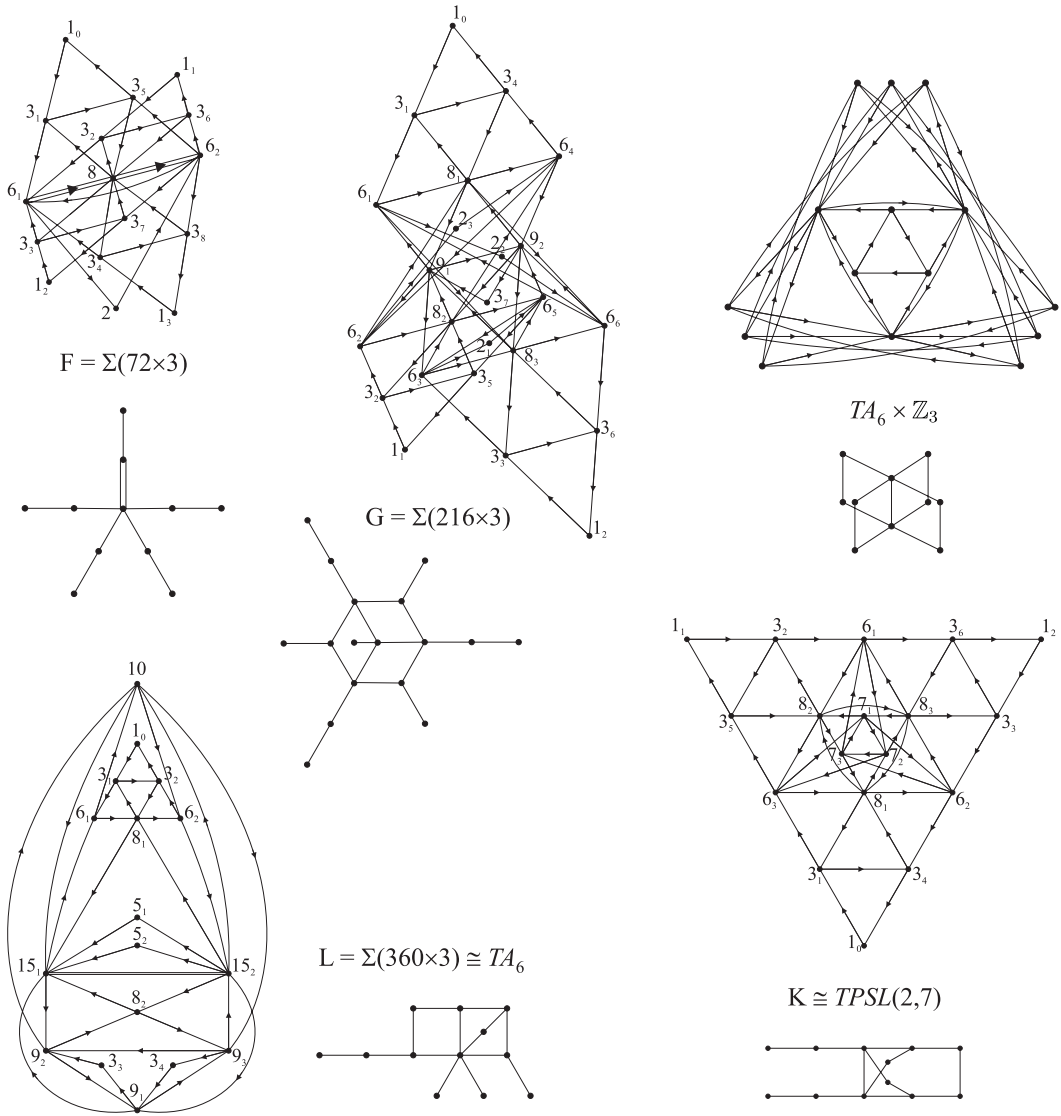}\\
 \caption{F $= \Sigma (72 \times 3)$, G $= \Sigma (216 \times 3)$, $BD_4 \times \mathbb{Z}_3$, L $= \Sigma (360 \times 3) \cong TA_6$, K $\cong TPSL(2,7)$, and their 0-1 parts}\label{fig:extended-E(12)graph}
\end{center}
\end{figure}

\begin{figure}[tbp]
\begin{center}
  \includegraphics[width=125mm]{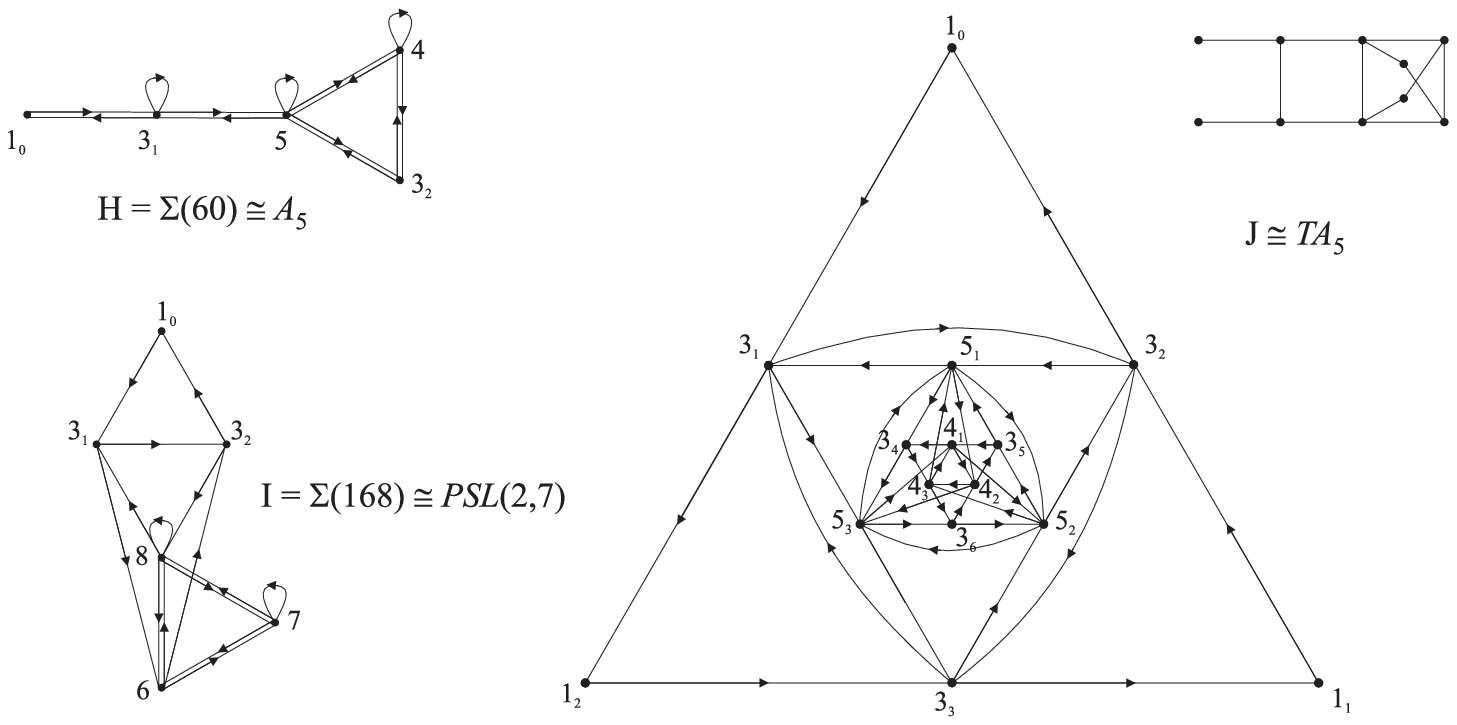}\\
 \caption{H $= \Sigma (60) \cong A_5$, I $= \Sigma (168) \cong PSL(2,7)$, and J $\cong TA_5$ and its 0-1 part}\label{fig:extended-HIJgraph}
\end{center}
\end{figure}

\vspace{-12mm}
Unlike in the case of $SU(2)$, these $SU(3)$ $\mathcal{ADE}$ graphs cannot be matched exactly to (the McKay graphs of the representation theory of) the finite subgroups of $SU(3)$. There is however a correspondence between some of the $SU(3)$ $\mathcal{ADE}$ graphs and McKay graphs for the finite subgroups of $SU(3)$, as shown in Table \ref{Table:ADE-subgroupsSU(3)}, but there are $SU(3)$ $\mathcal{ADE}$ graphs which do not have a corresponding finite subgroup of $SU(3)$ and vice-versa.
The $SU(3)$ $\mathcal{ADE}$ graph can be obtained from its corresponding McKay graph $\mathcal{G}_{\Gamma}$ by now removing more than one vertex, and all the edges that start or end at those vertices, as well as possibly some other edges, as was noted in \cite{di_francesco/zuber:1990}. To obtain the graph $\mathcal{E}_5^{(12)}$ from the McKay graph for the subgroup $\mathrm{K} \cong TPSL(2,7)$ an extra edge must also be added.
The McKay graphs are illustrated in Figures \ref{fig:extended-ADstargraph}-\ref{fig:extended-HIJgraph}. In Figure \ref{fig:extended-ADstargraph} we have illustrated two representations of $\mathbb{Z}_3 \times \mathbb{Z}_3$, where vertices with the same symbol are identified: the representation on the left shows the connection between $\mathbb{Z}_{m-2} \times \mathbb{Z}_{m-2}$ and the graph $\mathcal{A}^{(m)}$, whilst the other representation shows the connection between $\mathbb{Z}_{p} \times \mathbb{Z}_{3}$ and the graph $\mathcal{D}^{(2p+2)\ast}$. The 0-1 part of $\mathbb{Z}_{m} \times \mathbb{Z}_{n}$ only makes sense when both $m, n \equiv 0 \textrm{ mod } 3$.

\subsection{Spectrum of nimrep index values}

\begin{figure}[p]
\begin{center}
  \includegraphics[width=130mm]{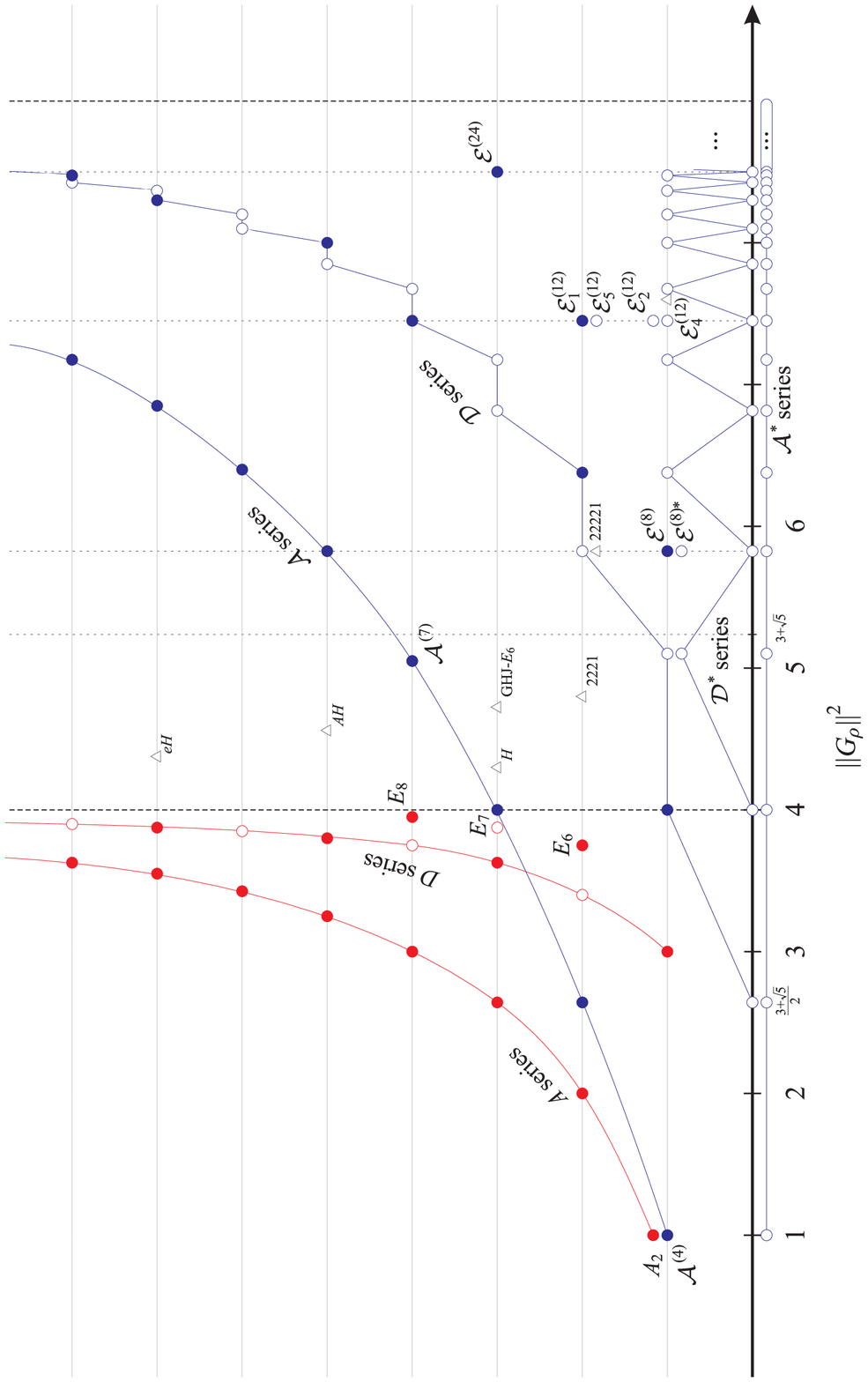}\\
  \caption{Plotting the $SU(n)$-supertransitivity and the norm squared of the \hbox{$SU(n)$ $N$-$M$ nimrep graphs $G_{\rho}$, $n=2,3$.}} \label{fig:plot_M-N_graphs}
\end{center}
\end{figure}

Jones defined a notion of supertransitivity for a subfactor, or equivalently, for its principal graph \cite{jones:2001, grossman/jones:2007}.
For $n=2,3$, we define the \emph{$SU(n)$-supertransitivity} of an $SU(n)$ $N$-$M$ graph $\mathcal{G}$ to be $k \in \mathbb{N}$, where $k$ is the largest integer such that $\bigoplus_v \ast \cdot \mathrm{End} \left( (\mathbb{C}\mathcal{G})_{k} \right) \cdot v = \bigoplus_{v'} \ast \cdot \mathrm{End} \left( (\mathbb{C}\mathcal{A})_{k} \right) \cdot v'$, where the direct sum on the left hand side is over all vertices $v$ of $\mathcal{G}$ and the direct sum on the right hand side is over all vertices $v'$ of $\mathcal{A}$, where $\mathcal{A}$ is the $SU(n)$ $A$-graph with the same graph norm as $\mathcal{G}$. Note that $SU(n)$-supertransitivity reduces to the usual supertransitivity of subfactor theory when $n=2$, in the case where $\mathcal{G}$ is the principal graph of a subfactor.

In Figure \ref{fig:plot_M-N_graphs} we have plotted the $SU(n)$-supertransitivity and norm of the $SU(n)$ $N$-$M$ graphs $G_{\rho}$, for $n=2,3$.
This figure was inspired by the figure of Scott Morrison in \cite{morrison:2011}.
Here a solid, hollow circle respectively denotes that the inclusion $N \subset M$ is of type I, II respectively.
Where more than one graph has the same norm and $SU(n)$-supertransitivity, one graph is stacked vertically on top of the other.
The norm of the $SU(2)$ $ADE$ graphs is of the form $[2]=2\cos(\pi/m)$, whilst the norm of the $SU(3)$ $\mathcal{ADE}$ graphs is of the form $[3]=1+2\cos(2\pi/m)$.
When $\mathcal{G}$ is one of the $SU(3)$ $N$-$M$ nimrep graphs $\mathcal{A}^{(m)}$ and $\mathcal{D}^{(3k)}$, the principal graph is given by the 0-1 part of $\mathcal{G}$ \cite[Theorems 3.3, 5.8 \& Corollary 3.4]{evans/kawahigashi:1994}. It is conjectured \cite[Section 5.7]{evans/pugh:2009ii} that the same is true for the other $SU(3)$ $N$-$M$ graphs which are realised by type I inclusions.
There are 4 $SU(3)$ $N$-$M$ graphs which are realised by type I inclusions and which have norm squared less than or equal to 4. These are the graphs $\mathcal{A}^{(m)}$, $m=4,5,6$, and $\mathcal{D}^{(6)}$. The 0-1 part of $\mathcal{A}^{(4)}$ is the Dynkin diagram $A_3$, the 0-1 part of $\mathcal{A}^{(5)}$ is $A_4$, the 0-1 part of $\mathcal{A}^{(6)}$ is the affine Dynkin diagram $E_6^{(1)}$, and the 0-1 part of $\mathcal{D}^{(6)}$ is the affine Dynkin diagram $A_1^{(1)}$.

The hollow triangles indicate certain other graphs $\mathcal{G}$ which are principal graphs for subfactors with index $||\mathcal{G}||^2$, where $SU(2)$-supertransitivity has been used for the vertical axis.
The graphs $H$, $AH$, $eH$, $GHJ\textrm{-}E_6$, $2221$ are all the principal graphs for subfactors with index between 4 and 5. The first three subfactors are the Haagerup $H$ with index $(5+\sqrt{13})/2$ \cite{haagerup:1994}, the Asaeda-Haagerup $AH$ with index $(5+\sqrt{17})/2$ \cite{asaeda/haagerup:1999} and the extended Haagerup $eH$ with index given by the largest root of $x^3 - 8 x^2 + 17 x - 5$ \cite{bigelow/morrison/peters/snyder:2009}. The next is the Goodman-de la Harpe-Jones subfactor for the inclusion of $SU(2)$ in $E_6$ with index $3+\sqrt{3}$ \cite{goodman/de_la_harpe/jones:1989}, and the final subfactor can be constructed from the conformal inclusion $(G_2)_3 \subset (E_6)_1$, with index $(5+\sqrt{21})/2$ \cite{izumi:2001}.
There are two other hollow triangles, one which indicates the graph $22221$, which has the same norm $1+\sqrt{2}$ as the $SU(3)$ graphs with Coxeter number 8. The other is unlabelled, and indicates the subfactor of \cite{grossman/snyder:2011} with index $4+\sqrt{13}$, which is a subfactor corresponding to a simple object in the category of Morita equivalences between unitary fusion categories $H_1$ and $H_3$, where $H_1$ and $H_2$ are unitary fusion categories coming from the Haagerup subfactor, and $H_3$ is a unitary fusion category with the Grothendieck ring as $H_2$ but which is not isomorphic to $H_2$ as a category.
We have marked the value of $3+\sqrt{5}$, which is twice the golden mean, on the graph. This value is the smallest index value which can be obtained as a non-trivial multiple of other index values, apart from the index value $4=2^2$.
Note that beyond 4, the square of the norm of all the $SU(3)$ $\mathcal{ADE}$ graphs has non-integer value (the norm squared for the graphs with Coxeter number 15 is 7.992 to 3 decimal places).
The Haagerup subfactor is the first of a series of subfactors \cite{izumi:2001, evans/gannon:2011} and the modular data of their doubles have been studied in \cite{evans/pinto:2006, evans/pinto:2011, evans/gannon:2011}

\section{Spectral Measures for nimreps} \label{sect:spectral_measures}

We wish to compute invariants for the module categories or braided subfactors associated to $SU(3)$ modular invariants. We will now consider spectral measures for the classifying graphs $G_{\rho}$ and the McKay graphs of subgroups of $SU(3)$. Such measures were computed by Banica and Bisch \cite{banica/bisch:2007} for the $ADE$ graphs and the McKay graphs of subgroups of $SU(2)$.

Suppose $A$ is a unital $C^{\ast}$-algebra with state $\varphi$.
If $b \in A$ is a normal operator then there exists a compactly supported probability measure $\mu_b$ on the spectrum $\sigma(b) \subset \mathbb{C}$ of $b$, uniquely determined by its moments
$\varphi(b^m b^{\ast n}) = \int_{\sigma(b)} z^m \overline{z}^n \mathrm{d}\mu_b (z)$,
for non-negative integers $m$, $n$.

We computed in \cite{evans/pugh:2009v, evans/pugh:2010i} such spectral measures and generating series when $b$ is the normal operator $\Delta = G_{\rho}$ acting on the Hilbert space of square summable functions on the graph, for the nimreps $G$ described above, i.e. $G_{\rho}$ is the adjacency matrix of the $ADE$ and $SU(3)$ $\mathcal{ADE}$ graphs, and the McKay graphs of subgroups of $SU(2)$ and $SU(3)$. We computed the spectral measure for the vacuum, i.e. the distinguished vertex of the graph which has lowest Perron-Frobenius weight. However the spectral measures for the other vertices of the graph could also be computed by the same methods. For the $SU(2)$ and $SU(3)$ graphs, the spectral measures distill onto very special subsets of the semicircle/circle ($SU(2)$) and discoid/torus ($SU(3)$), and the theory of nimreps allowed us to compute these measures precisely. Our methods gave an alternative approach to deriving the results of Banica and Bisch \cite{banica/bisch:2007} for $ADE$ graphs and subgroups of $SU(2)$, and explained the connection between their results for affine $ADE$ graphs and the Kostant polynomials.

In particular, for $SU(2)$, we can understand the spectral measures for the torus $T$ and $SU(2)$ as follows.
If $w_Z$ and $w_N$ are the self adjoint operators arising from the McKay graph of the fusion rules of the representation theory of $T$ and $SU(2)$, then the spectral measures in the vacuum state can be describe in terms of semicircular law, on the interval $[-2,2]$ which is the spectrum of either
as the image of the map $z \in T \rightarrow z + z^{-1}$ \cite[Sections 2 \& 3.1]{evans/pugh:2009v}:
$$\mathrm{dim}\left( \left(\otimes^k M_2 \right)^{\mathbb{T}} \right) \;\; = \;\; \varphi(w_Z^{2k}) \;\; = \;\; \frac{1}{\pi} \int_{-2}^2 x^{2k} \frac{1}{\sqrt{4-x^2}} \; \mathrm{d}x \, ,$$
$$\mathrm{dim}\left( \left(\otimes^k M_2 \right)^{SU(2)} \right) \;\; = \;\; \varphi(w_N^{2k}) \;\; = \;\; \frac{1}{2\pi} \int_{-2}^2 x^{2k} \sqrt{4-x^2} \; \mathrm{d}x \, .$$
The fusion matrix for $SO(3)$ is just $\mathbf{1} + \Delta$, where $\Delta$ is the fusion matrix for $SU(2)$, and thus is equal to the infinite $SU(3)$ $\mathcal{ADE}$ graph $\mathcal{A}^{(\infty)\ast}$. Thus the spectral measure $\mu$ in the vacuum state (over $[-1,3]$) for $SO(3)$ has semicircle distribution with mean 1 and variance 1, i.e. $\mathrm{d}\mu(x) = \sqrt{4 - (x-1)^2} \mathrm{d}x$ \cite[Section 7.3]{evans/pugh:2009v}.

\begin{figure}[tb]
\begin{center}
  \includegraphics[width=70mm]{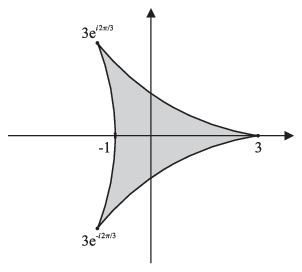}\\
  \caption{The (3-cusp) discoid $\mathfrak{D}$, the union of the deltoid and its interior.} \label{fig:hypocycloid-S}
\end{center}
\end{figure}

The spectral weight $\sqrt{4 - (x-1)^2}$ for $SU(2)$ arises from the Jacobian of a change of variable between the interval $[-2,2]$ and the circle.
Then for $\mathbb{T}^2$ and $SU(3)$, the 3-cusp discoid $\mathfrak{D}$ in the complex plane, illustrated in Figure \ref{fig:hypocycloid-S}, is the image of the two-torus under the map $\Phi: (\omega_1, \omega_2) \mapsto \omega_1 + \omega_2^{-1} + \omega_1^{-1} \omega_2$, which is the spectrum of the corresponding normal operators $v_Z$ and $v_N$ on the Hilbert spaces of the fusion graphs of $\mathbb{T}^2$ and $SU(3)$ respectively. The corresponding spectral measures are then described by a corresponding Jacobian or discriminant $\sqrt{27 - 18z\overline{z} + 4z^3 + 4\overline{z}^3 -z^2 \overline{z}^2}$ as \cite[Theorems 3 \& 5]{evans/pugh:2009v}:
$$\mathrm{dim}\left( \left(\otimes^k M_3 \right)^{\mathbb{T}^2} \right) \;\; = \;\; \varphi(|v_Z|^{2k}) \;\; = \;\; \frac{3}{\pi^2} \int_{\mathfrak{D}} |z|^{2k} \frac{1}{\sqrt{27 - 18z\overline{z} + 4z^3 + 4\overline{z}^3 -z^2 \overline{z}^2}} \; \mathrm{d}z \, ,$$
$$\mathrm{dim}\left( \left(\otimes^k M_3 \right)^{SU(3)} \right) \;\; = \;\; \varphi(|v_N|^{2k}) \;\; = \;\; \frac{1}{2\pi^2} \int_{\mathfrak{D}} |z|^{2k} \sqrt{27 - 18z\overline{z} + 4z^3 + 4\overline{z}^3 -z^2 \overline{z}^2} \; \mathrm{d}z \, ,$$
where $\mathrm{d}z:=\mathrm{d}\,\mathrm{Re}z \; \mathrm{d}\,\mathrm{Im}z$ denotes the Lebesgue measure on $\mathbb{C}$.

For vertices $\nu$ of $\mathcal{A}^{(m)}$ (generalized) Chebyshev polynomials $S_{\nu}(x,y)$ in two variables are defined by $S_{(0,0)}(x,y) = 1$, and $x S_{\nu}(x,y) = \sum_{\mu} \Delta_{\mathcal{A}}(\nu, \mu) S_{\mu}(x,y)$ and $y S_{\nu}(x,y) = \sum_{\mu} \Delta_{\mathcal{A}}^T(\nu, \mu) S_{\mu}(x,y)$.
For concrete values of the first few $S_{\mu}(x,y)$ see \cite[p. 610]{evans/kawahigashi:1998}.
Koornwinder \cite{koornwinder:1975} proved that the measure for $v_N$ above is the measure required to make these polynomials $S_{\mu}(z,\overline{z})$ orthogonal, i.e.
$$\frac{1}{2 \pi^2} \int_{\mathbb{T}^2} S_{\mu}(z,\overline{z}) \overline{S_{\nu}(z,\overline{z})} \sqrt{27 - 18z\overline{z} + 4z^3 + 4\overline{z}^3 -z^2 \overline{z}^2} \; \mathrm{d}z = \delta_{\mu,\nu}.$$
Gepner \cite{gepner:1991} extended this result to orthogonal measures for the analogues of the Chebyshev polynomials for $SU(n)$, for each $n \geq 2$.

Related constructions of spectral measures and spectra can be found in \cite{cartwright/mlotkowski:1994, cartwright/mlotkowski/steger:1994} for averaging operators on triangle buildings, that is, on simplicial complexes consisting of vertices, edges and triangles.

\begin{figure}[bt]
\begin{center}
  \includegraphics[width=55mm]{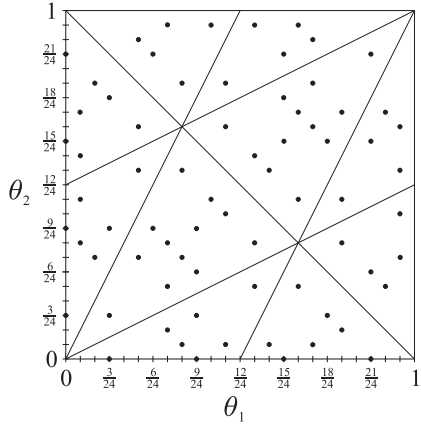}\\
 \caption{The points $(\theta_1,\theta_2)$ for $\mathcal{E}^{(8)}$.} \label{fig:poly-11}
\end{center}
\end{figure}

As an example of the spectral measure of one of the finite $\mathcal{ADE}$ graphs, the spectral measure for $\mathcal{E}^{(8)}$ over $\mathbb{T}^2$ is
\begin{equation}
\mathrm{d}\varepsilon = \frac{2-\sqrt{2}}{8} \; \mathrm{d}^{((8))} + \frac{2+\sqrt{2}}{8} \; \mathrm{d}^{((8/3))} + \frac{1}{2} \; \mathrm{d}^{(24/5,1/12)},
\end{equation}
where $\mathrm{d}^{((n))}$ is the uniform measure on the $S_3$-orbit of the points $(\tau, \tau)$, $(\overline{\omega} \, \overline{\tau}, \omega)$, $(\omega, \overline{\omega} \, \overline{\tau})$, for $n \in \mathbb{Q}$, and $\mathrm{d}^{(n,k)}$ is the uniform measure on the $S_3$-orbit of the points $(\tau \, e^{2 \pi i k}, \tau)$, $(\tau, \tau \, e^{2 \pi i k})$, $(\overline{\omega} \, \overline{\tau}, \omega \, e^{2 \pi i k})$, $(\omega \, e^{2 \pi i k}, \overline{\omega} \, \overline{\tau})$, $(\overline{\omega} \, \overline{\tau} \, e^{-2 \pi i k}, \omega \, e^{-2 \pi i k})$, $(\omega \, e^{-2 \pi i k}, \overline{\omega} \, \overline{\tau} \, e^{-2 \pi i k})$, for $n, k \in \mathbb{Q}$, $n > 2$, $0 \leq k \leq 1/n$.
It is a discrete measure over $(e^{2 \pi i \theta_1}, e^{2 \pi i \theta_2}) \in \mathbb{T}^2$ for the points $(\theta_1,\theta_2)$ illustrated in Figure \ref{fig:poly-11}.

The theory of nimreps carries over to the case of finite subgroups of $SU(3)$, since there is an $S$-matrix which simultaneously diagonalizes the representations of $\Gamma \subset SU(3)$ \cite{kawai:1989}.
The $m,n^{\mathrm{th}}$ moment $\varsigma_{m,n}$ is given by
$$\varsigma_{m,n} \; = \; \int_{\mathfrak{D}} z^m \overline{z}^n \mathrm{d}\mu(z) \; = \; \sum_{j=1}^s \frac{|\Gamma_j|}{|\Gamma|} \chi_{\rho} (\Gamma_j)^m \overline{\chi_{\rho} (\Gamma_j)}^n.$$
In order to determine the spectral measure over $\mathbb{T}^2$ for the McKay graph $\mathcal{G}_{\Gamma}$, we must first determine \emph{right inverse} maps $\Phi^{-1}:\mathfrak{D} \rightarrow \mathbb{T}^2$ such that $\Phi \circ \Phi^{-1} = \mathrm{id}$, which amounts to finding solutions $\omega_1$ to the cubic equation
\begin{equation} \label{eqn:cubic_w1}
\omega_1^3 - z \omega_1^2 + \overline{z} \omega_1 - 1 = 0.
\end{equation}
Solving (\ref{eqn:cubic_w1}) we obtain solutions $\omega_1 = \omega^{(k)}$, $k = 0,1,2$, given by
$$\omega^{(k)} = (z + 2^{-1/3} \epsilon_k P + 2^{1/3} \overline{\epsilon_k} (z^2-3\overline{z}) P^{-1})/3,$$
where $\epsilon_k = e^{2 \pi i k /3}$, $2^{1/3}$ takes a real value, and $P$ is the cube root $P = (27 - 9z\overline{z} + 2z^3 + 3 \sqrt{3} \sqrt{27 - 18z\overline{z} + 4z^3 + 4\overline{z}^3 - z^2\overline{z}^2})^{1/3}$ such that $P \in \{ r e^{i \theta} | \; 0 \leq \theta < 2 \pi/3 \}$. We notice that the Jacobian $J$ appears in the expression for $P$ as the discriminant of the cubic equation (\ref{eqn:cubic_w1}).
Right inverse maps $\Phi^{-1}_{k,l}:\mathfrak{D} \rightarrow \mathbb{T}^2$ are given by
$\Phi^{-1}_{k,l}(z) = (\omega^{(k)},\overline{\omega^{(l)}})$, for $z \in \mathfrak{D}$, where $k, l \in \{ 0,1,2 \}$ such that $k \neq l$.
The spectral measure of $\Gamma$ (over $\mathbb{T}^2$) is taken as the average over these six $\Phi^{-1}_{k,l}(z)$, which are the $S_3$-orbit of $\Phi^{-1}_{0,1}(z)$:
$$\varsigma_{m,n} \; = \; \frac{1}{6} \sum_{j=1}^s \sum_{\stackrel{k,l \in \{ 0,1,2 \}:}{\scriptscriptstyle{k \neq l}}} \frac{|\Gamma_j|}{|\Gamma|} (\omega^{(k,j)} + \overline{\omega^{(l,j)}} + \overline{\omega^{(k,j)}}\omega^{(l,j)})^m (\overline{\omega^{(k,j)}} + \omega^{(l,j)} + \omega^{(k,j)}\overline{\omega^{(l,j)}})^n,$$
where $\omega^{(p,j)}$, $j=1,\ldots,s$, are given by $\Phi_{k,l}^{-1}(\chi_{\rho}(\Gamma_j)) = (\omega^{(k,j)},\overline{\omega^{(l,j)}})$.

\begin{figure}[tb]
\begin{center}
  \includegraphics[width=55mm]{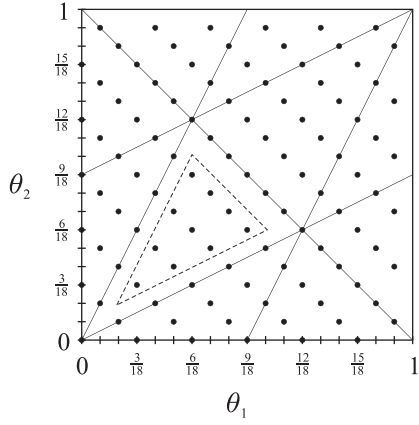}\\
 \caption{The points $(\theta_1,\theta_2)$ such that $(e^{2 \pi i \theta_1}, e^{2 \pi i \theta_2}) \in D_6$.} \label{fig:D6}
\end{center}
\end{figure}

As an example, consider the subgroup $\mathrm{E} = \Sigma(36 \times 3)$, the semidirect product of the ternary trihedral group $\Delta (3.3^2)$ with $\mathbb{Z}_4$, which is the subgroup corresponding to the graph $\mathcal{E}^{(8)}$.
The spectral measure over $\mathbb{T}^2$ for the group E is
\begin{equation}
\mathrm{d}\varepsilon = \frac{1}{48\pi^4} J^2 \, \mathrm{d}^{(4)} + \frac{1}{108\pi^4} J^2 \, \mathrm{d}^{(3)} + \frac{1}{3} \mathrm{d}^{(2)} - \frac{1}{18} \mathrm{d}^{(1)},
\end{equation}
where $\mathrm{d}^{(m)}$ is the spectral measure for $\mathcal{A}^{(m)}$, which is uniform measure over the points in $D_m = \{ (e^{2 \pi i q_1/3m}, e^{2 \pi i q_2/3m}) \in \mathbb{T}^2 | \; q_1,q_2 = 0, 1, \ldots, 3m-1; q_1 + q_2 \equiv 0 \textrm{ mod } 3 \}$, illustrated in Figure \ref{fig:D6} for $m=6$. Notice that the points in the interior of the fundamental domain $C$ (those enclosed by the dashed line) correspond to the vertices of the graph $\mathcal{A}^{(6)}$.

\section{$SU(3)$ module categories and the $A_2$-planar algebra} \label{sect:diagrammatic}

In this section we review some pictorial descriptions associated to the $SU(3)$ $\mathcal{ADE}$ graphs. These include a diagrammatic presentation of the $A_2$-Temperley-Lieb algebra in Section \ref{sect:A2TLalg}, the $A_2$-Temperley-Lieb category in Section \ref{sect:A2TLcat}, and the $A_2$-planar algebras of \cite{evans/pugh:2009iii} in Section \ref{sect:A2PA}. In Section \ref{sect:module_categories} we review the construction of the path algebra for an $\mathcal{ADE}$ graph $\mathcal{G}$ via a monoidal functor coming from a module category with classifying graph $\mathcal{G}$. This is related to the $A_2$-graph planar algebra construction in Section \ref{sect:A2PA}, and a partial decomposition of the $A_2$-graph planar algebras for the $\mathcal{ADE}$ graphs into irreducible $A_2$-planar modules is given in Section \ref{sect:A2PAmodules}.

\subsection{The $A_2$-Temperley-Lieb algebra} \label{sect:A2TLalg}

We begin with a diagrammatic presentation of the $A_2$-Temperley-Lieb algebra given in \cite{evans/pugh:2009iii} (c.f. also \cite{evans/pugh:2010ii}), using the $A_2$ spider of Kuperberg \cite{kuperberg:1996}.
For $SU(3)$, the vanishing of the $q$-antisymmetrizer (\ref{SU(N)q condition}) yields the relation
\begin{equation} \label{eqn:SU(3)q_condition}
\left( U_i - U_{i+2} U_{i+1} U_i + U_{i+1} \right) \left( U_{i+1} U_{i+2} U_{i+1} - U_{i+1} \right) = 0.
\end{equation}
The $A_2$-Temperley-Lieb algebra is the algebra generated by a family $\{ U_n \}$ of self-adjoint operators which satisfy the Hecke relations H1-H3 and the extra condition (\ref{eqn:SU(3)q_condition}).

\begin{figure}[bt]
\begin{center}
  \includegraphics[width=40mm]{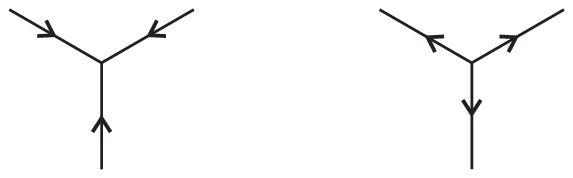}\\
 \caption{$A_2$ webs}\label{fig:A2-webs}
\end{center}
\end{figure}

In \cite{kuperberg:1996}, Kuperberg defined the notion of a spider, which is an axiomatization of the representation theory of groups and other group-like objects. The invariant spaces have bases given by certain planar graphs, called webs. In \cite{kuperberg:1996} certain spiders were defined in terms of generators and relations, isomorphic to the representation theories of rank two Lie algebras and the quantum deformations of these representation theories, which generalized a well-known construction for $A_1 = \textrm{su}(2)$ by Kauffman \cite{kauffman:1987}.
For the $A_2 = \textrm{su}(3)$ case, the $A_2$ webs are illustrated in Figure \ref{fig:A2-webs}, and the relations are given by K1-K3 below.
A complete list of relations for a diagrammatic formalism for $A_3$ has been conjectured in \cite{kim:2003}, and for general $A_l$, $l \geq 2$, in \cite{morrison:2007}.

A string $s$ is a sequence of signs $+$, $-$.
For two (possible empty) strings $s_1, s_2$, an $A_2$-$s_1,s_2$-tangle $T$ is a tangle on an rectangle with strings $s_1$, $s_2$ along the top, bottom edges respectively, generated by the $A_2$ webs such that every free end of $T$ is attached to a vertex along the top or bottom of the rectangle in a way that respects the orientation of the strings, every vertex has a string attached to it, and the tangle contains no elliptic faces. We call a vertex a source vertex if the string attached to it has orientation away from the vertex. Similarly, a sink vertex will be a vertex where the string attached has orientation towards the vertex. Along the top edge the points $+$ are source vertices and $-$ are sink vertices, while along the bottom edge the roles are reversed.
We define the vector space $\mathcal{V}^{A_2}_{s_1,s_2}$ to be the free vector space over $\mathbb{C}$ with basis all $A_2$-$s_1,s_2$-tangles.

We define $V^{A_2}_{s_1,s_2}$ to be the quotient of $\mathcal{V}^{A_2}_{s_1,s_2}$ by the Kuperberg ideal generated by the Kuperberg relations K1-K3 \cite{kuperberg:1996} below. That is, composition in $V^{A_2}_{s_1,s_2}$ is defined as follows. The composition $RS \in V^{A_2}_{s_1,s_3}$ of an $A_2$-$s_1,s_2$-tangle $R$ and an $A_2$-$s_2,s_3$-tangle $S$ is given by gluing $S$ vertically below $R$ such that the vertices at the bottom of $R$ and the top of $S$ coincide, removing these vertices, and isotoping the glued strings if necessary to make them smooth. Any closed loops which may appear are removed, contributing a factor of $\alpha$, as in relation K1 below. If any elliptic faces appear, they are removed using relations K2, K3 below. The composition is associative and is extended linearly to elements in $V^{A_2}_{s_1,s_2}$.
\begin{center}
\begin{minipage}[b]{11.5cm}
 \begin{minipage}[t]{3cm}
  \parbox[t]{2cm}{\begin{eqnarray*}\textrm{K1:}\end{eqnarray*}}
 \end{minipage}
 \begin{minipage}[t]{5.5cm}
  \begin{center}
  \mbox{} \\
 \includegraphics[width=20mm]{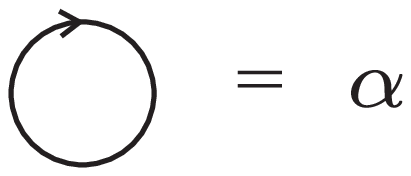}
  \end{center}
 \end{minipage}
 \begin{minipage}[t]{2cm}
  \mbox{} \\
  \parbox[t]{1cm}{}
 \end{minipage}
\end{minipage}
\begin{minipage}[b]{11.5cm}
 \begin{minipage}[t]{3cm}
  \parbox[t]{2cm}{\begin{eqnarray*}\textrm{K2:}\end{eqnarray*}}
 \end{minipage}
 \begin{minipage}[t]{5.5cm}
  \begin{center}
  \mbox{} \\
 \includegraphics[width=23mm]{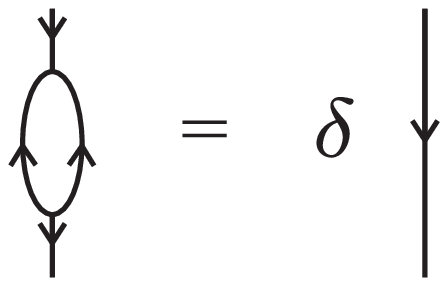}
  \end{center}
 \end{minipage}
 \begin{minipage}[t]{2cm}
  \mbox{} \\
  \parbox[t]{1cm}{}
 \end{minipage}
\end{minipage}
\begin{minipage}[b]{11.5cm}
 \begin{minipage}[t]{3cm}
  \parbox[t]{2cm}{\begin{eqnarray*}\textrm{K3:}\end{eqnarray*}}
 \end{minipage}
 \begin{minipage}[t]{5.5cm}
  \begin{center}
  \mbox{} \\
 \includegraphics[width=55mm]{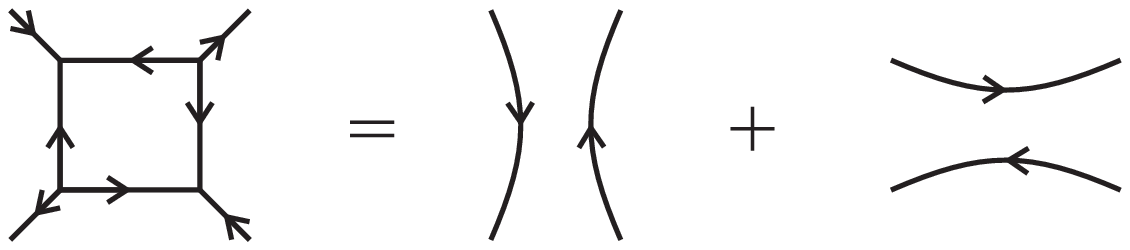}
  \end{center}
 \end{minipage}
 \begin{minipage}[t]{2cm}
  \mbox{} \\
  \parbox[t]{1cm}{}
 \end{minipage}
\end{minipage}
\end{center}
The local picture $\,$ \includegraphics[width=10mm]{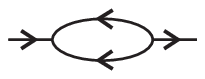} $\,$ is called a \emph{digon}, and $\,$ \includegraphics[width=8mm]{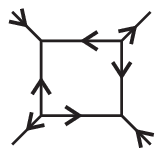} $\,$ an \emph{embedded square}.
There is a braiding on $\mathcal{V}^{A_2}_{s_1,s_2}$, defined locally by the following linear combinations of local diagrams in $\mathcal{V}^{A_2}_{s_1,s_2}$ (see \cite{kuperberg:1996, suciu:1997}), for any $q \in \mathbb{C}$:
\begin{center}
\begin{minipage}[b]{12cm}
 \begin{minipage}[t]{2.5cm}
  \mbox{} \\
  \parbox[t]{1cm}{}
 \end{minipage}
 \begin{minipage}[t]{5.5cm}
  \mbox{} \\
   \includegraphics[width=50mm]{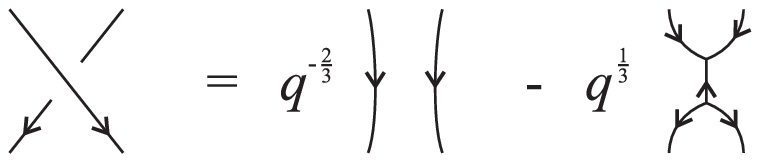}
 \end{minipage}
 \begin{minipage}[t]{3.5cm}
  \hfill
  \parbox[t]{2cm}{\begin{eqnarray}\label{braiding1}\end{eqnarray}}
 \end{minipage}
\end{minipage}
\end{center}

\noindent The braiding satisfies type II and type III Reidemeister moves, and a braiding fusion
relation \cite[Equations (8), (9)]{evans/pugh:2009iii}, provided $\delta = [2]_q$ and $\alpha = [3]_q$, where the quantum integer $[m]_q$ is defined by $[m]_q = (q^m - q^{-m})/(q - q^{-1})$.

Thus it is sufficient to work over $V^{A_2}_{(m,n),(m',n')} := V^{A_2}_{+^m -^n, +^{m'} -^{n'}}$, where $+^k -^l$ is the string of $k$ signs $+$ followed by $l$ signs $-$, since, for any arbitrary string $s$ with $m$ signs $+$ and $n$ signs $-$ and string $s'$ with $m'$ signs $+$ and $n'$ signs $-$, there is an isomorphism $\iota$ between $V^{A_2}_{s,s'}$ and $V^{A_2}_{(m,n), (m',n')}$ given by using the braiding to permute the order of the signs in $s$ to $+^m -^n$, and the inverse braiding to permute the order of the signs in $s'$ to $+^{m'} -^{n'}$.

\begin{figure}[tb]
\begin{center}
  \includegraphics[width=35mm]{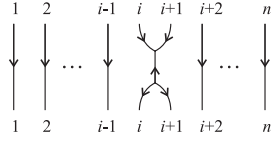}\\
 \caption{The tangle $W_i \in \mathcal{V}^{A_2}_{m}$.}\label{fig:W_i}
\end{center}
\end{figure}

A diagrammatic representation of the Hecke algebra for $SU(3)$ is given as follows: Let $W_i \in \mathcal{V}^{A_2}_{m} := \mathcal{V}^{A_2}_{(m,0),(m,0)}$ be the tangle illustrated in Figure \ref{fig:W_i}.
A $\ast$-operation can be defined on $\mathcal{V}^{A_2}_{m}$, where for an $m$-tangle $T$, $T^{\ast}$ is the $m$-tangle obtained by reflecting $T$ about a horizontal line halfway between the top and bottom vertices of the tangle, and reversing the orientations on every string. Then $\ast$ on $\mathcal{V}^{A_2}_{m}$ is the conjugate linear extension of $\ast$ on $m$-tangles. For $\delta \in \mathbb{R}$ (so $q \in \mathbb{R}$ or $q$ a root of unity), the $\ast$-operation leaves the Kuperberg ideal invariant due to the symmetry of the relations K1-K3.
For $m \in \mathbb{N} \cup \{ 0 \}$ we define the algebra $A_2\textrm{-}TL_m$ to be the quotient space $V^{A_2}_m := V^{A_2}_{(m,0),(m,0)}$.
The images of the $W_i$'s in $V^{A_2}_m$ are clearly self-adjoint, and satisfy the relations H1-H3 and (\ref{eqn:SU(3)q_condition}) \cite[Section 3.3]{evans/pugh:2009iii}.

\subsection{The $A_2$-Temperley-Lieb category} \label{sect:A2TLcat}

We now review the $A_2$-Temperley-Lieb category, following the construction given in \cite{evans/pugh:2010ii} (see also an earlier construction in \cite{cooper:2007}, and the Temperley-Lieb category in \cite{turaev:1994, yamagami:2003}).

The $A_2$-Temperley-Lieb category is defined by $A_2\textrm{-}TL = \mathrm{Mat}(C^{A_2})$, where $C^{A_2}$ is the tensor category whose objects are projections in $V^{A_2}_{(m,n)}$, and whose morphisms are $\mathrm{Hom}(p_1,p_2) = p_2 V^{A_2}_{(m_2,n_2),(m_1,n_1)} p_1$, for projections $p_i \in V^{A_2}_{(m_i,n_i)}$, $i=1,2$. We write $A_2\textrm{-}TL_{(m,n)} = V^{A_2}_{(m,n)}$, and $\rho$, $\overline{\rho}$ for the identity projections in $A_2\textrm{-}TL_{(1,0)}$, $A_2\textrm{-}TL_{(0,1)}$ respectively consisting of a single string with orientation downwards, upwards respectively. Then the identity diagram in $A_2\textrm{-}TL_{(m,n)}$, given by $m+n$ vertical strings where the first $m$ strings have downwards orientation and the next $n$ have upwards orientation, is expressed as $\rho^m \overline{\rho}^n$. It is a linear combination of simple projections $f_{(i,j)}$ for $i,j \geq 0$, $0 \leq i+j < m+n$ such that $i-j \cong m-n \textrm{ mod } 3$, and a simple projection $f_{(m,n)}$, where $f_{(1,0)} = \rho$, $f_{(0,1)} = \overline{\rho}$ and $f_{(0,0)}$ is the empty diagram \cite{evans/pugh:2010ii}.

The morphisms $\mathfrak{f}_{(p,l)} = \mathrm{id}_{f_{(p,l)}}$ are generalized Jones-Wenzl projections which satisfy the recursion relations \cite[(2.1.0)-(2.1.2)]{suciu:1997}:
\begin{center}
\begin{minipage}[b]{16cm}
 \begin{minipage}[t]{1cm}
  \hfill
  \parbox[t]{0.5cm}{}
 \end{minipage}
 \begin{minipage}[t]{13cm}
  \begin{center}
  \mbox{} \\
   \includegraphics[width=80mm]{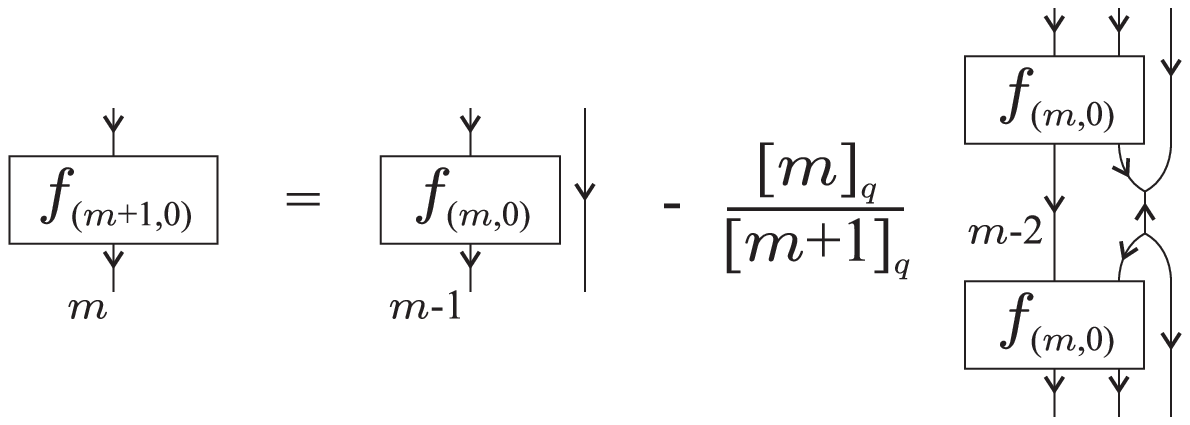}
  \end{center}
 \end{minipage}
 \begin{minipage}[t]{1.5cm}
  \hfill
  \vspace{2mm} \parbox[t]{1.5cm}{\begin{eqnarray}\label{eqn:f(k,0)}\end{eqnarray}}
 \end{minipage}
\end{minipage}
\end{center}

and \cite[(2.1.7)]{suciu:1997}:
\begin{center}
\begin{minipage}[b]{16cm}
 \begin{minipage}[t]{1cm}
  \hfill
  \parbox[t]{0.5cm}{}
 \end{minipage}
 \begin{minipage}[t]{13cm}
  \begin{center}
  \mbox{} \\
   \includegraphics[width=120mm]{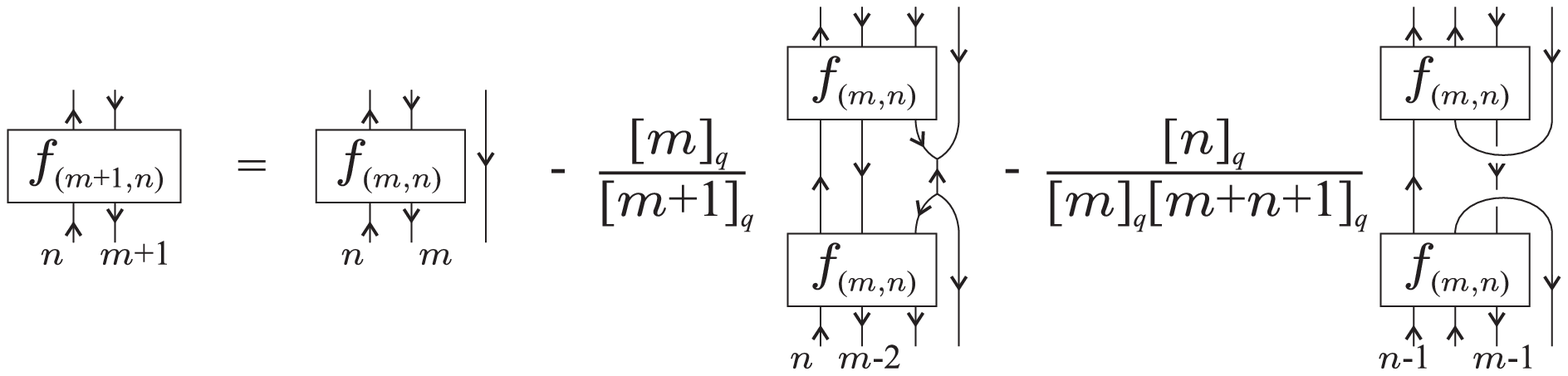}
  \end{center}
 \end{minipage}
 \begin{minipage}[t]{1.5cm}
  \hfill
  \vspace{2mm} \parbox[t]{1.5cm}{\begin{eqnarray}\label{eqn:f(m,n)}\end{eqnarray}}
 \end{minipage}
\end{minipage}
\end{center}

The morphisms $\mathfrak{f}_{(p,l)}$ also satisfy the properties \cite{suciu:1997}:
\begin{center}
\includegraphics[width=150mm]{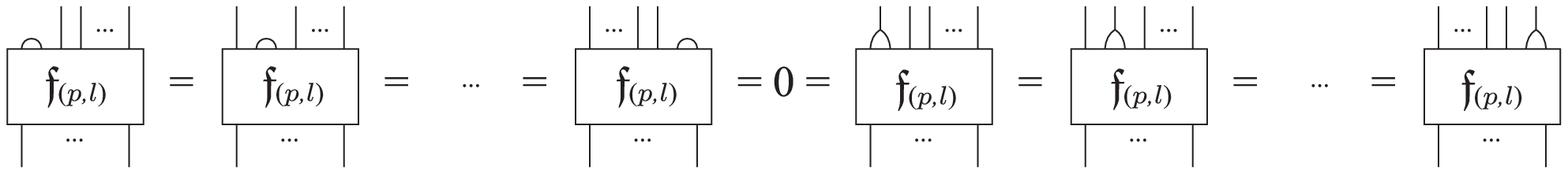} \\
\vspace{5mm}
\includegraphics[width=65mm]{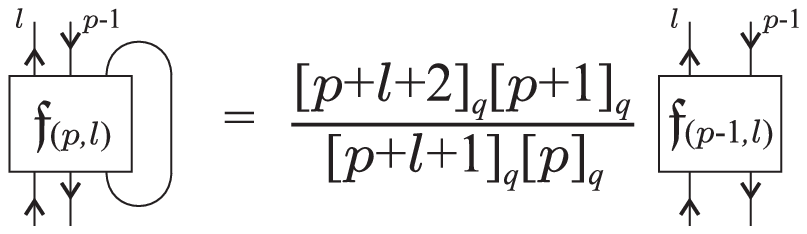}
$$\mathrm{tr}(\mathfrak{f}_{(p,l)}) = \frac{[p+1]_q[l+1]_q[p+l+2]_q}{[2]_q},$$
\end{center}
where $\mathrm{tr}(\mathfrak{f}_{(p,l)})$ is given by connecting the $i^{\mathrm{th}}$ string from the left along the top to the $i^{\mathrm{th}}$ string from the left along the bottom for each $i = 1,\ldots,p+l$.
For $p \geq p'$ and $l \geq l'$, these generalized Jones-Wenzl projections also satisfy the property $\mathfrak{f}_{(p,l)} (\mathrm{id}_{\rho^i \overline{\rho}^j} \otimes \mathfrak{f}_{(p',l')} \otimes \mathrm{id}_{\rho^{p-p'-i} \overline{\rho}^{l-l'-j}} = \mathfrak{f}_{(p,l)} = (\mathrm{id}_{\rho^i \overline{\rho}^j} \otimes \mathfrak{f}_{(p',l')} \otimes \mathrm{id}_{\rho^{p-p'-i} \overline{\rho}^{l-l'-j}}) \mathfrak{f}_{(p,l)}$, for any $0 \leq i \leq p-p'$, $0 \leq j \leq l-l'$. This property also holds if we conjugate either $\mathfrak{f}_{(p,l)}$ or $\mathrm{id}_{\rho^i \overline{\rho}^j} \otimes \mathfrak{f}_{(p',l')} \otimes \mathrm{id}_{\rho^{p-p'-i} \overline{\rho}^{l-l'-j}}$ by any braiding.

The $f_{(p,l)}$ satisfy the fusion rules for $SU(3)$ \cite{evans/pugh:2010ii}:
\begin{equation} \label{eqn:fusion_rule-f(k,l)}
f_{(p,l)} \otimes \rho \cong f_{(p,l-1)} \oplus f_{(p-1,l+1)} \oplus f_{(p+1,l)},
\quad
f_{(p,l)} \otimes \overline{\rho} \cong f_{(p-1,l)} \oplus f_{(p+1,l-1)} \oplus f_{(p,l+1)}.
\end{equation}

In the generic case, $\delta \geq 2$, the $A_2$-Temperley-Lieb category $A_2\textrm{-}TL$ is semisimple and for any pair of non-isomorphic simple projections $p_1$, $p_2$ we have $\langle p_1,p_2 \rangle = 0$.
We recover the infinite graph $\mathcal{A}^{(\infty)}$, where the vertices are labeled by the projections $f_{(p,l)}$ and the edges represent tensoring by $\rho$.

In the non-generic case where $q$ is an $k+3^{\mathrm{th}}$ root of unity, we have $\mathrm{tr}(\mathfrak{f}_{(p,l)}) = 0$ for $p+l=k+1$. By (\ref{eqn:f(m,n)}), $\mathfrak{f}_{(p',l')} = 0$ for all $p',l' \geq k+2$ if $\mathfrak{f}_{(p,l)} = 0$ for $p+l=k+1$. Thus the negligible morphisms are the ideal $\langle \mathfrak{f}_{(p,l)} | p+l=k+1 \rangle$ generated by $\mathfrak{f}_{(p,l)}$ such that $p+l=k+1$.
The quotient $A_2\textrm{-}TL^{(k)} := A_2\textrm{-}TL/ \langle \mathfrak{f}_{(p,l)} | p+l=k+1 \rangle$ is semisimple with simple objects $f_{(p,l)}$, $p,l \geq 0$ such that $p+l \leq k$ which satisfy the fusion rules (\ref{eqn:fusion_rule-f(k,l)}) where $f_{(p',l')}$ is understood to be zero if $p'<0$, $l'<0$ or $p'+l' \geq k+1$.
Thus we recover the graph $\mathcal{A}^{(k+3)}$, where the vertices are labeled by the projections $f_{(p,l)}$ and the edges represent tensoring by $\rho$.

\subsection{Ocneanu cells and $SU(3)$ module categories} \label{sect:module_categories}

\begin{figure}[tb]
\begin{center}
\includegraphics[width=100mm]{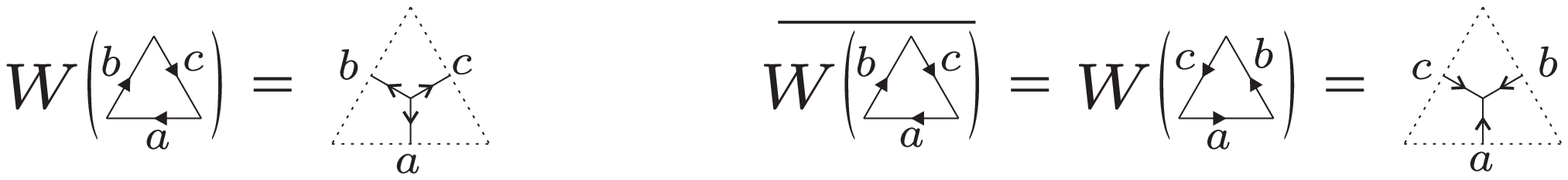}\\
 \caption{Cells associated to trivalent vertices} \label{fig:Oc-Kup}
\end{center}
\end{figure}

In this section we recall the definition of a cell system $W$ on a finite $\mathcal{ADE}$ graph $\mathcal{G}$.
Such a pair $(\mathcal{G},W)$ yields braided subfactor $N \subset M$ and a module category ${}_N \mathcal{X}_M$. We will describe the construction of the path algebra of $\mathcal{G}$ via a monoidal functor, which is essentially the module category ${}_N \mathcal{X}_M$, from the $A_2$-Temperley-Lieb category to $\mathrm{Fun}({}_N \mathcal{X}_M,{}_N \mathcal{X}_M)$.

Ocneanu \cite{ocneanu:2000ii} defined a cell system $W$ on $\mathcal{G}$, associating a complex number $W \left( \triangle_{i,j,k}^{(a,b,c)} \right)$, now called an Ocneanu cell, to each closed loop of length three $\triangle_{i,j,k}^{(a,b,c)}$ in $\mathcal{G}$ as in Figure \ref{fig:Oc-Kup}, where $a,b,c$ are edges on $\mathcal{G}$, and $i,j,k$ are the vertices on $\mathcal{G}$ given by $i=s(a)=r(c)$, $j=s(b)=r(a)$, $k=s(c)=r(b)$. These cells satisfy two properties, called Ocneanu's type I, II equations respectively, which are obtained by evaluating the Kuperberg relations K2, K3 respectively, using the identification in Figure \ref{fig:Oc-Kup}: \\
$(i)$ for any type I frame \includegraphics[width=16mm]{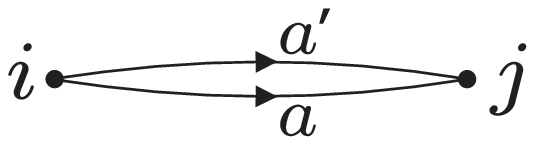} in $\mathcal{G}$ we have
\begin{equation} \label{eqn:typeI_frame}
\sum_{k,b_1,b_2} W \left( \triangle_{i,j,k}^{(a,b_1,b_2)} \right) \overline{W \left( \triangle_{i,j,k}^{(a',b_1,b_2)} \right)} = \delta_{a,a'} [2]_q \phi_i \phi_j
\end{equation}
$(ii)$ for any type II frame \includegraphics[width=30mm]{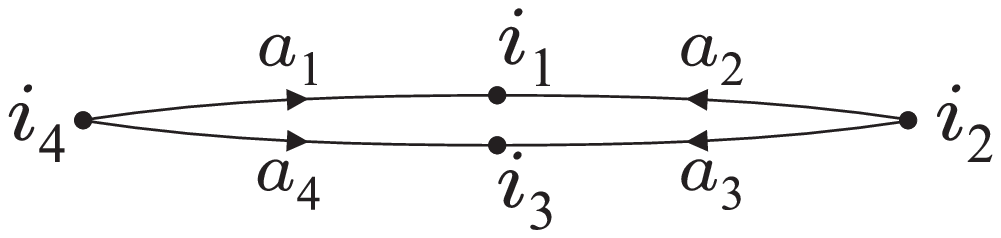} in $\mathcal{G}$ we have
\begin{eqnarray}
\lefteqn{ \sum_{k,b_j} \phi_k^{-1} W \left( \triangle_{i_2,i_1,k}^{(a_2,b_1,b_2)} \right) \overline{W \left( \triangle_{i_2,i_3,k}^{(a_3,b_3,b_2)} \right)} W \left( \triangle_{i_4,i_3,k}^{(a_4,b_3,b_4)} \right) \overline{W \left( \triangle_{i_4,i_1,k}^{(a_1,b_4,b_1)} \right)} } \nonumber \\
& \qquad & = \delta_{a_1,a_4} \delta_{a_2,a_3} \phi_{i_4} \phi_{i_1} \phi_{i_2} + \delta_{a_1,a_2} \delta_{a_3,a_4} \phi_{i_1} \phi_{i_2} \phi_{i_3} \hspace{25mm} \label{eqn:typeII_frame}
\end{eqnarray}
Here $(\phi_v)_v$ is the Perron-Frobenius eigenvector for the Perron-Frobenius eigenvalue $\alpha = [3]_q$ of $\mathcal{G}$.
The existence of these cells for the finite $\mathcal{ADE}$ graphs was claimed by Ocneanu \cite{ocneanu:2000ii}, and shown in \cite{evans/pugh:2009i} with the exception of the graph $\mathcal{E}_4^{(12)}$. These cells define a unitary connection on the graph $\mathcal{G}$ which satisfy the Yang-Baxter equation \cite[Lemma 3.2]{evans/pugh:2009i}. There is up to equivalence precisely one connection on the graphs $\mathcal{A}^{(m)}$, $\mathcal{A}^{(2m+1)\ast}$, $\mathcal{E}^{(8)}$, $\mathcal{E}^{(8)\ast}$, $\mathcal{E}_5^{(12)}$ and $\mathcal{E}^{(24)}$. For the graphs $\mathcal{A}^{(2m)\ast}$ and $\mathcal{E}_2^{(12)}$ there are precisely two inequivalent connections, which are obtained from each other by a $\mathbb{Z}_2$ symmetry of the graph. This $\mathbb{Z}_2$ symmetry is the conjugation of the graph in the case of $\mathcal{E}_2^{(12)}$. There is at least one connection for each graph $\mathcal{D}^{(m)}$, $m \not \equiv 0 \textrm{ mod } 3$, and at least two inequivalent connections for each graph $\mathcal{D}^{(3p)}$, which are the complex conjugates of each other. There is at least one connection for each graph $\mathcal{D}^{(2m+1)\ast}$, and at least two inequivalent connections for each graph $\mathcal{D}^{(2m)\ast}$, which are obtained from each other by a $\mathbb{Z}_2$ symmetry of the graph. There are also at least two inequivalent connections for the graph $\mathcal{E}_1^{(12)}$, which are obtained from each other by conjugation of the graph.

Let ${}_N \mathcal{X}_{N} = \{ \lambda_{(p,l)} | \; 0 \leq p,l,p+l \leq k < \infty \}$ be braided system of endomorphisms of $SU(3)_k$ on a factor $N$, and $N \subset M$ will be a braided inclusion with classifying graph $\mathcal{G} = G_{\rho}$, of $\mathcal{ADE}$ type, arising from the nimrep $G$ of ${}_N \mathcal{X}_N$ acting on ${}_N \mathcal{X}_M$. Then the module category gives rise to a monoidal functor $F$ from the $A_2$-Temperley-Lieb category $A_2\textrm{-}TL^{(k)}$ to $\mathrm{Fun}({}_N \mathcal{X}_M,{}_N \mathcal{X}_M)$, where $F$ is given by
\begin{equation} \label{eqn:functorF}
F(f_{(p,l)}) = \bigoplus_{i,j \in \mathcal{G}_0} G_{\lambda_{(p,l)}}(i,j) \, \mathbb{C}_{i,j},
\end{equation}
where $\lambda_{(p,l)}$ is an irreducible endomorphism in ${}_N \mathcal{X}_N$ identified with the generalized Jones-Wenzl projections $f_{(p,l)}$.
The $\mathbb{C}_{i,j}$ are 1-dimensional $R$-$R$ bimodules, where $R = (\mathbb{C}\mathcal{G})_0$. The category of $R$-$R$ bimodules has a natural monoidal structure given by $\otimes_R$.

We denote by $\mathcal{G}^{\mathrm{op}}$ the opposite graph of $\mathcal{G}$ obtained by reversing the orientation of every edge of $\mathcal{G}$.
Then we have that $F(\rho^m \overline{\rho}^n)$ is the $R$-$R$ bimodule with basis given by all paths of length $m+n$ on
$\mathcal{G}$, $\mathcal{G}^{\mathrm{op}}$, where the first $m$ edges are on $\mathcal{G}$ and the last $n$ edges are on $\mathcal{G}^{\mathrm{op}}$.
In particular $F(\rho^m) = (\mathbb{C}\mathcal{G})_m$, so that we have the graded algebra $\bigoplus_m F(\rho^m) = (\mathbb{C}\mathcal{G})$, the path algebra of $\mathcal{G}$.
The endomorphisms $\rho^m$ are not irreducible however, but decompose into direct sums of the generalized Jones-Wenzl projections $f_{(p,0)}$.
The natural algebra to consider is thus the graded algebra $\Sigma = \bigoplus_j F(f_{(j,0)})$, where the $p^{\mathrm{th}}$ graded part is $\Sigma_p = F(f_{(p,0)}) = F(\lambda_{(p,0)})$. The multiplication $\mu$ is defined by $\mu_{p,l} = F(\mathfrak{f}_{(p+l,0)}): \Sigma_p \otimes_R \Sigma_l \rightarrow \Sigma_{p+l}$, where $\mathfrak{f}_{(p,l)} = \mathrm{id}_{f_{(p,l)}}$.

To define the functor $F$ it is sufficient to define it on the morphisms in the $A_2$-Temperley-Lieb category.
If $a \in \mathcal{G}_1$ is an edge on $\mathcal{G}$, we denote by $\widetilde{a} \in \mathcal{G}^{\mathrm{op}}_1$ the corresponding edge with opposite orientation on $\mathcal{G}^{\mathrm{op}}$.
We define annihilation operators $c_l$, $c_r$ by:
\begin{equation} \label{eqn:annihilation-lr}
c_l(a\widetilde{b}) = \delta_{s(a),s(b)} \frac{\sqrt{\phi_{r(a)}}}{\sqrt{\phi_{s(a)}}} \, s(a), \qquad c_r(\widetilde{b}a) = \delta_{r(a),r(b)} \frac{\sqrt{\phi_{s(a)}}}{\sqrt{\phi_{r(a)}}} \, r(a),
\end{equation}
and creation operators $c_l^{\ast}$, $c_r^{\ast}$ as their adjoints, where $a$ is an edge on $\mathcal{G}$ and $\widetilde{b}$ an edge on $\mathcal{G}^{\mathrm{op}}$.
Define the following fork operators $\curlyvee$, $\overline{\curlyvee}$ by:
\begin{eqnarray}
\curlyvee(\widetilde{a}) & = & \frac{1}{\sqrt{\phi_{s(a)} \phi_{r(a)}}} \sum_{b_1,b_2} W(\triangle_{s(a),r(a),r(b_1)}^{(a,b_1,b_2)}) b_1 b_2, \label{eqn:Yfork(in)} \\
\overline{\curlyvee}(a) & = & \frac{1}{\sqrt{\phi_{s(a)} \phi_{r(a)}}} \sum_{b_1,b_2} \overline{W(\triangle_{s(a),r(a),r(b_1)}^{(a,b_1,b_2)})} \widetilde{b_1} \widetilde{b_2}, \label{eqn:Yfork(out)}
\end{eqnarray}
where $\triangle_{s(a),r(a),r(b_1)}^{(a,b_1,b_2)}$ denotes the closed loop of length 3 on $\mathcal{G}$ along the edge $a$, $b_1$ and $b_2$, and $W$ is a cell system on $\mathcal{G}$ constructed in \cite{evans/pugh:2009i}. We also define $\curlywedge = \overline{\curlyvee}^{\ast}$ and $\overline{\curlywedge} = \curlyvee^{\ast}$.
Then the functor $F$ is defined on the morphisms of $A_2\textrm{-}TL$ by assigning to the left, right caps the annihilation operators $c_l$, $c_r$ respectively, given by (\ref{eqn:annihilation-lr}), to the left, right cups the creation operators $c_l^{\ast}$, $c_r^{\ast}$ respectively, to the incoming, outgoing Y-forks the operators $\curlyvee$, $\overline{\curlyvee}$ respectively given in (\ref{eqn:Yfork(in)}), (\ref{eqn:Yfork(out)}), and to the incoming, outgoing inverted Y-forks the operators $\curlywedge$, $\overline{\curlywedge}$ respectively.

\begin{figure}[tb]
\begin{center}
   \includegraphics[width=70mm]{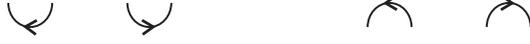}
\caption{left and right cups; left and right caps} \label{fig:cups&caps}
\end{center}
\end{figure}

\begin{figure}[tb]
\begin{center}
   \includegraphics[width=70mm]{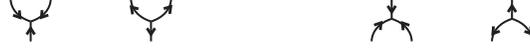}
\caption{incoming and outgoing Y-forks; incoming and outgoing inverted Y-forks} \label{fig:Y-forks}
\end{center}
\end{figure}

Another graded algebra $\Pi$ can be defined by $\Pi = \mathbb{C}\mathcal{G}/\langle \mathrm{Im} \big( F \big( \includegraphics[width=5mm]{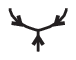} \big) \big) \rangle$, which is the quotient of $\mathbb{C}\mathcal{G}$ by the two-sided ideal generated by the image of the operators $\includegraphics[width=5mm]{fig_nakayama-Yfork}$ in $\mathbb{C}\mathcal{G}$. Its $p^{\mathrm{th}}$ graded part is $\Pi_p = (\mathbb{C}\mathcal{G})_p/\langle \mathrm{Im} \big( F \big( \includegraphics[width=5mm]{fig_nakayama-Yfork} \big) \big) \rangle_p$, where $\langle \mathrm{Im} \big( F \big( \includegraphics[width=5mm]{fig_nakayama-Yfork} \big) \big) \rangle_p$ is the restriction of $\langle \mathrm{Im} \big( F \big( \includegraphics[width=5mm]{fig_nakayama-Yfork} \big) \big) \rangle$ to $(\mathbb{C}\mathcal{G})_p$, which is equal to $\sum_{i=1}^{p-1}\mathrm{Im}(F(\mathfrak{U}_i))$, the union of the images on $\mathbb{C}\mathcal{G}_p$ of the morphisms $\mathfrak{U}_i = \mathrm{id}_{\delta^{-1} W_i}$.
This algebra is isomorphic to the graded algebra $\Sigma$ defined above \cite[Prop. 7.4.18]{cooper:2007} (see also \cite{evans/pugh:2010ii}).
The graded algebra $\Pi$ is the almost Calabi-Yau algebra of \cite{evans/pugh:2010ii}, which we will return to in Section \ref{sect:almostCYalg}. These almost Calabi-Yau algebras are the $SU(3)$ analogue of the preprojective algebras \cite{gelfand/ponomarev:1979} for an $ADE$ Dynkin diagram $\mathcal{G}$, which is the graded algebra defined by $\Pi' = \mathbb{C}\mathcal{G}/\langle \mathrm{Im} \left( F \left( \includegraphics[width=5mm]{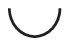} \right) \right) \rangle$, where $\langle \mathrm{Im} \left( F \left( \includegraphics[width=5mm]{fig_nakayama-cup} \right) \right) \rangle \subset \mathbb{C}\mathcal{G}$ is the two-sided ideal generated by the image of the creation operators $\includegraphics[width=5mm]{fig_nakayama-cup}$ in $\mathbb{C}\mathcal{G}$. The preprojective algebra $\Pi'$ appears in the subfactor literature: the $p^{\mathrm{th}}$ graded part $\Pi'_p$ of the preprojective algebra $\Pi'$ is isomorphic to the space of essential paths $\mathrm{EssPath}_p = \mathrm{ker}(F(\mathfrak{e}_1) \vee \cdots \vee F(\mathfrak{e}_{p-1}))$ of Ocneanu \cite{ocneanu:2000i}.

Let $\Lambda$ be the graded coalgebra $\Lambda = F(f_{(0,0)}) \oplus F(f_{(1,0)}) \oplus F(f_{(0,1)}) \oplus F(f_{(0,0)}) = (\mathbb{C}\mathcal{G})_0 \oplus (\mathbb{C}\mathcal{G})_1 \oplus (\mathbb{C}\mathcal{G}^{\mathrm{op}})_1 \oplus (\mathbb{C}\mathcal{G})_0$. The comultiplication is $\Delta$ is given by $\Delta_{1,1} = F \big( \includegraphics[width=5mm]{fig_nakayama-Yfork} \big): \Lambda_2 \rightarrow \Lambda_1 \otimes_R \Lambda_1$, $\Delta_{1,2} = F \big( \includegraphics[width=4mm]{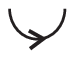} \big): \Lambda_3 \rightarrow \Lambda_1 \otimes_R \Lambda_2$, $\Delta_{2,1} = F \big( \includegraphics[width=4mm]{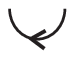} \big): \Lambda_3 \rightarrow \Lambda_2 \otimes_R \Lambda_1$, and the other comultiplications are trivial.
The pair $(\Pi,\Lambda)$ is almost Koszul, in the sense of \cite{brenner/butler/king:2002}, where the algebra $\Pi$ is a $(k,3)$-Koszul algebra \cite[Corollary 7.4.19]{cooper:2007} (see also \cite{evans/pugh:2010ii}).

The monoidal functor $F$ is closely related to the $SU(3)$ Goodman-de la Harpe-Jones construction of \cite{evans/pugh:2009ii} and its manifestation in the $SU(3)$-graph planar algebra construction \cite{evans/pugh:2009iv}, which we will review in the next section.

\subsection{$A_2$-planar algebras} \label{sect:A2PA}

We will now review the basics of $A_2$-planar algebras from \cite{evans/pugh:2009iii}. These $A_2$-planar algebras were useful to understand the double complexes of finite dimensional algebras which arise in the context of $SU(3)$ subfactors and modular invariants.
They are a direct generalization of the planar algebras of Jones \cite{jones:planar}, which have been extensively studied since their introduction ten years ago. Jones' planar algebras naturally contain the Temperley-Lieb algebra which encodes the representation theory of quantum $SU(2)$.
Our $A_2$-planar algebras naturally encode the representation theory of quantum $SU(3)$, or in the dual Hecke picture, the $A_2$-Temperley-Lieb algebra.

\begin{figure}[htb]
\begin{center}
  \includegraphics[width=60mm]{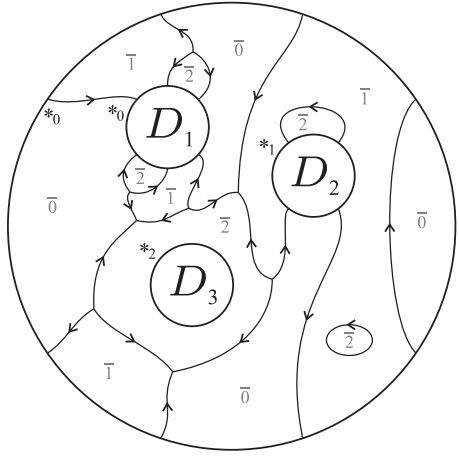}\\
 \caption{$A_2$-planar $\sigma$-tangle for $\sigma = -+-+-+-+$}\label{fig:A2_planar_tangle}
\end{center}
\end{figure}

Let $\sigma = \sigma_1 \cdots \sigma_m$ be a sign string, $\sigma_j \in \{ \pm \}$, such that the difference between the number of `$+$' and `$-$' is 0 mod 3. An $A_2$-planar $\sigma$-tangle will be the unit disc $D=D_0$ in $\mathbb{C}$ together with a finite (possibly empty) set of disjoint sub-discs $D_1, D_2, \ldots, D_n$ in the interior of $D$. Each disc $D_k$, $k \geq 0$, will have $m_k \geq 0$ vertices on its boundary $\partial D_k$, whose orientations are determined by sign strings $\sigma^{(k)} = \sigma^{(k)}_1 \cdots \sigma^{(k)}_{m_k}$ where `$+$' denotes a sink and `$-$' a source. The disc $D_k$ will be said to have pattern $\sigma^{(k)}$.
Inside $D$ we have an $A_2$-tangle where the endpoint of any string is either a trivalent vertex (see Figure \ref{fig:A2-webs}) or one of the vertices on the boundary of a disc $D_k$, $k=0, \ldots, n$, or else the string forms a closed loop. Each vertex on the boundaries of the $D_k$ is the endpoint of exactly one string, which meets $\partial D_k$ transversally.
An example of an $A_2$-planar $\sigma$-tangle is illustrated in Figure \ref{fig:A2_planar_tangle} for $\sigma = -+-+-+-+$.

The regions inside $D$ have as boundaries segments of the $\partial D_k$ or the strings. These regions are labelled $\overline{0}$, $\overline{1}$ or $\overline{2}$, called the colouring, such that if we pass from a region $R$ of colour $\overline{a}$ to an adjacent region $R'$ by passing to the right over a vertical string with downwards orientation, then $R'$ has colour $\overline{a+1}$ (mod 3). We mark the segment of each $\partial D_k$ between the last and first vertices with $\ast_{b_k}$, $b_k \in \{0,1,2\}$, so that the region inside $D$ which meets $\partial D_k$ at this segment is of colour $\overline{b_k}$, and the choice of these $\ast_{b_k}$ must give a consistent colouring of the regions.
For each $\sigma$ we have three types of tangle, depending on the colour $\overline{b}$ of the marked segment, or of the marked region near $\partial D$ for $\sigma = \varnothing$.

An $A_2$-planar $\sigma$-tangle $T$ with an internal disc $D_l$ with pattern $\sigma_l = \sigma'$ can be composed with an $A_2$-planar $\sigma'$-tangle $S$ with external disc $D'$ and $\ast_{D'}=\ast_{D_l}$, giving a new $\sigma$-tangle $T \circ_l S$, by inserting the $A_2$-tangle $S$ inside the inner disc $D_l$ of $T$ such that the vertices on the outer disc of $S$ coincide with those on the disc $D_l$ and the regions marked by $\ast$ also coincide. The boundary of the disc $D_l$ is removed, and the strings smoothed if necessary.
Let $\widetilde{\mathcal{P}}$ be the collection of all diffeomorphism classes of such $A_2$-planar tangles, with composition defined as above.
The $A_2$-planar operad $\mathcal{P}$ is the quotient of $\widetilde{\mathcal{P}}$ by the Kuperberg relations K1-K3.

An $A_2$-planar algebra is then defined to be an algebra over this operad, i.e. a family
$P = \left( P_{\sigma}^{\overline{a}}| \, \textrm{sign strings } \sigma, a \in \{0,1,2\} \right)$
of vector spaces
with the following property: for every $\sigma$-tangle $T \in \mathcal{P}_{\sigma}$ with outer disc marked by $\ast_b$, and with $n$ internal discs $D_j$ pattern $\sigma_k$, outer disc marked by $\ast_{b_k}$ and labelled by elements $x_j \in P_{k_j}$, $j=1,\ldots,n$, there is associated a linear map $Z(T): \otimes_{k=1}^n P_{\sigma_k}^{\overline{b_k}} \longrightarrow P_{\sigma}^{\overline{b}}$ which is compatible with the composition of tangles in the following way. If $S$ is a $\sigma_k$-tangle with internal discs $D_{n+1}, \ldots, D_{n+m}$, where $D_k$ has pattern $\sigma_k$, then the composite tangle $T \circ_l S$ is a $\sigma$-tangle with $n+m-1$ internal discs $D_k$, $k = 1,2, \ldots l-1, l+1, l+2, \ldots, n+m$. From the definition of an operad, associativity means that the following diagram commutes:
\begin{equation} \label{eqn:compatability_condition_for_Z(T)}
\xymatrix{
{\left( \bigotimes_{\stackrel{k=1}{\scriptscriptstyle{k \neq l}}}^n P_{\sigma_k}^{\overline{b_k}} \right) \otimes \left( \bigotimes_{k=n+1}^{n+m} P_{\sigma_k}^{\overline{b_k}} \right)} \ar[d]_{\mathrm{id} \otimes Z(S)} \ar[dr]^(.6){Z(T \circ_l S)} \\
{\bigotimes_{k=1}^n P_{\sigma_k}^{\overline{b_k}}} \ar[r]_{Z(T)} & P_{\sigma}^{\overline{b}} }
\end{equation}
so that $Z(T \circ_l S) = Z(T')$, where $T'$ is the tangle $T$ with $Z(S)$ used as the label for disc $D_l$. We also require $Z(T)$ to be independent of the ordering of the internal discs, that is, independent of the order in which we insert the labels into the discs.
When $\sigma = \varnothing$, we will often write $P_{\varnothing}^{\overline{a}}$ as $P_{\overline{a}}$, and we adopt the convention that the empty tensor product is the complex numbers $\mathbb{C}$.

Let $\sigma^{\ast}$ be the sign string obtained by reversing the string $\sigma$ and flipping all its signs.
When each $P_{\sigma\sigma^{\ast}}$ is a $\ast$-algebra, the adjoint $T^{\ast} \in \mathcal{P}_{\sigma^{\ast}}$ of a tangle $T \in \mathcal{P}_{\sigma}$ is defined by reflecting the whole tangle about the horizontal line that passes through its centre and reversing all orientations. The labels $x_k \in P_{\sigma_k}$ of $T$ are replaced by labels $x_k^{\ast}$ in $T^{\ast}$, where $x_k^{\ast}$ is the unique element in $P_{\sigma_k^{\ast}}$ such that $m(x_k,y)^{\ast} = m(y^{\ast},x_k^{\ast})$ for some tangle $y \in \mathcal{P}_{\sigma_k}$ with no internal discs, and where $m(\cdot,\cdot) \in \mathcal{P}_{\sigma_k \sigma_k^{\ast}}$ is the tangle defined in \cite[Section 4.4] {evans/pugh:2009iii}. For any linear combination of tangles in $\mathcal{P}_{\sigma}$ the involution is the conjugate linear extension.
Then $P$ is an $A_2$-planar $\ast$-algebra if each $P_{\sigma\sigma^{\ast}}$ is a $\ast$-algebra, and for a $\sigma$-tangle $T$ with internal discs $D_k$ with patterns $\sigma_k$, labelled by $x_k \in P_{\sigma_k}$, $k=1, \ldots, n$, we have
$Z(T)^{\ast} = Z(T^{\ast})$,
where the labels of the discs in $T^{\ast}$ are $x_k^{\ast}$, and where the definition of $Z(T)^{\ast}$ is extended to linear combinations of $\sigma$-tangles by conjugate linearity.
Note a typographical error in the definition of an $A_2$-planar $\ast$-algebra in \cite[Section 4.4]{evans/pugh:2009iii}, with $P_{\sigma\sigma^{\ast}}$ above incorrectly given as $P_{\sigma}$.

In \cite{evans/pugh:2009iv} we introduced an $A_2$-graph planar algebra construction for an $SU(3)$ $\mathcal{ADE}$ graph $\mathcal{G}$. The $A_2$-graph planar algebra $P^{\mathcal{G}}$ of an $SU(3)$ $\mathcal{ADE}$ graph $\mathcal{G}$ is the path algebra (in the operator algebraic sense) on $\mathcal{G}$ and $\mathcal{G}^{\mathrm{op}}$.
The presenting map $Z:\mathcal{P} \rightarrow P^{\mathcal{G}}$ is defined uniquely \cite[Theorem 5.1]{evans/pugh:2009iv}, up to isotopy, by first isotoping the strings of $T$ in such a way that the diagram $T$ may be divided into horizontal strips so that each horizontal strip only contains the following elements: a (left or right) cup, a (left or right) cap, an (incoming or outgoing) Y-fork, or an (incoming or outgoing) inverted Y-fork, see Figures \ref{fig:cups&caps} and \ref{fig:Y-forks}. Then $Z$ assigns to the left, right caps the annihilation operators $c_l$, $c_r$ respectively, given by (\ref{eqn:annihilation-lr}), to the left, right cups the creation operators $c_l^{\ast}$, $c_r^{\ast}$ respectively, to the incoming, outgoing Y-forks the operators $\curlyvee$, $\overline{\curlyvee}$ respectively given in (\ref{eqn:Yfork(in)}), (\ref{eqn:Yfork(out)}), and to the incoming, outgoing inverted Y-forks the operators $\curlywedge$, $\overline{\curlywedge}$ respectively.

Then we have a tower of algebras $P^{\mathcal{G}}_{0,0} \subset P^{\mathcal{G}}_{0,1} \subset P^{\mathcal{G}}_{0,2} \subset \cdots \;$, where the inclusion $P^{\mathcal{G}}_m \subset P^{\mathcal{G}}_{m+1}$ is given by the $m$-$m+1$ part of the graph $\mathcal{G}$.
There is a positive definite inner product defined from the trace on $P^{\mathcal{G}}$.
We have the inclusion $A_2\textrm{-}PTL_{0,m} := Z(\mathcal{V}^{A_2}_m) \subset P^{\mathcal{G}}_{0,m}$ for each $m$, and we have a double sequence
$$\begin{array}{ccccccc}
A_2\textrm{-}PTL_{0,0} & \subset & A_2\textrm{-}PTL_{0,1} & \subset & A_2\textrm{-}PTL_{0,2} & \subset & \cdots \\
\cap & & \cap & & \cap & & \\
P^{\mathcal{G}}_{0,0} & \subset & P^{\mathcal{G}}_{0,1} & \subset & P^{\mathcal{G}}_{0,2} & \subset & \cdots
\end{array}$$
Then $A_2\textrm{-}PTL = Z(\mathcal{V}^{A_2})$ is the embedding of the $A_2$-Temperley-Lieb algebra into the path algebra of $\mathcal{G}$, which is used to construct the $A_2$-Goodman-de la Harpe-Jones subfactors \cite{evans/pugh:2009ii}.
Let $\overline{P^{\mathcal{G}}}$ denote the GNS-completion of $P^{\mathcal{G}}$ with respect to the trace.
Then for $q = Z(\ast_{\mathcal{G}})$ the minimal projection in $P^{\mathcal{G}}_0$ corresponding to the vertex $\ast_{\mathcal{G}}$ of $\mathcal{G}$, we have an inclusion $q \overline{A_2\textrm{-}PTL} \subset q \overline{P^{\mathcal{G}}} q$ which gives the $A_2$-Goodman-de la Harpe-Jones subfactor $N_{\mathcal{A}} \subset N_{\mathcal{G}}$, where $N_{\mathcal{A}}' \cap N_{\mathcal{G}} = q P^{\mathcal{G}}_{0,0} q = \mathbb{C}$ and the sequence $\{ q A_2\textrm{-}PTL_{0,m} \subset P^{\mathcal{G}}_{0,m} \}_m$ is a periodic sequence of commuting squares of period 3, in the sense of Wenzl in \cite{wenzl:1988}.
Thus we obtain a commuting square of inclusions:
\begin{equation} \label{eqn:NA-NG_commuting_square}
\begin{array}{ccc}
N_{A} & \subset & N_{\mathcal{G}} \\
\cap & & \cap \\
M_{A} & \subset & M_{\mathcal{G}}
\end{array}
\end{equation}
which allows us to compute the dual canonical endomorphism $\theta$ of the $A_2$-Goodman-de la Harpe-Jones subfactor, from which we constructed the nimrep graph $\mathcal{G} = G_{\rho}$ \cite{evans/pugh:2009ii}, see Section \ref{sect:SU(3)MI}.

Since the functor $F$ in Section \ref{sect:module_categories} is defined by the annihilation and creation operators given by (\ref{eqn:annihilation-lr}), and the incoming, outgoing (inverted) Y-fork operators given in (\ref{eqn:Yfork(in)}), (\ref{eqn:Yfork(out)}), we see that $F$ is equivalent to the presenting map $Z$ above. The embedding $P^A \subset P^{\mathcal{G}}$ is given by the image under $F$ of the morphisms in the $A_2$-Temperley-Lieb category.

\subsection{$A_2$-planar modules} \label{sect:A2PAmodules}

We now review modules over $A_2$-planar algebras, which were defined in \cite{evans/pugh:2009iv}, including a partial decomposition of the $A_2$-graph planar algebras of Section \ref{sect:A2PA} into irreducible Hilbert $A_2$-$PTL$-modules. This is the analogue of the notion of Hilbert $TL$-modules over a planar algebra \cite{jones:2001} and the decomposition of bipartite graph planar algebras for the $ADE$ graphs into irreducible Hilbert $TL$-modules \cite{reznikoff:2005}.

A module over and $A_2$-planar algebra $P$, called a $P$-module, is a graded vector space $V=(V_{\sigma,\overline{a}} \, | \, \textrm{sign strings } \sigma, a \in \{0,1,2\})$ with an action of $P$.
An $A_2$-annular $(\sigma,\sigma')$-tangle is an $A_2$-$\sigma$-tangle with a choice of distinguished internal disc which has pattern $\sigma'$.
Then given an $A_2$-annular $(\sigma,\sigma')$-tangle $T$ in $\mathcal{P}$ with outer disc marked by $\ast_a$, with a distinguished ($V$ input) internal disc $D_1$ marked by $\ast_{a'}$ and with pattern $\sigma'$, and with other ($P$ input) internal discs $D_p$, $p=2, \ldots, n$, marked by $\ast_{a_p}$, with patterns $\sigma^{(p)}$, there is a linear map $Z(T): V_{\sigma',\overline{a'}} \otimes \left( \otimes_{p=2}^n P^{\overline{a_p}}_{\sigma^{(p)}} \right) \rightarrow V_{\sigma,\overline{a}}$. The map $Z(T)$ satisfies the same compatability condition (\ref{eqn:compatability_condition_for_Z(T)}) for the composition of tangles as $P$ itself.

For an $A_2$-$C^{\ast}$-planar algebra $P$, that is, $P$ is a non-degenerate finite-dimensional $A_2$-planar $\ast$-algebra with positive definite partition function, a $P$-module $H$ is called a Hilbert $P$-module if each $H_{\sigma}$ is a finite dimensional Hilbert space with an invariant inner-product $\langle \cdot, \cdot \rangle$.
All irreducible Hilbert $A_2$-$PTL$-modules of lowest weight zero were determined in \cite{evans/pugh:2009iv}. We now review these modules.
For any irreducible Hilbert $A_2$-$PTL$-module $H$ of weight zero, the dimension of $H_{\varnothing,\overline{a}}$ is either 0 or 1.

\begin{figure}[htb]
\begin{center}
  \includegraphics[width=100mm]{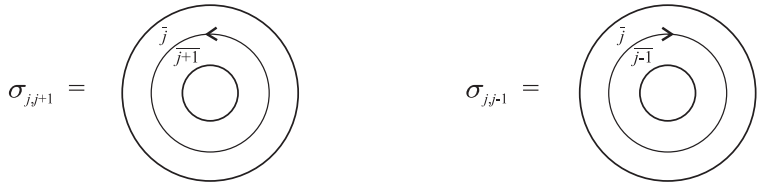}
 \caption{$\sigma_{j,j+1}$ and $\sigma_{j,j-1}$} \label{fig:ATL4}
\end{center}
\end{figure}

The $A_2$-$PTL$-module $V^{\beta}$ where $\beta \in \mathbb{C}$, $|\beta| \leq \alpha = [3]_q$, is the $A_2$-$PTL$-module of weight zero in which $\sigma^{(a)}_{(k_1,k_2)} = \beta^{k_1} \overline{\beta}^{k_2}$ for $a=0,1,2$, where
$$\sigma^{(i,\varepsilon)}_{(k_1,k_2)} = (\sigma_{i,i\varepsilon1} \sigma_{i\varepsilon1,i\varepsilon2} \cdots \sigma_{i\varepsilon k_1-\varepsilon1,i\varepsilon k_1})(\sigma_{i\varepsilon k_1,i\varepsilon k_1-\varepsilon1} \sigma_{i\varepsilon k_1-\varepsilon1,i\varepsilon k_1-\varepsilon2} \cdots \sigma_{i\varepsilon k_1-\varepsilon k_2\varepsilon1,i\varepsilon k_1-\varepsilon k_2})$$
for $i \in \{ 0,1,2 \}$, $\varepsilon \in \{ \pm \}$, and non-negative integers $k_1 \equiv k_2 \textrm{ mod } 3$, where the $0$-tangles $\sigma_{j,j \pm 1}$ are illustrated in Figure \ref{fig:ATL4}, $j \in \{ 0,1,2 \}$.

When $\beta = \alpha$, $V_{\sigma}^{\alpha} = A_2\textrm{-}PTL_{\sigma}$. For $\alpha > 3$ (which corresponds to $\delta > 2$), the inner product is positive definite and $H_{\sigma}^{\alpha} = V_{\sigma}^{\alpha}$ is a Hilbert $A_2$-$PTL$-module.
When $0 < |\beta| < \alpha$, $V^{\beta}$ is the module such that $V^{\beta}_{\sigma}$ has as basis the set of all $(\sigma,\overline{0})$-tangles with no contractible circles and at most two non-contractible circles. The action of $A_2$-$PTL$ on $V^{\beta}$, $0 < |\beta| < \alpha$, is given as follows. Let $T \in A_2\textrm{-}PTL(\sigma',\sigma)$ and $R \in A_2\textrm{-}PTL_{\sigma}$. We form the tangle $TR$ and reduce it using K1-K3, so that $TR = \sum_j \delta^{b_j} \alpha^{c_j} (TR)_j$, for some basis $A_2$-annular $(\sigma',\overline{0})$-tangles $(TR)_j$, where $b_j$, $c_j$ are non-negative integers. Let $\sharp^a_j$, $\sharp^c_j$ denote the number of non-contractible circles in the tangle $(TR)_j$ which have anti-clockwise, clockwise orientation respectively. We define integers $d_j$, $f_j$ and $g_j$ as follows: $d_j = \textrm{min}(\sharp^a_j, \sharp^c_j)$, $f_j = \sharp^a_j - \sharp^c_j - \gamma_{f_j}$ if $\sharp^a_j \geq \sharp^c_j$ and $f_j = 0$ otherwise, and $g_j = \sharp^c_j - \sharp^a_j - \gamma_{g_j}$ if $\sharp^a_j \leq \sharp^c_j$ and $g_j = 0$ otherwise, where $\gamma_{f_j}, \gamma_{g_j} \in \{ 0,1,2 \}$ such that $f_j, g_j \equiv 0 \textrm{ mod }3$. Then we set $T(R) = \sum_j \delta^{b_j} \alpha^{c_j} \beta^{d_j + f_j} \overline{\beta}^{d_j + g_j} (\widehat{TR})_j$, where $(\widehat{TR})_j$ is the tangle $(TR)_j$ with $d_j + f_j$ anti-clockwise non-contractible circles removed, and $d_j + g_j$ clockwise ones removed.
Given two tangles $S,T \in V^{\beta}_{\sigma}$, we reduce $T^{\ast} S$ using K1-K3 so that $T^{\ast} S = \sum_j \delta^{b_j} \alpha^{c_j} (T^{\ast}S)_j$, where $(T^{\ast}S)_j$ are basis tangles in $A_2$-$PTL_{\varnothing}$. Then the inner-product on $V^{\beta}$ is defined by $\langle S,T \rangle = \sum_j \delta^{b_j} \alpha^{c_j} \beta^{d_j + f_j} \overline{\beta}^{d_j + g_j}$, where $d_j$, $f_j$ and $g_j$ are defined for each $(T^{\ast}S)_j$ as above.

When $\beta = 0$, $V^{0,\overline{a}}$ is the module such that $V^{0,\overline{a}}_{\sigma}$ has as basis the set of all $(\sigma,\overline{a})$-tangles with no contractible or non-contractible circles at all, and is equipped with an $A_2$-$PTL$-module structure of lowest weight zero as follows. Let $T \in A_2\textrm{-}PTL(\sigma',\sigma)$ and $R \in V^{0,\overline{a}}_{\sigma}$. We form $TR$ and reduce it using K1-K3, so that $TR = \sum_j \delta^{b_j} \alpha^{c_j} (TR)_j$ as in the case $0 < |\beta| < \alpha$. We define $(\widehat{TR})_j$ to be zero if there are any non-contractible circles in $(TR)_j$, and $(TR)_j$ otherwise. Then $T(R) = \sum_j \delta^{b_j} \alpha^{c_j} (\widehat{TR})_j$.
Given two tangles $S,T \in V^{0,\overline{a}}_{\sigma}$, we reduce $T^{\ast} S$ using K1-K3 so that $T^{\ast} S = \sum_j \delta^{b_j} \alpha^{c_j} (T^{\ast}S)_j$ for basis $(\overline{a}:\overline{a})$-tangles $(T^{\ast}S)_j$. Then the inner-product on $V^{0,\overline{a}}$ is defined by $\langle S,T \rangle = \sum_j \delta^{b_j} \alpha^{c_j} \langle S,T \rangle_j$, where $\langle S,T \rangle_j$ is defined to be 0 if there are any non-contractible circles in $(T^{\ast}S)_j$, and 1 otherwise.

For $0 < \alpha \leq 3$, if the inner product is positive semi-definite on the $A_2$-$PTL$-module $V$ of weight zero, then the Hilbert module $H$ is the quotient of $V$ by the subspace of vectors of length zero; otherwise $H$ does not exist.

All zero-weight submodules of the $A_2$-graph planar algebra $P^{\mathcal{G}}$ for an $SU(3)$ $\mathcal{ADE}$ graph $\mathcal{G}$ were determined in \cite[Proposition 5.3]{evans/pugh:2009iv}, analogously to the result of Reznikoff \cite[Prop. 13]{reznikoff:2005} for bipartite graph planar algebras.
For the graphs $\mathcal{A}^{(m)}$,
$$P^{\mathcal{A}^{(m)}} \supset \bigoplus_{(l_1,l_2)} H^{\beta_{(l_1,l_2)}},$$
for $m \not \equiv 0 \textrm{ mod } 3$, whilst for $m=3p$, $p \geq 2$,
$$P^{\mathcal{A}^{(3p)}} \supset \bigoplus_{(l_1,l_2)} H^{\beta_{(l_1,l_2)}} \oplus H^{0,\overline{0}},$$
where in both cases the summation is over all $(l_1,l_2) \in \{ (m_1,m_2) | \; 3m_2 \leq m-3, 3m_1 + 3m_2 < 2m-6 \}$, i.e. each $\beta_{(l_1,l_2)}$ is a cubic root of an eigenvalue of $\Delta_{01} \Delta_{12} \Delta_{20}$.
For the $\mathcal{D}$ graphs, we have
$$P^{\mathcal{D}^{(3p)}} \supset \bigoplus_{(l_1,l_2)} H^{\beta_{(l_1,l_2)}} \oplus 3H^{0,\overline{0}},$$
for $p \geq 2$, where the summation is over all $(l_1,l_2) \in \{ (m_1,m_2) | m_2 \leq p-1, m_1 + m_2 < 2p-2, m_1 - m_2 \equiv 0 \textrm{ mod } 3 \}$, whilst for $m \not \equiv 0 \textrm{ mod } 3$,
$$P^{\mathcal{D}^{(m)}} \supset \bigoplus_{(l_1,l_2)} H^{\beta_{(l_1,l_2)}},$$
where the summation is over all $(l_1,l_2) \in \{ (m_1,m_2) | 3m_2 \leq m-3, 3m_1 + 3m_2 < 2m-6 \}$.
The path algebras for $\mathcal{A}^{(m)\ast}$ and $\mathcal{D}^{(m)\ast}$ are identified under the map which sends the vertices $i_l$, $j_l$ and $k_l$ of $\mathcal{D}^{(m)\ast}$ to the vertex $l$ of $\mathcal{A}^{(m)\ast}$, $l=1,2,\ldots, \lfloor l/2 \rfloor$. We have
$$P^{\mathcal{A}^{(m)\ast}} = P^{\mathcal{D}^{(m)\ast}} \supset \bigoplus_{(l_1,l_2)} H^{\beta_{(l_1,l_2)}},$$
where the summation is over all $(l_1,l_2) \in \{ (p,p) | p=0,1,\ldots, \lfloor (m-3)/2 \rfloor \}$. Similarly, the path algebras for $\mathcal{E}^{(8)}$ and $\mathcal{E}^{(8)\ast}$ are identified, and
$$P^{\mathcal{E}^{(8)}} = P^{\mathcal{E}^{(8)\ast}} \supset H^{\beta_{(0,0)}} \oplus H^{\beta_{(3,0)}} \oplus H^{\beta_{(0,3)}} \oplus H^{\beta_{(2,2)}}.$$
For the graphs $\mathcal{E}_i^{(12)}$, $i=1,2,3$, we have
$$P^{\mathcal{E}_i^{(12)}} \supset H^{\beta_{(0,0)}} \oplus 2 H^{\beta_{(2,2)}} \oplus H^{\beta_{(4,4)}}.$$
For the remaining exceptional graphs we have
\begin{eqnarray*}
P^{\mathcal{E}_4^{(12)}} & \supset & H^{\beta_{(0,0)}} \oplus H^{\beta_{(2,2)}} \oplus H^{\beta_{(4,4)}} \oplus 2H^{0,\overline{0}}, \\
P^{\mathcal{E}_5^{(12)}} & \supset & H^{\beta_{(0,0)}} \oplus H^{\beta_{(3,0)}} \oplus H^{\beta_{(0,3)}} \oplus H^{\beta_{(2,2)}} \oplus H^{\beta_{(4,4)}} \oplus H^{0,\overline{1}} \oplus H^{0,\overline{2}}, \\
P^{\mathcal{E}^{(24)}} & \supset & H^{\beta_{(0,0)}} \oplus H^{\beta_{(6,0)}} \oplus H^{\beta_{(0,6)}} \oplus H^{\beta_{(4,4)}} \oplus H^{\beta_{(7,4)}} \oplus H^{\beta_{(4,7)}} \oplus H^{\beta_{(6,6)}} \oplus H^{\beta_{(10,10)}}. \qquad
\end{eqnarray*}

The following irreducible modules of non-zero weight were also constructed in \cite{evans/pugh:2009iv}.
For $m,n \neq 0$, let $V^{(m,n),\gamma}_{+^m -^n} = \mathrm{span}(\varphi_{(m,n)}^l | \; l = 0,1,\ldots,m+n-1)$, where the $+^m -^n$-tangle $\varphi_{(m,n)}$ is illustrated in Figure \ref{fig-varphi(m,n)} and $\gamma \in \mathbb{T}$. Then $V^{(m,n),\gamma} = A_2\mathrm{-}PTL(V^{(m,n),\gamma}_{+^m -^n})$ is an $A_2$-$PTL$-module in which $\varphi_{(m,n)}^{m+n} = \gamma^{(m+2n)/3} \mathbf{1}_{+^m -^n}$.
For all $m,n \neq 0$, we can choose a faithful trace $\mathrm{tr}'$ on $A^{(m,n)} := A_2\textrm{-}PTL_{+^m -^n}/A_2\textrm{-}PTL_{+^m -^n}^{(m,n)}$,
where $A_2$-$PTL_{\sigma}^{(t_1,t_2)}$ the ideal in $A_2$-$PTL_{\sigma}$ of all labelled $A_2$-annular $\sigma$-tangles with rank $(t_1',t_2') < (t_1,t_2)$, which we extend to a trace $\mathrm{tr}$ on $A_2$-$ATL_{+^m -^n}$ by $\mathrm{tr} = \mathrm{tr}' \circ \pi$, where $\pi$ is the quotient map $\pi:A_2\textrm{-}ATL_{+^m -^n} \rightarrow A^{(m,n)}$. The inner product on $A_2$-$ATL(\sigma:+^m -^n)$ is defined by $\langle S,T \rangle = \mathrm{tr}(T^{\ast}S)$ for any $S,T \in A_2\textrm{-}PTL(\sigma:+^m -^n)$.
If $\gamma^{(m+2n)k/3} \neq 1$ for any $k \in \mathbb{N}$, $V^{(m,n),\gamma}$ is irreducible.

For $m=0$, one can obtain an irreducible $A_2$-$PTL$-module $V^{(3,0),\gamma}$ is the module in which $\varphi_{(3,0)} = \gamma \mathbf{1}_{+^3}$, where $\varphi_{(3,0)}$ is illustrated in Figure \ref{fig-varphi(m,n)}. One can obtain an irreducible $A_2$-$PTL$-module $V^{(0,3),\gamma}$ similarly.

\begin{figure}[htb]
\begin{center}
  \includegraphics[width=120mm]{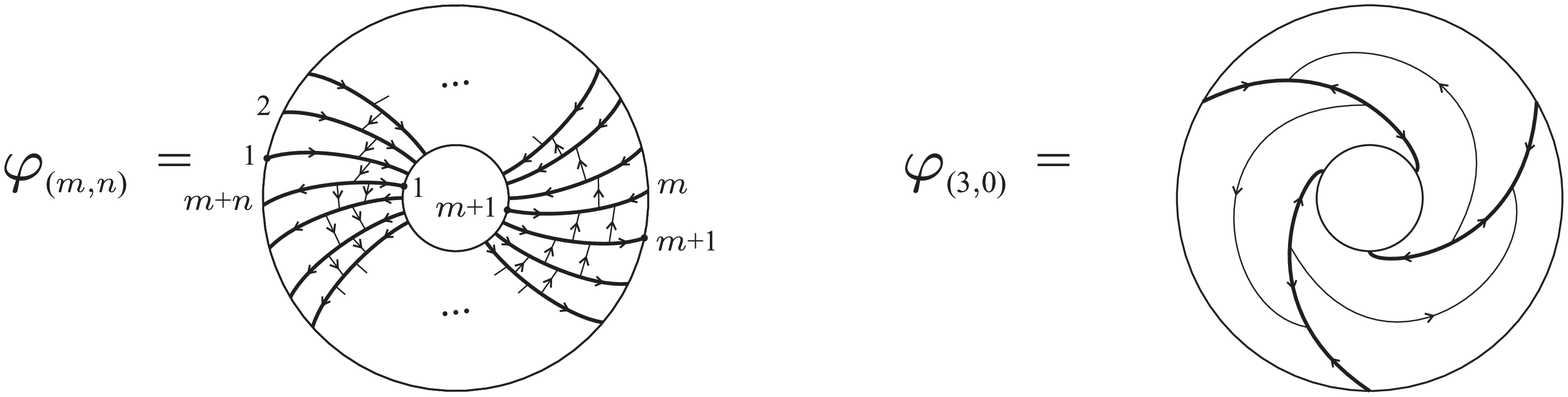}
 \caption{$+^m -^n$-tangle $\varphi_{(m,n)}$ and $+^3$-tangle $\varphi_{(3,0)}$} \label{fig-varphi(m,n)}
\end{center}
\end{figure}

The following inclusion was conjectured in \cite{evans/pugh:2009iv} for the graph $\mathcal{D}^{(6)}$:
$$P^{\mathcal{D}^{(6)}} \supset H^{\beta_{(0,0)}} \oplus H^{0,\overline{0}} \oplus H^{(2,2),\gamma_1,\varepsilon_1} \oplus H^{(2,2),\gamma_2,\varepsilon_2},$$
where $\varepsilon_1, \varepsilon_2 \in \{ \pm 1 \}$, and either $\gamma_1, \gamma_2 \in \mathbb{R}$ or else $\gamma_1 = \overline{\gamma}_2$.
Turning again to the graph $\mathcal{E}^{(8)}$, the following inclusion was conjectured in \cite{evans/pugh:2009iv}:
$$P^{\mathcal{E}^{(8)}} \supset H^{\beta_{(0,0)}} \oplus H^{\beta_{(3,0)}} \oplus H^{\beta_{(0,3)}} \oplus H^{\beta_{(2,2)}} \oplus H^{(3,0),\varepsilon_1} \oplus H^{(0,3),\varepsilon_1} \oplus H^{(2,2),\gamma_1,\varepsilon_2i} \oplus H^{(2,2),\gamma_2,\varepsilon_3i},$$
where $\varepsilon_i \in \{ \pm 1 \}$, $i = 1,2,3$, and either $\gamma_1, \gamma_2 \in \mathbb{R}$ or else $\gamma_1 = \overline{\gamma}_2$.

\section{Almost Calabi-Yau algebras} \label{sect:almostCYalg}

In this section we review the almost Calabi-Yau algebra $A=A(\mathcal{G},W)$ associated to a finite $SU(3)$ $\mathcal{ADE}$ graph $\mathcal{G}$ which carries a cell system $W$, including its matrix Hilbert series of dimensions. These almost Calabi-Yau algebras are essentially the algebras $\Sigma \cong \Pi$ of Section \ref{sect:module_categories}. Then in Section \ref{sect:Nakayama} we describe the Nakayama automorphism for $A$ and the computation of a finite resolution of $A$ as an $A$-$A$ bimodule.

For a finite subgroup $\Gamma \subset SU(2)$, the study of quotient Kleinian singularities $\mathbb{C}^2/\Gamma$ and their resolution has been assisted with the study of the structure of certain noncommutative algebras.
Minimal resolutions of Kleinian singularities can be described via the moduli space of representations of the preprojective algebra associated to the action of $\Gamma$. This leads to general programme to understand singularities via a noncommutative algebra $A$, often called a noncommutative resolution, whose centre corresponds to the coordinate ring of the singularity \cite{vandenBergh:2004}. The algebra should be finitely generated over its centre, and the desired favourable resolution is the moduli space of representations of $A$, whose category of finitely generated modules is derived equivalent to the category of coherent sheaves of the resolution.
In the case of a quotient singularity $\mathbb{C}^3/\Gamma$ for a finite subgroup $\Gamma$ of $SU(3)$, the corresponding noncommutative algebra $A$ is a Calabi-Yau algebra of dimension 3.

Calabi-Yau algebras arise naturally in the study of Calabi-Yau manifolds, providing a noncommutative version of conventional Calabi-Yau geometry.
An algebra $A$ is Calabi-Yau of dimension $n$ if the bounded derived category of the abelian category of finite dimensional $A$-modules is a Calabi-Yau category of dimension $n$. In this case the global dimension of $A$ is $n$ \cite{bocklandt:2008}.
The derived category of coherent sheaves over an $n$-dimensional Calabi-Yau manifold is a Calabi-Yau category of dimension $n$ and they appear naturally in the study of boundary conditions of the $B$-model in superstring theory over the manifold. For more on Calabi-Yau algebras, see e.g. \cite{bocklandt:2008, ginzburg:2006}.

In \cite[Remark 4.5.7]{ginzburg:2006} Ginzburg introduced, in his terminology, $q$-deformed Calabi-Yau algebras. In the case where $q$ is not a root of unity, these algebras are Calabi-Yau algebras of dimension 3.
We study these algebras in the case where $q$ is a root of unity.

Bocklandt \cite{bocklandt:2008} showed that any graded Calabi-Yau algebra of dimension 3 is isomorphic to the path algebra of a directed graph with relations derived from a potential.
We will review the construction of such a path algebra for any finite directed graph $\mathcal{G}$.
Let $[\mathbb{C}\mathcal{G}, \mathbb{C}\mathcal{G}]$ denote the subspace of $\mathbb{C}\mathcal{G}$ spanned by all commutators of the form $xy - yx$, for $x,y \in \mathbb{C}\mathcal{G}$. If $x,y$ are paths in $\mathbb{C}\mathcal{G}$ such that $r(x) = s(y)$ but $r(y) \neq s(x)$, then $xy - yx = xy$, so in the quotient $\mathbb{C}\mathcal{G} / [\mathbb{C}\mathcal{G}, \mathbb{C}\mathcal{G}]$ the path $xy$ will be zero. Then any non-cyclic path, i.e. any path $x$ such that $r(x) \neq s(x)$, will be zero in $\mathbb{C}\mathcal{G} / [\mathbb{C}\mathcal{G}, \mathbb{C}\mathcal{G}]$. If $x = a_1 a_2 \cdots a_p$ is a cyclic path in $\mathbb{C}\mathcal{G}$, then $a_1 a_2 \cdots a_p - a_p a_1 \cdots a_{p-1} = 0$ in $\mathbb{C}\mathcal{G} / [\mathbb{C}\mathcal{G}, \mathbb{C}\mathcal{G}]$, so $x$ is identified with every cyclic permutation of the edges $a_j$, $j=1,\ldots,p$. Then the commutator quotient $\mathbb{C}\mathcal{G} / [\mathbb{C}\mathcal{G}, \mathbb{C}\mathcal{G}]$ may be identified, up to cyclic permutation of the arrows, with the vector space spanned by cyclic paths in $\mathcal{G}$.
One defines a derivation $\partial_a : \mathbb{C}\mathcal{G} / [\mathbb{C}\mathcal{G}, \mathbb{C}\mathcal{G}] \rightarrow \mathbb{C}\mathcal{G}$, and a potential $\Phi \in \mathbb{C}\mathcal{G} / [\mathbb{C}\mathcal{G}, \mathbb{C}\mathcal{G}]$, which is some linear combination of cyclic paths in $\mathcal{G}$. Then the graded algebra
$$A(\mathbb{C}\mathcal{G}, \Phi) = \mathbb{C}\mathcal{G} / \langle \rho_a \rangle,$$
is the quotient of the path algebra by the two-sided ideal generated by the relations $\rho_a = \partial_a \Phi \in \mathbb{C}\mathcal{G}$, for all edges $a$ of $\mathcal{G}$.
When $\Phi$ is homogeneous, the Hilbert series $H_A$ for $A(\mathbb{C}\mathcal{G},\Phi)$ is $H_A(t) = \sum_{p=0}^{\infty} H_{ji}^p t^p$, where the $H_{ji}^p$ are matrices which count the dimension of the subspace $\{ i x j | \; x \in A(\mathbb{C}\mathcal{G},\Phi)_p \}$, where $A(\mathbb{C}\mathcal{G},\Phi)_p$ is the subspace of $A(\mathbb{C}\mathcal{G},\Phi)$ of all paths of length $p$, and $i,j \in A(\mathbb{C}\mathcal{G},\Phi)_0$.

In the case where $\mathcal{G}$ is one of the $ADE$ Dynkin diagrams or affine Dynkin diagrams, an orientation is chosen for each edge of $\mathcal{G}$. The double $\overline{\mathcal{G}}$ of $\mathcal{G}$ has the same vertices as $\mathcal{G}$, but its edges are the union of the edges of $\mathcal{G}$ and its opposite $\mathcal{G}^{\mathrm{op}}$. For any closed loop $\gamma_1 \gamma_2 \cdots \gamma_p$ of length $p$, $p >1$, we define derivatives $\partial_i: \mathbb{C}\overline{\mathcal{G}} / [\mathbb{C}\overline{\mathcal{G}}, \mathbb{C}\overline{\mathcal{G}}] \rightarrow \mathbb{C}\overline{\mathcal{G}}$ for each vertex $i \in \mathfrak{V}_{\mathcal{G}}$ of $\mathcal{G}$ by $\partial_i (\gamma_1 \gamma_2 \cdots \gamma_p) = \sum_j \gamma_j \gamma_{j+1} \cdots \gamma_p \gamma_1 \cdots \gamma_{j-1}$, where the summation is over all $1 \leq j \leq p$ such that $s(\gamma_j) = i$. Then on paths $\gamma \widetilde{\gamma} \in \mathbb{C}\overline{\mathcal{G}} / [\mathbb{C}\overline{\mathcal{G}}, \mathbb{C}\overline{\mathcal{G}}]$, we have
$$\partial_i (\gamma \widetilde{\gamma}) = \left\{ \begin{array}{cl}
                                \gamma \widetilde{\gamma} & \textrm{ if } s(\gamma) = i, \\
                                \widetilde{\gamma} \gamma & \textrm{ if } r(\gamma) = i, \\
                                0 & \textrm{ otherwise.}
                                \end{array} \right.$$
We define a potential $\Phi$ by $\Phi = \sum_{\gamma} \gamma \widetilde{\gamma}$, where the summation is over all edges of $\overline{\mathcal{G}}$. Then $\Pi = \mathbb{C}\overline{\mathcal{G}} / (\partial_i \Phi: i \in \mathfrak{V}_{\mathcal{G}})$ is the preprojective algebra of $\mathcal{G}$ \cite{gelfand/ponomarev:1979}.
For an affine Dynkin diagram $\mathcal{G}$, the preprojective algebra $A(\mathbb{C}\mathcal{G},\Phi)$ is a Koszul algebra, or Calabi-Yau algebra of dimension 2 \cite[Theorem 3.2]{bocklandt:2008}, and it has Hilbert series given by \cite[Theorem 2.3a]{malkin/ostrik/vybornov:2006}
$$H(t) = \frac{1}{1 - \Delta t + t^2}.$$
When $\mathcal{G}$ is one of the Dynkin diagrams or the tadpole graph $T_n$, the preprojective algebra $A(\mathbb{C}\mathcal{G},\Phi)$ is a $(k,2)$-Koszul algebra, or almost Koszul algebra, \cite[Corollary 4.3]{brenner/butler/king:2002}.
There is a `correction' term in the numerator of the Hilbert series, so that \cite[Theorem 2.3b]{malkin/ostrik/vybornov:2006}
$$H(t) = \frac{1 + Pt^h}{1 - \Delta t + t^2},$$
where $P$ is a permutation corresponding to some involution of the vertices of the graph, and $h = k+2$ is the Coxeter number of the graph, where $k$ is the level of $SU(2)$.

Note that there is an insignificant difference in the definition above of the preprojective algebra $A(\mathbb{C}\mathcal{G},\Phi)$, which appears in the algebraic literature, and the appearance of the preprojective algebra $\Pi' \cong \mathrm{EssPath}$ in the theory of subfactors, as mentioned in Section \ref{sect:module_categories}. The Perron-Frobenius weights appear in the definition of $\Pi' = \mathbb{C}\mathcal{G}/\langle \mathrm{Im} \left( F \left( \includegraphics[width=5mm]{fig_nakayama-cup} \right) \right) \rangle$, as in (\ref{eqn:annihilation-lr}), whereas the usual definition in the algebraic literature could be written as $A(\mathbb{C}\mathcal{G},\Phi) = \mathbb{C}\mathcal{G}/\langle \mathrm{Im} \left( F \left( \includegraphics[width=5mm]{fig_nakayama-cup} \right) \right) \rangle$, where the corresponding annihilation operators are given by $c_l(a\widetilde{b}) = \delta_{s(a),s(b)} \, s(a)$, $c_r(\widetilde{b}a) = \delta_{r(a),r(b)} \, r(a)$, without the Perron-Frobenius weights, where $a$ is an edge on $\mathcal{G}$ and $\widetilde{b}$ an edge on $\mathcal{G}^{\mathrm{op}}$.
However, both these definitions give isomorphic preprojective algebras $\Pi \cong \Pi'$, as in \cite[Section 5]{evans/pugh:2010ii}.

We now consider $\mathcal{G}$ an $SU(3)$ $\mathcal{ADE}$ graph or the McKay graph of a finite subgroup $\Gamma \subset SU(3)$. One defines a derivation $\partial_a$ by $\partial_{a} (a_1 \cdots a_n) = \sum_{j} a_{j+1} \cdots a_n a_1 \cdots a_{j-1}$, where the summation is over all indices $j$ such that $a_j = a$. A homogeneous potential $\Phi$ is defined by \cite[Remark 4.5.7]{ginzburg:2006}:
\begin{equation} \label{eqn:potential-Phi}
\Phi = \sum_{a,b,c \in \mathcal{G}_1} W(\triangle^{(a,b,c)}_{s(a),s(b),s(c)}) \cdot \triangle^{(a,b,c)}_{s(a),s(b),s(c)} \quad \in \mathbb{C} \mathcal{G} / [\mathbb{C} \mathcal{G}, \mathbb{C} \mathcal{G}],
\end{equation}
for a cell system $W$, and we denote the algebra $A(\mathbb{C} \mathcal{G}, \Phi)$ by $A(\mathcal{G},W)$. In the case of a finite subgroup $\Gamma \subset SU(3)$ the cell system $W$ is given in \cite[Section 4.4]{ginzburg:2006}.
When $\mathcal{G}$ is the McKay graph of a subgroup of $SU(3)$, $A = A(\mathcal{G},W)$ is a Calabi-Yau algebra of dimension 3 \cite[Theorem 4.4.6]{ginzburg:2006}, and by \cite[Theorem 4.6]{bocklandt:2008} its Hilbert series is given by:
\begin{equation} \label{eqn:H(t)-CYd}
H_A(t) = \frac{1}{1 - \Delta_{\mathcal{G}} t + \Delta_{\mathcal{G}}^T t^2 - t^3}.
\end{equation}

When $\mathcal{G}$ is a finite $SU(3)$ $\mathcal{ADE}$ graph which carries a cell system $W$, and thus yields a braided subfactor, we call $A = A(\mathcal{G},W)$ an almost Calabi-Yau algebra.
Applying the functor $F$ to the morphisms $\mathfrak{U}_p$ in $A_2\mathrm{-}TL$ we obtain a representation of the Hecke algebra on $\mathbb{C}\mathcal{G}$.
If $\langle \rho_a \rangle_p$ denotes the restriction of the ideal $\langle \rho_a \rangle$ in $(\mathbb{C}\mathcal{G})_p$, which is isomorphic to $\sum_{i=1}^{p-1}\mathrm{Im}(F(\mathfrak{U}_i))$, then $A(\mathcal{G},W)_p \cong (\mathbb{C}\mathcal{G})_p / \sum_{i=1}^{p-1}\mathrm{Im}(F(\mathfrak{U}_i)) = \Pi_p \cong \Sigma_p$, where $\Pi$, $\Sigma$ are the graded algebras defined in Section \ref{sect:module_categories}.

The algebra $A = A(\mathcal{G},W)$ has Hilbert series given by \cite[Theorem 3.1]{evans/pugh:2010ii}
\begin{equation} \label{eqn:H(t)-qCYd}
H_A (t) = \frac{1 - P t^h}{1 - \Delta_{\mathcal{G}} t + \Delta_{\mathcal{G}}^T t^2 - t^3},
\end{equation}
where $P$ is a permutation matrix corresponding to a $\mathbb{Z}_3$ symmetry of the graph, and $h = k+3$ is the Coxeter number of $\mathcal{G}$, where $k$ is the level of $SU(3)$.
The permutation matrix $P$ on the vertices of $\mathcal{G}$ is defined by the unique permutation $\nu$ on the graph $\mathcal{G}$ given by the image under $F$ of the generalized Jones-Wenzl projection $f_{(h-3,0)}$, which is described as follows.
If the permutation matrix $P$ in (\ref{eqn:H(t)-qCYd}) is the identity matrix, then the permutation $\nu$ on the graph $\mathcal{G}$ is just the identity. For the other graphs, the permutation $\nu$ is given on the vertices of $\mathcal{G}$ by the permutation matrix $P$, and on $\mathcal{G}_1$ by the unique permutation on the edges of $\mathcal{G}$ such that $s(\nu(a)) = \nu(s(a))$ and $r(\nu(a)) = \nu(r(a))$ (note that there are no double edges on the graphs $\mathcal{G}$ for which $P$ is non-trivial).
The permutation matrix $P$ is the identity for $\mathcal{D}^{(m)}$, $\mathcal{A}^{(m)\ast}$, $m \geq 5$, $\mathcal{E}^{(8)\ast}$, $\mathcal{E}_l^{(12)}$, $l=1,2,4,5$, and $\mathcal{E}^{(24)}$. For the remaining graphs $\mathcal{A}^{(m)}$, $\mathcal{D}^{(m) \ast}$ and $\mathcal{E}^{(8)}$, let $P_0$ be the permutation matrix corresponding to the clockwise rotation of the graph by $2 \pi /3$. Then
$$ P = \left\{
\begin{array}{cl} P_0^2 & \mbox{ for } \quad \mathcal{A}^{(m)}, m \geq 4, \\
                  P_0 & \mbox{ for } \quad \mathcal{E}^{(8)}, \\
                  P_0^{2m} & \mbox{ for } \quad \mathcal{D}^{(m) \ast}, m \geq 5.
\end{array} \right.$$

\subsection{A finite resolution of $A=A(\mathcal{G},W)$ as an $A$-$A$ bimodule} \label{sect:Nakayama}

The algebra $A=A(\mathcal{G},W)$ is a Frobenius algebra, that is, there is a linear function $f:A \rightarrow \mathbb{C}$ such that $(x,y):=f(xy)$ is a non-degenerate bilinear form (this is equivalent to the statement that $A$ is isomorphic to its dual $\widehat{A} = \mathrm{Hom}(A,\mathbb{C})$ as left (or right) $A$-modules). There is an automorphism $\beta$ of $A$, called the Nakayama automorphism of $A$ (associated to $f$), such that $(x,y) = (y,\beta(x))$. Then there is an $A$-$A$ bimodule isomorphism $\widehat{A} \rightarrow {}_1 A_{\beta}$ \cite{yamagata:1996}.
We define a non-degenerate form $A$ by setting $f$ to be the function which is 0 on every element of $A$ of length $< h-3$, and 1 on $u_{i\nu(i)}$ for some $i \in \mathcal{G}_1$, where $\mathcal{G}_1$ denotes the edges of $\mathcal{G}$ and $u_{j\nu(j)}$ denotes a generator of the one-dimensional top-degree space $j \cdot A_{h-3} \cdot \nu(j)$, where $\nu$ is the permutation of the vertices of $\mathcal{G}$ given by the permutation matrix $P$ in (\ref{eqn:H(t)-qCYd}). Then using the relation $(x,y) = (y,\beta(x))$ this determines the value of $f$ on $u_{j\nu(j)}$, for all other $j \in \mathcal{G}_1$. We normalize the $u_{j\nu(j)}$ such that $f(u_{j\nu(j)}) = 1$ for all $j \in \mathcal{G}_1$.
Then following similar methods to \cite{brenner/butler/king:2002} we showed that the Nakayama automorphism $\beta$ of $A$ is defined on $\mathcal{G}$ by $\beta = \nu$ \cite[Theorem 4.6]{evans/pugh:2010ii}.

It was shown in \cite[Theorem 5.1]{evans/pugh:2010ii} that $A$ has the following finite resolution as an $A$-$A$ bimodule:
\begin{equation} \label{seq:q-CY3_resolution}
0 \rightarrow {}_1 A_{\beta^{-1}} \stackrel{\mu_4}{\longrightarrow} A \otimes_R A \stackrel{\mu_3}{\longrightarrow} A \otimes_R \widehat{V} \otimes_R A \stackrel{\mu_2}{\longrightarrow} A \otimes_R V \otimes_R A \stackrel{\mu_1}{\longrightarrow} A \otimes_R A \stackrel{\mu_0}{\longrightarrow} A \rightarrow 0
\end{equation}
where $R$, $V$ are the $A$-$A$ bimodules generated by the vertices, edges of $\mathcal{G}$ respectively, the $A$-$A$ bimodule $\widehat{V}$ is the dual space of $V$,
with the connecting $A$-$A$ bimodule maps given by
\begin{eqnarray*}
\mu_0(1 \otimes 1) & = & 1, \\
\mu_1(1 \otimes a \otimes 1) & = & a \otimes 1 - 1 \otimes a, \\
\mu_2(1 \otimes \widetilde{a} \otimes 1) & = & \sum_{b,b' \in \mathcal{G}_1} W_{abb'} (b \otimes b' \otimes 1 + 1 \otimes b \otimes b'), \\
\mu_3(1 \otimes 1) & = & \sum_{a \in \mathcal{G}_1} a \otimes \widetilde{a} \otimes 1 - \sum_{a \in \mathcal{G}_1} 1 \otimes \widetilde{a} \otimes a, \\
\mu_4(1) & = & \sum_j w_j \otimes w_j^{\ast},
\end{eqnarray*}
where $a \in \mathcal{G}_1$, the $w_j$ are a homogeneous basis for $A$, $w_j^{\ast}$ is its dual basis under the non-degenerate form on $A$, and we denote by $W_{abb'}$ the cell $W(\triangle^{(a,b,b')}_{s(a),s(b),s(b')})$.
This \emph{almost Calabi-Yau} condition should be compared with the Calabi-Yau condition \cite[Theorem 4.3]{bocklandt:2008}, which says that $A$ is a Calabi-Yau algebra of dimension 3 if an only if it has the following resolution as an $A$-$A$ bimodule:
$$0 \rightarrow A \otimes_R A \rightarrow A \otimes_R \widehat{V} \otimes_R A \rightarrow A \otimes_R V \otimes_R A \rightarrow A \otimes_R A \rightarrow A \rightarrow 0.$$
In particular, this is the case for $A = A(\mathcal{G}_{\Gamma},W)$, for the McKay graph of a finite subgroup $\Gamma \subset SU(3)$ \cite[Theorem 4.4.6]{ginzburg:2006}.

The resolution (\ref{seq:q-CY3_resolution}) of the almost Calabi-Yau algebra $A = A(\mathcal{G},W)$ will yield a projective resolution of $A$ as an $A$-$A$ bimodule.
The objective in deriving this resolution is to provide a basis for the computation of the Hochschild (co)homology and cyclic homology of the algebras $A(\mathcal{G},W)$ for the $SU(3)$ $\mathcal{ADE}$ graphs. Beginning with a pair $(\mathcal{G},W)$ given by a cell system $W$ on an $SU(3)$ $\mathcal{ADE}$ graph $\mathcal{G}$, we construct a subfactor $N \subset M$ which yields a nimrep which recovers the graph $\mathcal{G}$ as described in Section \ref{sect:braided_subfactors}. Then we can construct the algebra $A(\mathcal{G},W)$ whose Hochschild (co)homology and cyclic homology only depends on the original pair $(\mathcal{G},W)$, or equivalently, on the subfactor $N \subset M$. Thus the Hochschild (co)homology and cyclic homology of $A$ should be regarded as invariants for the subfactor $N \subset M$.
This is work in progress.

\paragraph{Acknowledgements}

Both authors were supported by the Marie Curie Research Training Network MRTN-CT-2006-031962 EU-NCG, and the second author was supported by a scholarship from the School of Mathematics, Cardiff University.

\end{document}